\documentclass[a4paper,11pt,reqno]{amsart}
\usepackage[greek,english]{babel}
\usepackage[utf8]{inputenc}
\usepackage[dvipsnames]{xcolor}
\definecolor{darkred}{RGB}{192,0,0}
\definecolor{darkblue}{RGB}{0,0,192}
\definecolor{MyGreen}{RGB}{0,112,0}
	\definecolor{amber}{rgb}{1.0, 0.75, 0.0}
	\definecolor{darkorange}{rgb}{1.0, 0.27, 0.0}
\usepackage[T1]{fontenc}
\usepackage{amsmath}

\usepackage[title,titletoc]{appendix}

\usepackage{newtxtext}
\usepackage{fancyhdr,amsthm,mathtools,stmaryrd,amscd,tikz-cd,amssymb}
\usepackage{caption}
\usepackage[left=0.95in, right=0.95in, top=0.95in, bottom=0.95in, includefoot, headheight=13.6pt]{geometry}
\usepackage{csquotes}
\usepackage{amssymb, graphicx}
\usepackage{pmboxdraw}
\newcommand{\widesim}[2][2]{
  \mathrel{\overset{#2}{\scalebox{#1}[1]{$\sim$}}}
}
\usepackage{bbm}
\usepackage{dsfont}
\usepackage{nicematrix}
\usepackage{tikz,pgfplots}
\usepackage{tikz-cd}
\usetikzlibrary{cd}
\usetikzlibrary{arrows}
\usetikzlibrary{matrix}
\usetikzlibrary{positioning}
\usetikzlibrary{intersections,calc}
\usepackage{tikz-3dplot}
\tikzcdset{arrow style=math font}
\pgfplotsset{compat=1.17}
\usepackage{pgfplots}
\usepackage{comment}
\usepackage{numprint}
\usepackage{url}
\usepackage{subcaption}

\usepackage[export]{adjustbox}
\usepackage{hyperref}
\hypersetup{
    colorlinks=true,
    linkcolor=orange!80!black,
    filecolor=magenta,      
    urlcolor=RoyalBlue,
    citecolor=RoyalBlue,
}
\usepackage{mathrsfs}  
\usepackage{enumitem}
\usepackage{array}
\usepackage{braket}
\usepackage{physics}

\allowdisplaybreaks

\usepackage[subscriptcorrection]{newtxmath}

\DeclareMathAlphabet{\mathbb}{U}{msb}{m}{n}
\DeclareMathAlphabet{\mathfrak}{U}{euf}{m}{n}

\usepackage[font=small,labelfont=bf,margin=1.5cm]{caption}

\usepackage[nameinlink]{cleveref}

\numberwithin{equation}{section}
\theoremstyle{plain}
\newtheorem{theorem}[equation]{Theorem}

\newtheorem*{theorem*}{Theorem}

\newtheorem{lemma}[equation]{Lemma}
\newtheorem*{lemma*}{Lemma}
\newtheorem{corollary}[equation]{Corollary}
\newtheorem{proposition}[equation]{Proposition}

\newtheorem{thm/conj}[equation]{Theorem/Conjecture}
\newtheorem{conjecture}[equation]{Conjecture}

\theoremstyle{definition}
\newtheorem{examplex}[equation]{Example}
\newenvironment{example}
{\pushQED{\qed}\examplex}
{\popQED\endexamplex}
\newtheorem{definition}[equation]{Definition}
\newtheorem{remark}[equation]{Remark}

\counterwithin{figure}{section}

\theoremstyle{remark}

\newcommand{\bfa}{\mathbf{a}}
\newcommand{\bfb}{\mathbf{b}}
\newcommand{\bfc}{\mathbf{c}}

\newcommand{\bfm}{\mathbf{m}}

\newcommand{\bfu}{\mathbf{u}}
\newcommand{\bfv}{\mathbf{v}}
\newcommand{\bfw}{\mathbf{w}}

\newcommand{\bbC}{\mathbb{C}}

\newcommand{\bbK}{\mathbb{K}}

\newcommand{\bbP}{\mathbb{P}}

\newcommand{\bbR}{\mathbb{R}}

\newcommand{\bbZ}{\mathbb{Z}}

\newcommand{\calA}{\mathcal{A}}
\newcommand{\calB}{\mathcal{B}}
\newcommand{\calC}{\mathcal{C}}

\newcommand{\calE}{\mathcal{E}}
\newcommand{\calF}{\mathcal{F}}

\newcommand{\calL}{\mathcal{L}}
\newcommand{\calM}{\mathcal{M}}

\newcommand{\calO}{\mathcal{O}}
\newcommand{\calP}{\mathcal{P}}
\newcommand{\calQ}{\mathcal{Q}}
\newcommand{\calR}{\mathcal{R}}
\newcommand{\calS}{\mathcal{S}}
\newcommand{\calT}{\mathcal{T}}
\newcommand{\calU}{\mathcal{U}}
\newcommand{\calV}{\mathcal{V}}
\newcommand{\calW}{\mathcal{W}}

\newcommand{\rmr}{\mathrm{r}}

\newcommand{\eps}{\varepsilon}

\newcommand{\C}{\mathbb{C}}

\DeclareMathOperator{\coim}{coim}
\DeclareMathOperator{\id}{id}
\DeclareMathOperator{\im}{Im}
\DeclareMathOperator{\Ker}{Ker}
\DeclareMathOperator{\Mat}{Mat}
\DeclareMathOperator{\rk}{\mathbf{R}}
\DeclareMathOperator{\Seg}{Seg}

\usepackage[textwidth=0.8in]{todonotes}
\Crefname{examples}{Example}{Examples}
\Crefname{rmks}{Remark}{Remarks}

\usepackage{titlesec}
\titleformat{name=\section}{
 \vspace{10pt}\scshape\fontfamily{ptm}\raggedright\LARGE}{}{0em}{\hspace{-0.4pt}\LARGE \thesection.\hspace{0.5em}}[\color{black}\titlerule \vspace{5pt}]

\titleformat{name=\section,numberless}
{\vspace{10pt}\scshape\fontfamily{ptm}\raggedright\LARGE}{}{0em}
  {\hspace{-0.4pt}\LARGE}[\color{black}\titlerule \vspace{5pt}]

\titleformat{name=\subsection}{\vspace{7pt}}{\bfseries\Large\thetitle.}{0.5em}{\Large\bfseries}

\titleformat{name=\subsection*,numberless}{\vspace{7pt}}{\bfseries\Large\thetitle.}{0.5em}{\Large\bfseries}

\renewcommand{\sectionmark}[1]{\markleft{\normalfont \scshape\fontfamily{ptm}\selectfont V.~Borovik, C.~Flavi, P.~Pielasa, A.~Shatsila, and J.~Song}}\renewcommand{\subsectionmark}[1]{\markright{\normalfont \scshape\fontfamily{ptm}\selectfont On the construction of a counterexample to Strassen's conjecture}}
  \keywords{Strassen’s conjecture, tensor decomposition, tensor rank, augmented tensor}

\subjclass{Primary 15A69, 15A72; Secondary 14N07, 14N15}

\allowdisplaybreaks

\title{On the construction of a counterexample to \\ Strassen's rank additivity conjecture}
\author{Viktoriia Borovik}
\author{Cosimo Flavi}
\author{Pawe{\l} Pielasa}
\author{Anatoli Shatsila}
\author{Jeyoung Song}
\address{{\normalfont (Viktoriia Borovik)}
\normalfont \scshape\fontfamily{ptm}\selectfont Max Planck Institute for Mathematics in the Sciences,
 \normalfont{Inselstrasse 22,
04103 Leipzig, Germany.}}
\email{borovik@mis.mpg.de}
\address{{\normalfont (Cosimo Flavi)}
\normalfont \scshape\fontfamily{ptm}\selectfont Wydział Matematyki, Informatyki i Mechaniki, Uniwersytet Warszawski,
 \normalfont{ul.~Stefana Banacha
2, 02-097 Warsaw, Poland.}}
\email{c.flavi@uw.edu.pl}
\address{{\normalfont (Pawe{\l} Pielasa)}
\normalfont \scshape\fontfamily{ptm}\selectfont Wydział Matematyki, Informatyki i Mechaniki, Uniwersytet Warszawski,
 \normalfont{ul.~Stefana Banacha
2, 02-097 Warsaw, Poland.}}
\email{p.pielasa@student.uw.edu.pl}
\address{{\normalfont (Anatoli Shatsila)}
\normalfont \scshape\fontfamily{ptm}\selectfont Wydział Matematyki i Informatyki, Uniwersytet Jagielloński,
 \normalfont{ul.~Profesora Stanisława Łojasiewicza 6, 30-348 Kraków, Poland.}}
\email{anatoli.shatsila@doctoral.uj.edu.pl}
\address{{\normalfont (Jeyoung Song)}
\normalfont \scshape\fontfamily{ptm}\selectfont Department of Mathematics and Statistical Science, University of Idaho,
 \normalfont{875 Perimeter Drive, MS 1103, Moscow, Idaho 83844-1103, United States of America.}}
\email{jeyoungsong@uidaho.edu}
\date{}

\usepackage[dotinlabels]{titletoc}

\dottedcontents{section}[1.1em]{\hypersetup{linkcolor=black}}{1.2em}{5pt}

\begin{document}
\renewcommand{\abstractname}{\normalfont \scshape\fontfamily{ptm}\selectfont{Abstract}}
\renewcommand{\contentsname}{\normalfont \scshape\fontfamily{ptm}\selectfont{Contents}}
\begin{abstract}
The rank additivity conjecture, first formulated by Volker Strassen in 1973, states that the rank of the direct sum of two independent tensors is equal to the sum of their individual ranks. In the last decades, this conjecture has been a central topic in tensor rank theory and its implications for computational complexity.
In 2019, Yaroslav Shitov disproved this conjecture in its general form 
by showing the existence of a counter-example using a dimension counting argument.
In this paper, we provide an overview of the Strassen problem and Shitov's work and revisit his counterexample with a~detailed explanation, offering an alternative proof. 
\end{abstract}
\maketitle
\thispagestyle{empty}
\tableofcontents

\section{Introduction}
One of the central questions of complexity theory is: given two computational tasks, can they be performed simultaneously using fewer resources than if they were performed separately? A particularly important case of this question arises in the context of \emph{matrix multiplication}, see \cite{CHILO18, CU03, CW87, LeGall14, LM18, LO15, Pan17, Seyn18, Wil12}. Suppose we are given four matrices over some field $\mathbb{K}$: $M_1 \in \mathrm{Mat}_{i_1 \times j_1}(\mathbb{K})$, $N_1 \in \mathrm{Mat}_{j_1 \times k_1}(\mathbb{K})$, and $M_2 \in \mathrm{Mat}_{i_2 \times j_2}(\mathbb{K})$, $N_2 \in \mathrm{Mat}_{j_2 \times k_2}(\mathbb{K})$. If we want to simultaneously compute the products $M_1N_1$ and $M_2N_2$, what is the minimum number of scalar multiplications required? Is it simply the sum of the complexities of computing these products separately? For an overview of complexity theory problems and their connections to algebraic geometry and representation theory we refer the reader to the books \cite{BCS97, Lan17}.

This leads naturally to the notion of \emph{tensor rank}, since we are essentially asking about the tensor rank of the \emph{matrix multiplication tensor}, which is a bilinear map
$$M_{\langle i,j,k\rangle}\colon \mathrm{Mat}_{i\times j}(\mathbb{K}) \times \mathrm{Mat}_{j\times k}(\mathbb{K}) \to \mathrm{Mat}_{i\times k}(\mathbb{K}).$$ The space $\mathrm{Mat}_{s\times r}(\mathbb{K})$ is isomorphic to $\mathbb{K}^{s r}$ as a vector space, so we can reinterpret $M_{\langle i,j,k\rangle}$ as a tensor:
\[ M_{\langle i,j,k\rangle} \in (\mathbb{K}^{i\times j})^{*} \otimes (\mathbb{K}^{j\times k})^{*} \otimes \mathbb{K}^{i\times k}. \]
Let $A$, $B$, and $C$ be finite-dimensional vector spaces over $\mathbb{K}$. A tensor $T \in A \otimes B \otimes C$ is called \emph{simple} if it can be written as $a \otimes b \otimes c$ for some elements $a \in A$, $b \in B$, and $c \in C$. The \emph{rank} $\mathbf{R}(T)$ of a tensor $T\in A\otimes B\otimes C$ is the minimum number of simple tensors needed to express $T$ as a linear combination of those tensors, that is,
    \[
    \rk(T)\coloneqq\min\Set{r\in\mathbb{N}|T=\sum_{i=1}^r a_i\otimes b_i\otimes c_i:a_i\in A,\, b_i\in B,\, c_i\in C,\, \forall i=1,\dots,r}.
    \]
In general, the rank serves as a measure of complexity: the higher the rank, the more complex the tensor.
In particular, the minimum number of scalar multiplications required to compute $M_1N_1$ or $M_2N_2$ corresponds to $\mathbf{R}(M_{\langle i_1,j_1,k_1\rangle})$ and $\mathbf{R}(M_{\langle i_2,j_2,k_2\rangle})$, respectively.
With these definitions, the question of simultaneous matrix multiplication becomes: 
\begin{center}
    \emph{Is the minimum number of scalar multiplications required to simultaneously compute $M_1N_1$ and $M_2N_2$ equal to $\mathbf{R}(M_{\langle i_1,j_1,k_1\rangle}) + \mathbf{R}(M_{\langle i_2,j_2,k_2\rangle})$?}
\end{center}
 More generally, the above question applies to arbitrary tensors. The positive answer to this question is known in the literature as \emph{Strassen's direct sum conjecture}.
\begin{conjecture}[Strassen's direct sum conjecture]\label{conj: Strassen}
Given 
$T_1 \in A_1 \otimes B_1 \otimes C_1$ and $T_2 \in A_2 \otimes B_2 \otimes C_2$, let 
$T = T_1 \oplus T_2$ be their direct sum in $(A_1 \oplus A_2) \otimes (B_1 \oplus B_2) \otimes (C_1 \oplus C_2)$. Then
    $$\mathbf{R}(T) = \mathbf{R}(T_1) + \mathbf{R}(T_2).$$
\end{conjecture}
Volker Strassen originally posed this question in \cite{Strassen1973}, asking whether the multiplicative complexity of the union of two bilinear systems depending on disjoint sets of variables is equal to the sum of the complexities of the individual systems, with special attention to matrix multiplication. This question was generalized in many directions, for example, to higher-order tensors, symmetric or skew-symmetric tensors. Moreover, one can state an analogous additivity conjecture for the \emph{border rank} of tensors. 

Additivity results for some families of bilinear forms were shown in \cite{FEIGWinograd, validityConj}; see also \cite[Section 4]{LM17}. They were generalized in \cite{10.1137/19M1243099, Rup24} for small 3-order tensors.
However, for an analog of the~\Cref{conj: Strassen} for the border rank, the answer is generally negative. Sch\"onhage provided a family of counterexamples in \cite{borderrankconj}. In \cite{10.1137/19M1243099} it was shown that if the dimensions of the spaces $A_1 \oplus A_2$, $B_1 \oplus B_2$ and $C_1 \oplus C_2$ are at most 4, then there is no counterexample for the additivity of the border rank. 
%In contrast, no counterexample at all is known for the symmetric version of \Cref{conj: Strassen}. Some positive results in this setting appeared in \cite{CARLINI20153149, teitler, carlini17, CARLINI2017630, CMM18}. 

Despite small positive examples for Strassen's question, it is now known that the conjecture is false in the general case. This was expected on the basis of Sch\"onhage's counterexamples for the border rank version, although finding an explicit counterexample for the tensor rank version was considered difficult due to the high dimensionality of the desired tensor spaces, see \cite[Theorem 4.1]{LM17}. \Cref{conj: Strassen} remained open for a long time until Yaroslav Shitov constructed a counterexample in 2019. 
\begin{theorem}[{\cite[Theorem 4]{shitov}}] \label{maintheorem}
    There exist vector spaces $\calA_1$, $\calB_1$, $\calC_1$, $\calA_2$, $\calB_2$, $\calC_2$, and tensors $\mathcal{T}_1 \in \mathcal{A}_1 \otimes \mathcal{B}_1 \otimes \mathcal{C}_1$ and $\mathcal{T}_2 \in \mathcal{A}_2 \otimes \mathcal{B}_2 \otimes \mathcal{C}_2$ such that  $$\mathbf{R}(\mathcal{T}_1 \oplus \mathcal{T}_2) < \mathbf{R}(\mathcal{T}_1) + \mathbf{R}(\mathcal{T}_2).$$
\end{theorem}
Recently, Shitov constructed a counterexample to the symmetric version of \Cref{conj: Strassen} (see \cite[Theorem 3.9]{Shi24}). More precisely, denoting by $\mathbf{SR}(\calT)$ the symmetric rank of a tensor $\calT$, he proved the existence of two symmetric tensors $\mathcal{T}_1, \mathcal{T}_2$ such that $$\mathbf{SR}(\mathcal{T}_1 \oplus \mathcal{T}_2) < \mathbf{R}(\mathcal{T}_1) + \mathbf{R}(\mathcal{T}_2) \leq \mathbf{SR}(\mathcal{T}_1) + \mathbf{SR}(\mathcal{T}_2).$$ %where $\mathbf{SR}$ denotes the symmetric rank. 
Before this work, some partial positive results in this setting appeared in  \cite{CARLINI20153149, teitler, carlini17, CARLINI2017630, CMM18}.

The purpose of this paper is to present a highly elaborated and rigorous version of Shitov's proof of \Cref{maintheorem}. 
The explicit proof is presented in \Cref{Sec3} and rests on two main results, appearing in \cite[Claims 5 and 6]{shitov}. The first of these results is stated in \Cref{lemma1}.
We have retained the structure and ideas of its original proof presented by Shitov in \cite[Claim 5]{shitov}, providing a more geometric point of view by using the general language of tensors, with additional refinements and clarifications. The proof, entirely presented in \Cref{Sec4}, is quite technical and requires some tools connected to tensor decomposition theory, such as augmented tensors, clones, and modifications of tensors, which are introduced in \Cref{Sec2}, including relevant definitions, background material, and illustrative examples and figures to develop intuition. 

The strategy used to prove \Cref{lemma1} is based on a construction via induction. The first step only involves considering a one-dimensional vector space. 
This is the base case of induction, which  corresponds to the simplest version of \Cref{lemma1}, stated in \Cref{prop:induction_case1}. Iterating the process of \Cref{prop:induction_case1} 
several times, it is possible to prove the inductive step, which completes the proof.

In \Cref{Sec5} we provide the original part of our work. We present an alternative version of \cite[Claim 6]{shitov}, stated in \Cref{lemma2}.  
Our strategy is based on methods from basic differential geometry and linear algebra, which allow us to guarantee the existence of two specific tensors that can be used together with \Cref{lemma1} to obtain the proof of \Cref{maintheorem}.
Our argument is inspired by a version presented by Joseph~M.~Landsberg in \cite[Section 2.5]{Lan19}, jointly developed
with Mateusz Michałek and Tim Seynnaeve. In that work, the authors provide a new version of \cite[Claim 5]{shitov}. However, an inaccuracy emerged in the analysis of the final part of the procedure  (see \Cref{rem:inaccuracy_Lan}), which we clarify in this new version of the proof.

 To improve the readability of the main arguments, we collect proofs of technical lemmas and propositions used in the main text in \hyperref[SecA]{Appendix~A}. They are mostly computational in nature and are included separately.

\markboth{\normalfont \scshape\fontfamily{ptm}\selectfont V.~Borovik, C.~Flavi, P.~Pielasa, A.~Shatsila, and J.~Song}{\normalfont \scshape\fontfamily{ptm}\selectfont On the construction of a counterexample to Strassen's conjecture}
The results contained in \cite{shitov} are valid for tensors over an arbitrary infinite field $\bbK$. For our purposes, we restrict ourselves to the field of complex numbers, that is, $\mathbb{K} = \C$.
Throughout the paper, we denote any finite-dimensional vector space by a capital letter, while its corresponding lowercase letter in bold denotes its dimension. For instance,
$A,B,C$ are vector spaces, and their dimensions are denoted as
\[\bfa\coloneqq\dim A,\quad \bfb\coloneqq\dim B,\quad\bfc\coloneqq\dim C.\]
In particular, we may often assume that $A=\C^\bfa$, $B=\C^\bfb$, $C=\C^\bfc$. Other notations include:
\begin{enumerate}[label=(\arabic*), left= 0pt, widest=*,nosep]
    \item We suppose that every tensor denoted by~$T$ belongs to the tensor product $A\otimes B\otimes C$, unless differently stated. 
    \item For any linear map $f\colon V\to W$, the space $\coim f\coloneqq V/\ker f$ denotes the coimage of $f$.
    \item For any $n\in\mathbb{Z}_{>0}$, we denote the elements of the canonical basis of $\bbC^n$ by $e_1,\dots,e_n$. In \Cref{Sec4} we often denote a (canonical) basis of a vector space by $\Sigma$.
    \item $T_1 \boxtimes T_2$ denotes the Kronecker product of the two tensors $T_1$ and $T_2$.
\end{enumerate}

\subsection*{Acknowledgements}
The first, third, fourth, and fifth authors would like to express their gratitude to Mateusz Michałek for inspiring them
to write this paper
during their
stay at the University of Konstanz, as well as for his helpful discussions and meaningful remarks. 
The second author would like to thank Giorgio Ottaviani for introducing him to the topic, and for the helpful discussions related to it.
The authors would like to thank Fulvio Gesmundo, Joachim Jelisiejew, Joseph M.~Landsberg, Mateusz Micha{\l}ek, Pierpaola Santarsiero, and Yaroslav Shitov for very helpful comments.
Part of this paper was written while the second author was a research fellow at the {University of Florence}.
The second author is a member of the research group \textit{Gruppo Nazionale per le Strutture Algebriche Geometriche e Affini} (GNSAGA) of \textit{Istituto Nazionale di Alta Matematica} (INdAM) and has been supported by the scientific project \textit{Multilinear Algebraic Geometry} of the program \textit{Progetti di ricerca di Rilevante Interesse Nazionale} (PRIN), Grant Assignment Decree No.~973, adopted on 06/30/2023 by the Italian Ministry of University and Research (MUR), and by the project \textit{Thematic Research Programmes}, Action I.1.5 of the program \textit{Excellence Initiative -- Research University} (IDUB) of the Polish Ministry of Science and Higher Education.

\section{Preliminaries}\label{Sec2}
In this section, we provide some notions related to tensors that we need to analyze the counterexample provided by Yaroslav~Shitov. We mainly consider tensors of order $3$. 

By the definition of rank, a non-zero tensor $T$ is simple if and only if we have $\rk(T)= 1$.
The set of projective points corresponding to these tensors is the Segre variety, which is the image of the map
\[
\Seg\colon \bbP A\times \bbP B\times \bbP C\to \bbP(A\otimes B \otimes C),
\]
which is defined by
\[
\Seg([a],[b],[c])=[a\otimes b\otimes c],
\]
for every non-zero $a\in A$, $b\in B$, and $c\in C$. Denoting the Segre variety by  $\Sigma_3\coloneqq \Seg(\bbP A\times \bbP B\times \bbP C)$, we then have
\[
\Sigma_3=\Set{[T]\in \bbP(A\otimes B\otimes C)|\rk (T)= 1}.
\]
It is also possible to generalize this definition to tensors of higher rank.  In particular, for any $r\in\mathbb{Z}_{>0}$, the \textit{$r$-th secant variety} of the Segre variety $\Sigma_3$ is the algebraic closure of the set of projective points associated to tensors having rank at most $r$, that is,
\[
\sigma_r(\Sigma_3)\coloneqq\overline{\Set{[T]\in \bbP(A\otimes B\otimes C)|\rk (T)\leq r}}.
\]
Equivalently, $\sigma_r(\Sigma_3)$ is the set of tensors of \emph{border rank} at most $r$, see \cite{Lan12}.
The notion of border rank, which was introduced for the first time in \cite{BCRL79} and then named in this way in \cite{BLR80}, is then strictly connected to secant varieties. In general, the notion of secant varieties can be extended to any projective variety. In particular, one of the main related problems is determining the dimension of the secant varieties, which is, in general, very hard.
For further details, there are several texts in the literature that can be consulted (e.g.,~\cite{Lan12,BCC+18}).

An important classical result that is largely used for getting the dimensions of secant varieties is Terracini's lemma, provided by Alessandro Terracini in 1911 (see \cite{Ter11}), concerning the dimensions of tangent spaces of secant varieties at general points. We recall that by a \textit{general} element in a subset of a vector space we refer to all the elements contained in a Zariski open subset. Since any Zariski open set is dense, this term is used to say that a property holds for \textit{almost any element} of the considered set. 

The relevant fact is that for any $r\in\mathbb{Z}_{>0}$, the general tensor in $\sigma_r(\Sigma_3)$ has rank $r$. In particular, the \textit{generic rank} for tensors in the space $A\otimes B\otimes C$, that is, the rank of a general tensor in $A\otimes B\otimes C$, is the minimum number $r\in\mathbb{Z}_{>0}$ such that the $r$-th secant variety of $\Sigma_3$ fills the whole space, that is, the value
\[
{\rmr}_{\mathrm{gen}}\coloneqq \min\Set{r\in\mathbb{Z}_{>0}|\sigma_r(\Sigma_3)=\bbP(A\otimes B\otimes C)}.
\]

For the case where $A=B=C=\bbC^m$, with $m\geq 4$, the dimensions of the secant varieties of the Segre variety $\Sigma_{3,m}\coloneqq \Seg\bigl((\bbP\bbC^m)^3\bigr)$ have been completely determined. In particular, in \cite[Theorem 4.4]{Lic85}, it is proved that
\begin{equation}\label{formula:secant_variety_Segre}
\dim\bigl(\sigma_r(\bbP\Sigma_{3,m})\bigr)=\min\{r(3m-2)-1,m^3-1\},
\end{equation}
where $m^3-1=\dim \mathbb{P}\bigl((\bbC^m)^{\otimes 3}\bigr)$. Therefore, the value of the generic rank follows immediately, as shown in \cite[Corollary 4.5]{Lic85},
and it is equal to the value
\begin{equation}\label{formula:generic_rank_cubes}
\rmr_{\mathrm{gen}}^m=\biggl\lceil\dfrac{m^3}{3m-2}\biggr\rceil.
\end{equation}

Despite  the fact that the generic rank for tensors in $(\bbC^m)^{\otimes 3}$ has been computed, determining the rank of a given tensor remains a hard issue. However, several methods to deal with this problem have been developed over the years. One of these is the substitution method. Before stating this, for our purposes, we  present in the following subsections the notions of flattenings, modifications of tensors, and clones of tensors.

\subsection{Flattenings and modifications of tensors}
\begin{figure}[hbt!]
\centering
\tdplotsetmaincoords{68}{110}
\begin{tikzpicture}[tdplot_main_coords,line join=round,scale=0.65]
\begin{scope}
\foreach \y in {0,1,...,3}{
\foreach \x in {0,1,...,3}{
\begin{scope}[shift={(-3+\y,\x,-2)}]
\pgfsetfillopacity{0.3}
\pgfsetstrokeopacity{0.4}
\pgfsetlinewidth{0.5pt}

\filldraw[fill=gray!70,draw=gray] (0,0,0) -- (1,0,0) -- (1,0,1) -- (0,0,1) -- (0,0,0);
\filldraw[fill=gray!70,draw=gray] (0,0,0) -- (1,0,0) -- (1,1,0) -- (0,1,0) -- (0,0,0);
\filldraw[fill=gray!70,draw=gray] (0,0,0) -- (0,1,0) -- (0,1,1) -- (0,0,1) -- (0,0,0);
\filldraw[fill=gray!50,draw=gray] (1,0,0) -- (1,1,0) -- (1,1,1) -- (1,0,1) -- (1,0,0);
\filldraw[fill=gray!60,draw=gray] (0,0,1) -- (1,0,1) -- (1,1,1) -- (0,1,1) -- (0,0,1);
\filldraw[fill=gray!70,draw=gray] (0,1,0) -- (1,1,0) -- (1,1,1) -- (0,1,1) -- (0,1,0);

\end{scope}}}

\foreach \y in {0,1,...,3}{
\foreach \x in {0,1,...,3}{
\begin{scope}[shift={(-3+\y,\x,-1)}]
\pgfsetfillopacity{0.3}
\pgfsetstrokeopacity{0.4}
\pgfsetlinewidth{0.5pt}

\filldraw[fill=gray!70,draw=gray] (0,0,0) -- (1,0,0) -- (1,0,1) -- (0,0,1) -- (0,0,0);
\filldraw[fill=gray!70,draw=gray] (0,0,0) -- (1,0,0) -- (1,1,0) -- (0,1,0) -- (0,0,0);
\filldraw[fill=gray!70,draw=gray] (0,0,0) -- (0,1,0) -- (0,1,1) -- (0,0,1) -- (0,0,0);
\filldraw[fill=gray!50,draw=gray] (1,0,0) -- (1,1,0) -- (1,1,1) -- (1,0,1) -- (1,0,0);
\filldraw[fill=gray!60,draw=gray] (0,0,1) -- (1,0,1) -- (1,1,1) -- (0,1,1) -- (0,0,1);
\filldraw[fill=gray!70,draw=gray] (0,1,0) -- (1,1,0) -- (1,1,1) -- (0,1,1) -- (0,1,0);

\end{scope}}}

\foreach \x in {0,1,...,3}{
\begin{scope}[shift={(-3,\x,0)}]
\pgfsetfillopacity{0.3}
\pgfsetstrokeopacity{0.4}
\pgfsetlinewidth{0.5pt}

\filldraw[fill=gray!70,draw=gray] (0,0,0) -- (1,0,0) -- (1,0,1) -- (0,0,1) -- (0,0,0);
\filldraw[fill=gray!70,draw=gray] (0,0,0) -- (1,0,0) -- (1,1,0) -- (0,1,0) -- (0,0,0);
\filldraw[fill=gray!70,draw=gray] (0,0,0) -- (0,1,0) -- (0,1,1) -- (0,0,1) -- (0,0,0);
\filldraw[fill=gray!50,draw=gray] (1,0,0) -- (1,1,0) -- (1,1,1) -- (1,0,1) -- (1,0,0);
\filldraw[fill=gray!60,draw=gray] (0,0,1) -- (1,0,1) -- (1,1,1) -- (0,1,1) -- (0,0,1);
\filldraw[fill=gray!70,draw=gray] (0,1,0) -- (1,1,0) -- (1,1,1) -- (0,1,1) -- (0,1,0);

\end{scope}}

\begin{scope}[shift={(-2,0,0)}]
\pgfsetfillopacity{0.8}
\pgfsetlinewidth{0.5pt}

\filldraw[fill=MyGreen!90,draw=MyGreen] (0,0,0) -- (1,0,0) -- (1,0,1) -- (0,0,1) -- (0,0,0);
\filldraw[fill=MyGreen!90,draw=MyGreen] (0,0,0) -- (1,0,0) -- (1,1,0) -- (0,1,0) -- (0,0,0);
\filldraw[fill=MyGreen!90,draw=MyGreen] (0,0,0) -- (0,1,0) -- (0,1,1) -- (0,0,1) -- (0,0,0);
\filldraw[fill=MyGreen!60,draw=MyGreen] (1,0,0) -- (1,1,0) -- (1,1,1) -- (1,0,1) -- (1,0,0);
\filldraw[fill=MyGreen!75,draw=MyGreen] (0,0,1) -- (1,0,1) -- (1,1,1) -- (0,1,1) -- (0,0,1);
\filldraw[fill=MyGreen!90,draw=MyGreen] (0,1,0) -- (1,1,0) -- (1,1,1) -- (0,1,1) -- (0,1,0);

\end{scope}

\foreach \x in {1,...,3}{
\begin{scope}[shift={(-2,\x,0)}]
\pgfsetfillopacity{0.3}
\pgfsetstrokeopacity{0.4}
\pgfsetlinewidth{0.5pt}

\filldraw[fill=gray!70,draw=gray] (0,0,0) -- (1,0,0) -- (1,0,1) -- (0,0,1) -- (0,0,0);
\filldraw[fill=gray!70,draw=gray] (0,0,0) -- (1,0,0) -- (1,1,0) -- (0,1,0) -- (0,0,0);
\filldraw[fill=gray!70,draw=gray] (0,0,0) -- (0,1,0) -- (0,1,1) -- (0,0,1) -- (0,0,0);
\filldraw[fill=gray!50,draw=gray] (1,0,0) -- (1,1,0) -- (1,1,1) -- (1,0,1) -- (1,0,0);
\filldraw[fill=gray!60,draw=gray] (0,0,1) -- (1,0,1) -- (1,1,1) -- (0,1,1) -- (0,0,1);
\filldraw[fill=gray!70,draw=gray] (0,1,0) -- (1,1,0) -- (1,1,1) -- (0,1,1) -- (0,1,0);

\end{scope}}

\foreach \x in {0,1}{
\begin{scope}[shift={(-1,\x,0)}]
\pgfsetfillopacity{0.8}
\pgfsetlinewidth{0.5pt}

\filldraw[fill=MyGreen!90,draw=MyGreen] (0,0,0) -- (1,0,0) -- (1,0,1) -- (0,0,1) -- (0,0,0);
\filldraw[fill=MyGreen!90,draw=MyGreen] (0,0,0) -- (1,0,0) -- (1,1,0) -- (0,1,0) -- (0,0,0);
\filldraw[fill=MyGreen!90,draw=MyGreen] (0,0,0) -- (0,1,0) -- (0,1,1) -- (0,0,1) -- (0,0,0);
\filldraw[fill=MyGreen!60,draw=MyGreen] (1,0,0) -- (1,1,0) -- (1,1,1) -- (1,0,1) -- (1,0,0);
\filldraw[fill=MyGreen!75,draw=MyGreen] (0,0,1) -- (1,0,1) -- (1,1,1) -- (0,1,1) -- (0,0,1);
\filldraw[fill=MyGreen!90,draw=MyGreen] (0,1,0) -- (1,1,0) -- (1,1,1) -- (0,1,1) -- (0,1,0);

\end{scope}}

\begin{scope}[shift={(-1,2,0)}]
\pgfsetfillopacity{0.3}
\pgfsetstrokeopacity{0.4}
\pgfsetlinewidth{0.5pt}

\filldraw[fill=gray!70,draw=gray] (0,0,0) -- (1,0,0) -- (1,0,1) -- (0,0,1) -- (0,0,0);
\filldraw[fill=gray!70,draw=gray] (0,0,0) -- (1,0,0) -- (1,1,0) -- (0,1,0) -- (0,0,0);
\filldraw[fill=gray!70,draw=gray] (0,0,0) -- (0,1,0) -- (0,1,1) -- (0,0,1) -- (0,0,0);
\filldraw[fill=gray!50,draw=gray] (1,0,0) -- (1,1,0) -- (1,1,1) -- (1,0,1) -- (1,0,0);
\filldraw[fill=gray!60,draw=gray] (0,0,1) -- (1,0,1) -- (1,1,1) -- (0,1,1) -- (0,0,1);
\filldraw[fill=gray!70,draw=gray] (0,1,0) -- (1,1,0) -- (1,1,1) -- (0,1,1) -- (0,1,0);

\end{scope}

\foreach \x in {0,1,...,2}{
\begin{scope}[shift={(0,\x,0)}]
\pgfsetfillopacity{0.8}
\pgfsetlinewidth{0.5pt}

\filldraw[fill=MyGreen!90,draw=MyGreen] (0,0,0) -- (1,0,0) -- (1,0,1) -- (0,0,1) -- (0,0,0);
\filldraw[fill=MyGreen!90,draw=MyGreen] (0,0,0) -- (1,0,0) -- (1,1,0) -- (0,1,0) -- (0,0,0);
\filldraw[fill=MyGreen!90,draw=MyGreen] (0,0,0) -- (0,1,0) -- (0,1,1) -- (0,0,1) -- (0,0,0);
\filldraw[fill=MyGreen!60,draw=MyGreen] (1,0,0) -- (1,1,0) -- (1,1,1) -- (1,0,1) -- (1,0,0);
\filldraw[fill=MyGreen!75,draw=MyGreen] (0,0,1) -- (1,0,1) -- (1,1,1) -- (0,1,1) -- (0,0,1);
\filldraw[fill=MyGreen!90,draw=MyGreen] (0,1,0) -- (1,1,0) -- (1,1,1) -- (0,1,1) -- (0,1,0);

\end{scope}}

\begin{scope}[shift={(-1,3,0)}]
\pgfsetfillopacity{0.3}
\pgfsetstrokeopacity{0.4}
\pgfsetlinewidth{0.5pt}

\filldraw[fill=gray!70,draw=gray] (0,0,0) -- (1,0,0) -- (1,0,1) -- (0,0,1) -- (0,0,0);
\filldraw[fill=gray!70,draw=gray] (0,0,0) -- (1,0,0) -- (1,1,0) -- (0,1,0) -- (0,0,0);
\filldraw[fill=gray!70,draw=gray] (0,0,0) -- (0,1,0) -- (0,1,1) -- (0,0,1) -- (0,0,0);
\filldraw[fill=gray!50,draw=gray] (1,0,0) -- (1,1,0) -- (1,1,1) -- (1,0,1) -- (1,0,0);
\filldraw[fill=gray!60,draw=gray] (0,0,1) -- (1,0,1) -- (1,1,1) -- (0,1,1) -- (0,0,1);
\filldraw[fill=gray!70,draw=gray] (0,1,0) -- (1,1,0) -- (1,1,1) -- (0,1,1) -- (0,1,0);

\end{scope}

\begin{scope}[shift={(0,3,0)}]
\pgfsetfillopacity{0.3}
\pgfsetstrokeopacity{0.4}
\pgfsetlinewidth{0.5pt}

\filldraw[fill=gray!70,draw=gray] (0,0,0) -- (1,0,0) -- (1,0,1) -- (0,0,1) -- (0,0,0);
\filldraw[fill=gray!70,draw=gray] (0,0,0) -- (1,0,0) -- (1,1,0) -- (0,1,0) -- (0,0,0);
\filldraw[fill=gray!70,draw=gray] (0,0,0) -- (0,1,0) -- (0,1,1) -- (0,0,1) -- (0,0,0);
\filldraw[fill=gray!50,draw=gray] (1,0,0) -- (1,1,0) -- (1,1,1) -- (1,0,1) -- (1,0,0);
\filldraw[fill=gray!60,draw=gray] (0,0,1) -- (1,0,1) -- (1,1,1) -- (0,1,1) -- (0,0,1);
\filldraw[fill=gray!70,draw=gray] (0,1,0) -- (1,1,0) -- (1,1,1) -- (0,1,1) -- (0,1,0);

\end{scope}

\end{scope}
\begin{scope}[shift={(0,10,0)}]

\begin{scope}[shift={(-2,0,0)}]
\pgfsetfillopacity{0.8}
\pgfsetlinewidth{0.5pt}

\filldraw[fill=MyGreen!90,draw=MyGreen] (0,0,0) -- (1,0,0) -- (1,0,1) -- (0,0,1) -- (0,0,0);
\filldraw[fill=MyGreen!90,draw=MyGreen] (0,0,0) -- (1,0,0) -- (1,1,0) -- (0,1,0) -- (0,0,0);
\filldraw[fill=MyGreen!90,draw=MyGreen] (0,0,0) -- (0,1,0) -- (0,1,1) -- (0,0,1) -- (0,0,0);
\filldraw[fill=MyGreen!60,draw=MyGreen] (1,0,0) -- (1,1,0) -- (1,1,1) -- (1,0,1) -- (1,0,0);
\filldraw[fill=MyGreen!75,draw=MyGreen] (0,0,1) -- (1,0,1) -- (1,1,1) -- (0,1,1) -- (0,0,1);
\filldraw[fill=MyGreen!90,draw=MyGreen] (0,1,0) -- (1,1,0) -- (1,1,1) -- (0,1,1) -- (0,1,0);
\end{scope}

\begin{scope}[shift={(-2,1,0)}]
\pgfsetfillopacity{0.3}
\pgfsetstrokeopacity{0.4}
\pgfsetlinewidth{0.5pt}

\filldraw[fill=gray!70,draw=gray] (0,0,0) -- (1,0,0) -- (1,0,1) -- (0,0,1) -- (0,0,0);
\filldraw[fill=gray!70,draw=gray] (0,0,0) -- (1,0,0) -- (1,1,0) -- (0,1,0) -- (0,0,0);
\filldraw[fill=gray!70,draw=gray] (0,0,0) -- (0,1,0) -- (0,1,1) -- (0,0,1) -- (0,0,0);
\filldraw[fill=gray!50,draw=gray] (1,0,0) -- (1,1,0) -- (1,1,1) -- (1,0,1) -- (1,0,0);
\filldraw[fill=gray!60,draw=gray] (0,0,1) -- (1,0,1) -- (1,1,1) -- (0,1,1) -- (0,0,1);
\filldraw[fill=gray!70,draw=gray] (0,1,0) -- (1,1,0) -- (1,1,1) -- (0,1,1) -- (0,1,0);

\end{scope}

\begin{scope}[shift={(-2,2,0)}]
\pgfsetfillopacity{0.3}
\pgfsetstrokeopacity{0.4}
\pgfsetlinewidth{0.5pt}

\filldraw[fill=gray!70,draw=gray] (0,0,0) -- (1,0,0) -- (1,0,1) -- (0,0,1) -- (0,0,0);
\filldraw[fill=gray!70,draw=gray] (0,0,0) -- (1,0,0) -- (1,1,0) -- (0,1,0) -- (0,0,0);
\filldraw[fill=gray!70,draw=gray] (0,0,0) -- (0,1,0) -- (0,1,1) -- (0,0,1) -- (0,0,0);
\filldraw[fill=gray!50,draw=gray] (1,0,0) -- (1,1,0) -- (1,1,1) -- (1,0,1) -- (1,0,0);
\filldraw[fill=gray!60,draw=gray] (0,0,1) -- (1,0,1) -- (1,1,1) -- (0,1,1) -- (0,0,1);
\filldraw[fill=gray!70,draw=gray] (0,1,0) -- (1,1,0) -- (1,1,1) -- (0,1,1) -- (0,1,0);

\end{scope}

\foreach \x in {0,1}{
\begin{scope}[shift={(-1,\x,0)}]
\pgfsetfillopacity{0.8}
\pgfsetlinewidth{0.5pt}

\filldraw[fill=MyGreen!90,draw=MyGreen] (0,0,0) -- (1,0,0) -- (1,0,1) -- (0,0,1) -- (0,0,0);
\filldraw[fill=MyGreen!90,draw=MyGreen] (0,0,0) -- (1,0,0) -- (1,1,0) -- (0,1,0) -- (0,0,0);
\filldraw[fill=MyGreen!90,draw=MyGreen] (0,0,0) -- (0,1,0) -- (0,1,1) -- (0,0,1) -- (0,0,0);
\filldraw[fill=MyGreen!60,draw=MyGreen] (1,0,0) -- (1,1,0) -- (1,1,1) -- (1,0,1) -- (1,0,0);
\filldraw[fill=MyGreen!75,draw=MyGreen] (0,0,1) -- (1,0,1) -- (1,1,1) -- (0,1,1) -- (0,0,1);
\filldraw[fill=MyGreen!90,draw=MyGreen] (0,1,0) -- (1,1,0) -- (1,1,1) -- (0,1,1) -- (0,1,0);

\end{scope}}

\begin{scope}[shift={(-1,2,0)}]
\pgfsetfillopacity{0.3}
\pgfsetstrokeopacity{0.4}
\pgfsetlinewidth{0.5pt}

\filldraw[fill=gray!70,draw=gray] (0,0,0) -- (1,0,0) -- (1,0,1) -- (0,0,1) -- (0,0,0);
\filldraw[fill=gray!70,draw=gray] (0,0,0) -- (1,0,0) -- (1,1,0) -- (0,1,0) -- (0,0,0);
\filldraw[fill=gray!70,draw=gray] (0,0,0) -- (0,1,0) -- (0,1,1) -- (0,0,1) -- (0,0,0);
\filldraw[fill=gray!50,draw=gray] (1,0,0) -- (1,1,0) -- (1,1,1) -- (1,0,1) -- (1,0,0);
\filldraw[fill=gray!60,draw=gray] (0,0,1) -- (1,0,1) -- (1,1,1) -- (0,1,1) -- (0,0,1);
\filldraw[fill=gray!70,draw=gray] (0,1,0) -- (1,1,0) -- (1,1,1) -- (0,1,1) -- (0,1,0);

\end{scope}

\foreach \x in {0,1,2}{
\begin{scope}[shift={(0,\x,0)}]
\pgfsetfillopacity{0.8}
\pgfsetlinewidth{0.5pt}

\filldraw[fill=MyGreen!90,draw=MyGreen] (0,0,0) -- (1,0,0) -- (1,0,1) -- (0,0,1) -- (0,0,0);
\filldraw[fill=MyGreen!90,draw=MyGreen] (0,0,0) -- (1,0,0) -- (1,1,0) -- (0,1,0) -- (0,0,0);
\filldraw[fill=MyGreen!90,draw=MyGreen] (0,0,0) -- (0,1,0) -- (0,1,1) -- (0,0,1) -- (0,0,0);
\filldraw[fill=MyGreen!60,draw=MyGreen] (1,0,0) -- (1,1,0) -- (1,1,1) -- (1,0,1) -- (1,0,0);
\filldraw[fill=MyGreen!75,draw=MyGreen] (0,0,1) -- (1,0,1) -- (1,1,1) -- (0,1,1) -- (0,0,1);
\filldraw[fill=MyGreen!90,draw=MyGreen] (0,1,0) -- (1,1,0) -- (1,1,1) -- (0,1,1) -- (0,1,0);
\end{scope}}
\end{scope}
\end{tikzpicture}
\caption{\label{fig:support_tensor} Graphical representation of the tensor from \Cref{exam:support} and its support. Pictorially we can view $\coim T_I$ as a subspace of $I^*$ by choosing a (possibly non-unique) coordinate subspace of $I^*$ mapping isomorphically under the map $I^* \twoheadrightarrow \coim T_I$.}
\end{figure}
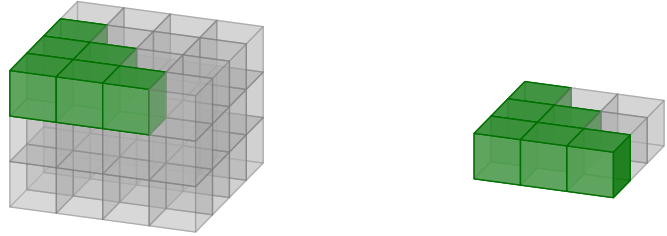

For a given 3-order tensor $T \in A\otimes B \otimes C$, the corresponding linear maps 
\[
T_{A}\colon A^*\to B\otimes C,\quad T_{B}\colon B^*\to A\otimes C,\quad T_{C}\colon C^*\to A\otimes B
\]
are called \textit{flattenings}. We  denote the images of these linear maps simply as
\[
T(A^*)\coloneqq T_A(A^*),\quad T(B^*)\coloneqq T_B(B^*),\quad T(C^*)\coloneqq T_C(C^*).
\]
The ranks of the flattenings $\mathbf{R}(T_I) = \dim T(I^*)$, for $I=A, B, C$, provide trivial lower and upper bounds for the tensor rank of~$T$. The proof of the following statement is discussed in \cite[Remark 10]{2011_Carlini}. 
\begin{proposition}\label{prop:trivial_lower_upper}
    For a tensor $T\in A\otimes B\otimes C$, let $\mathrm{r}_I\coloneqq\rk(T_I)$ be a flattening rank for any $I=A,B,C$.
    Then, we have
    \[
    \max\{\mathrm{r}_A,\mathrm{r}_B,\mathrm{r}_C\}\leq\rk (T)\leq\min\{\mathrm{r}_A\mathrm{r}_B,\mathrm{r}_A\mathrm{r}_C,\mathrm{r}_B\mathrm{r}_C\}.
    \]
\end{proposition}    
\begin{definition}
A tensor $T\in A\otimes B\otimes C$ is called \textit{$A$-concise} (\textit{$B$-concise} or \textit{$C$-concise}) if the map $T_A$ ($T_B$ or $T_C$ respectively) is injective.
Moreover, $T$ is called \textit{concise} if it is simultaneously $A$-concise, $B$-concise, and $C$-concise.    
\end{definition}
Recall that for any subspace $A'\subset A^*$, we have $(A')^*=A/(A')^{\perp}$. 
\begin{definition}\label{defn:support}
The \textit{support} of a tensor $T$ is the restriction
\[
T|_{(\coim T_A)\otimes (\coim T_B)\otimes (\coim T_C)}\in (\coim T_A)^*\otimes (\coim T_B)^*\otimes (\coim T_C)^*.
\]
Two tensors $T_1\in A_1\otimes B_1\otimes C_1$ and $T_2\in A_2\otimes B_2\otimes C_2$, with $A_j\subset A$, $B_j\subset B$, and $C_j\subset C$, for $j=1,2$, are \textit{equivalent} if they have the same support.
\end{definition}
\begin{remark}
    Observe that the support of a tensor $T\in A\otimes B\otimes C$ describes the restriction up to isomorphisms. In particular, considering the projection $\pi_I\colon I^*\to \coim T_I$, for each $I=A,B,C$, we can restrict the tensor $T$ to its support by choosing any linear spaces $A'\subseteq A^*$, $B'\subseteq B^*$, and $C'\subseteq C^*$ such that $I'\cong \coim T_I$ induced by $\pi_I$, obtaining $T|_{A'\otimes B'\otimes C'}$. In particular, two tensors $T_1\in A_1\otimes B_1\otimes C_1$ and $T_2\in A_2\otimes B_2\otimes C_2$ are equivalent if $T_1|_{A'\otimes B'\otimes C'}=T_2|_{A'\otimes B'\otimes C'}$.
\end{remark}
\begin{example}\label{exam:support}
Let $T\in \C^3\otimes \C^4\otimes \C^4$ be defined as 
%&\hphantom{{}={}}
\begin{align*}
T&= e_1\otimes e_1\otimes e_1+ 2 e_1\otimes e_2\otimes e_1 + e_1\otimes 2e_3\otimes e_1+e_1\otimes e_1\otimes e_2
+2e_1\otimes e_2\otimes e_2+ e_1\otimes e_1\otimes e_3\\
&=e_1\otimes (e_1+2e_2+2e_3)\otimes e_1+ e_1\otimes (e_1+2e_2)\otimes e_2+e_1\otimes e_1\otimes e_3.
\end{align*}
The support of the tensor $T$, illustrated in \autoref{fig:support_tensor} on the left, is the restriction 
$T|_{ \bbC\otimes \C^3 \otimes\bbC^3}$,
shown in the same figure on the right.
\end{example}
Now, we introduce the notion of modifications of tensors, which is essential in order to prove the existence of the counterexample guaranteed by \Cref{maintheorem}.
\begin{definition} \label{def: modofocation}
 Let $W_C\subset A\otimes B$ be a subspace.    Any sum
\[
    T+\sum_{i=1}^{n} w_i\otimes c_i,
\]
where $w_i\in W_C$ and $c_i \in C$ for every $i=1,\dots,n$,
is called a \textit{modification} of $T$ by $W_C$. The set of modifications of $T$ by $W_C$ is denoted by $T\bmod W_C$. If $W_B\subset A\otimes C$, a \textit{modification of $T$} by $W_B$ and $W_C$ is an element of the set
\[
T\bmod(W_B,W_C)\coloneqq\Set{T''\in T'\bmod W_B|T'\in T\bmod W_C}.
\]
If $W_A\subset B\otimes C$, a \textit{modification of $T$} by $W_A$, $W_B$ and $W_C$ is an element of the set
\[
T\bmod(W_A,W_B,W_C)\coloneqq\Set{T''\in T'\bmod W_A |T'\in T\bmod(W_B,W_C)}.
\]
The \textit{rank of $T$ modulo $W_A$, $W_B$, and $W_C$} is the smallest rank of a modification of $T$ by $W_A$, $W_B$, and $W_C$, and it is denoted by $\min\rk\bigl(T\bmod (W_A,W_B,W_C)\bigr)$.
\end{definition}
\begin{example}\label{exam:modification}
Let $T\in\bbC^8\otimes\bbC^5\otimes\bbC^3$, defined as
\[
T= \biggl(\sum_{i=1}^8 e_i\biggr)\otimes e_1\otimes e_2,
\]
and let $V=\langle e_1\otimes e_1, e_2\otimes e_2\rangle\subset \mathbb{C}^5 \otimes \mathbb{C}^3$. 
A tensor
\begin{align*}
T'={} & T+(e_3+e_5)\otimes e_1\otimes e_1+(e_6+e_7)\otimes e_2\otimes e_2\\
={}&(e_1+e_2+e_4+e_8)\otimes e_1\otimes e_2+(e_3+e_5)\otimes e_1\otimes(e_1+e_2)+(e_6+e_7)\otimes(e_1+e_2)\otimes e_2
\end{align*}
is a modification of $T$ by $V$ (see \Cref{fig:modification}).
\end{example}
\begin{figure}[hbt!]
\centering
\tdplotsetmaincoords{68}{110}
\begin{tikzpicture}[tdplot_main_coords,line join=round,scale=0.4]
\begin{scope}

\foreach \z in {3,6}{
\pgfsetfillopacity{0.6}
\pgfsetlinewidth{0.5pt}
\draw[amber,line cap=round] (0,0,\z) -- (3,0,\z);
\draw[amber,line cap=round] (0,0,\z) -- (0,5,\z);
\foreach \y in {0,1,...,5}{
\draw[amber,line cap=round] (0,\y,\z) -- (3,\y,\z);
\draw[amber,line cap=round] (0,\y,\z) -- (0,\y,\z+1);
}
\foreach \x in {0,1,...,3}{
\draw[amber,line cap=round] (\x,0,\z) -- (\x,5,\z);
\draw[amber,line cap=round] (\x,0,\z) -- (\x,0,\z+1);
}

\fill[amber!80] (0,0,\z) -- (3,0,\z) -- (3,5,\z) -- (0,5,\z) -- (0,0,\z);

\foreach \y in {0,1,2,3,4}{
\fill[amber!60] (0,\y,\z) -- (3,\y,\z) -- (3,\y,\z+1) -- (0,\y,\z+1) -- (0,\y,\z) ;
}
\foreach \y in {0,1,2}{
\fill[amber!40] (\y,0,\z) -- (\y,5,\z) -- (\y,5,\z+1) -- (\y,0,\z+1) -- (\y,0,\z) ;
}

\fill[amber!80] (0,0,\z+1) -- (3,0,\z+1) -- (3,5,\z+1) -- (0,5,\z+1) -- (0,0,\z+1);
\fill[amber!60] (0,5,\z) -- (3,5,\z) -- (3,5,\z+1) -- (0,5,\z+1) -- (0,5,\z) ;
\fill[amber!40] (3,0,\z) -- (3,5,\z) -- (3,5,\z+1) -- (3,0,\z+1) -- (3,0,\z) ;

\draw[amber,line cap=round] (0,5,\z) -- (3,5,\z);
\draw[amber,line cap=round] (3,0,\z) -- (3,5,\z);
\foreach \y in {0,1,...,5}{
\draw[amber,line cap=round] (0,\y,\z+1) -- (3,\y,\z+1);
\draw[amber,line cap=round] (3,\y,\z) -- (3,\y,\z+1);
}
\foreach \x in {0,1,...,3}{
\draw[amber,line cap=round] (\x,0,\z+1) -- (\x,5,\z+1);
\draw[amber,line cap=round] (\x,5,\z) -- (\x,5,\z+1);
}
}
\pgfsetfillopacity{1}
\node (b) at (2.8,-2.3,3.7) {$e_1\otimes e_1$};
\node (c) at (2.8,-2.3,6.7) {$e_2\otimes e_2$};

\end{scope}

\begin{scope}[xshift=10cm,yshift=1cm]

\pgfsetfillopacity{0.6}
\pgfsetlinewidth{0.5pt}
\foreach \y in {0,1,...,3}{
\draw[darkblue,line cap=round] (0,\y,0) -- (3,\y,0);
\draw[darkblue,line cap=round] (0,\y,0) -- (0,\y,8);
\draw[darkblue,line cap=round] (\y,0,0) -- (\y,5,0);
\draw[darkblue,line cap=round] (\y,0,0) -- (\y,0,8);
\draw[darkblue,line cap=round] (0,0,\y) -- (3,0,\y);
\draw[darkblue,line cap=round] (0,0,\y) -- (0,5,\y);
}
\foreach \y in {4,5}{
\draw[darkblue,line cap=round] (0,\y,0) -- (3,\y,0);
\draw[darkblue,line cap=round] (0,\y,0) -- (0,\y,8);
\draw[darkblue,line cap=round] (0,0,\y) -- (3,0,\y);
\draw[darkblue,line cap=round] (0,0,\y) -- (0,5,\y);
}
\foreach \y in {6,7,8}{
\draw[darkblue,line cap=round] (0,0,\y) -- (3,0,\y);
\draw[darkblue,line cap=round] (0,0,\y) -- (0,5,\y);
}

\foreach \y in {0,1,...,7}{
\fill[darkblue!70] (0,0,\y) -- (3,0,\y) -- (3,5,\y) -- (0,5,\y) -- (0,0,\y);
}
\foreach \y in {0,1,...,4}{
\fill[darkblue!50] (0,\y,0) -- (3,\y,0) -- (3,\y,8) -- (0,\y,8) -- (0,\y,0) ;
}
\foreach \y in {0,1,...,2}{
\fill[darkblue!30] (\y,0,0) -- (\y,5,0) -- (\y,5,8) -- (\y,0,8) -- (\y,0,0) ;
}

\fill[darkblue!70] (0,0,8) -- (3,0,8) -- (3,5,8) -- (0,5,8) -- (0,0,8);
\fill[darkblue!50] (0,5,0) -- (3,5,0) -- (3,5,8) -- (0,5,8) -- (0,5,0) ;
\fill[darkblue!30] (3,0,0) -- (3,5,0) -- (3,5,8) -- (3,0,8) -- (3,0,0) ;

\draw[darkblue,line cap=round] (0,5,0) -- (3,5,0);
\draw[darkblue,line cap=round] (3,0,0) -- (3,5,0);
\foreach \y in {0,1,...,3}{
\draw[darkblue,line cap=round] (0,\y,8) -- (3,\y,8);
\draw[darkblue,line cap=round] (3,\y,0) -- (3,\y,8);
\draw[darkblue,line cap=round] (\y,0,8) -- (\y,5,8);
\draw[darkblue,line cap=round] (\y,5,0) -- (\y,5,8);
\draw[darkblue,line cap=round] (0,5,\y) -- (3,5,\y);
\draw[darkblue,line cap=round] (3,0,\y) -- (3,5,\y);
}
\foreach \y in {4,5}{
\draw[darkblue,line cap=round] (0,\y,8) -- (3,\y,8);
\draw[darkblue,line cap=round] (3,\y,0) -- (3,\y,8);
\draw[darkblue,line cap=round] (0,5,\y) -- (3,5,\y);
\draw[darkblue,line cap=round] (3,0,\y) -- (3,5,\y);
}
\foreach \y in {6,7,8}{
\draw[darkblue,line cap=round] (0,5,\y) -- (3,5,\y);
\draw[darkblue,line cap=round] (3,0,\y) -- (3,5,\y);
}
\pgfsetfillopacity{1}
\node (d) at (3.7,-1.1,8.2) {$T$};
\end{scope}

\begin{scope}[xshift=20cm,yshift=1cm]

\pgfsetfillopacity{0.6}
\pgfsetlinewidth{0.5pt}
\draw[darkblue,line cap=round] (0,0,0) -- (3,0,0);
\draw[darkblue,line cap=round] (0,0,0) -- (0,5,0);
\foreach \y in {0,1,...,5}{
\draw[darkblue,line cap=round] (0,\y,0) -- (3,\y,0);
\draw[darkblue,line cap=round] (0,\y,0) -- (0,\y,2);
\draw[MyGreen,line cap=round] (0,\y,2) -- (0,\y,3);
\draw[darkblue,line cap=round] (0,\y,3) -- (0,\y,4);
\draw[MyGreen,line cap=round] (0,\y,4) -- (0,\y,7);
\draw[darkblue,line cap=round] (0,\y,7) -- (0,\y,8);
}
\foreach \y in {0,1,...,3}{
\draw[darkblue,line cap=round] (\y,0,0) -- (\y,5,0);
\draw[darkblue,line cap=round] (\y,0,0) -- (\y,0,2);
\draw[MyGreen,line cap=round] (\y,0,2) -- (\y,0,3);
\draw[darkblue,line cap=round] (\y,0,3) -- (\y,0,4);
\draw[MyGreen,line cap=round] (\y,0,4) -- (\y,0,7);
\draw[darkblue,line cap=round] (\y,0,7) -- (\y,0,8);
}
\foreach \y in {0,1}{
\draw[darkblue,line cap=round] (0,0,\y) -- (3,0,\y);
\draw[darkblue,line cap=round] (0,0,\y) -- (0,5,\y);
}
\foreach \y in {2,3,...,7}{
\draw[MyGreen,line cap=round] (0,0,\y) -- (3,0,\y);
\draw[MyGreen,line cap=round] (0,0,\y) -- (0,5,\y);
}
\draw[darkblue,line cap=round] (0,0,8) -- (3,0,8);
\draw[darkblue,line cap=round] (0,0,8) -- (0,5,8);

\foreach \y in {0,1}{
\fill[darkblue!70] (0,0,\y) -- (3,0,\y) -- (3,5,\y) -- (0,5,\y) -- (0,0,\y);
}
\foreach \y in {0,1,...,4}{
\fill[darkblue!50] (0,\y,0) -- (3,\y,0) -- (3,\y,2) -- (0,\y,2) -- (0,\y,0) ;
}
\foreach \y in {0,1,2}{
\fill[darkblue!30] (\y,0,0) -- (\y,5,0) -- (\y,5,2) -- (\y,0,2) -- (\y,0,0) ;
}

\fill[MyGreen!80] (0,0,2) -- (3,0,2) -- (3,5,2) -- (0,5,2) -- (0,0,2);
\foreach \y in {0,1,...,4}{
\fill[MyGreen!60] (0,\y,2) -- (3,\y,2) -- (3,\y,3) -- (0,\y,3) -- (0,\y,2) ;
}
\foreach \y in {0,1,2}{
\fill[MyGreen!40] (\y,0,2) -- (\y,5,2) -- (\y,5,3) -- (\y,0,3) -- (\y,0,2) ;
}

\fill[darkblue!70] (0,0,3) -- (3,0,3) -- (3,5,3) -- (0,5,3) -- (0,0,3);

\foreach \y in {0,1,...,4}{
\fill[darkblue!50] (0,\y,3) -- (3,\y,3) -- (3,\y,4) -- (0,\y,4) -- (0,\y,3) ;
}
\foreach \y in {0,1,2}{
\fill[darkblue!30] (\y,0,3) -- (\y,5,3) -- (\y,5,4) -- (\y,0,4) -- (\y,0,3) ;
}

\foreach \y in {4,5,6}{
\fill[MyGreen!80] (0,0,\y) -- (3,0,\y) -- (3,5,\y) -- (0,5,\y) -- (0,0,\y);
}
\foreach \y in {0,1,...,4}{
\fill[MyGreen!60] (0,\y,4) -- (3,\y,4) -- (3,\y,7) -- (0,\y,7) -- (0,\y,4) ;
}
\foreach \y in {0,1,2}{
\fill[MyGreen!40] (\y,0,4) -- (\y,5,4) -- (\y,5,7) -- (\y,0,7) -- (\y,0,4) ;
}

\foreach \y in {7}{
\fill[darkblue!70] (0,0,\y) -- (3,0,\y) -- (3,5,\y) -- (0,5,\y) -- (0,0,\y);
}

\foreach \y in {0,1,...,4}{
\fill[darkblue!50] (0,\y,7) -- (3,\y,7) -- (3,\y,8) -- (0,\y,8) -- (0,\y,7) ;
}

\foreach \y in {0,1,2}{
\fill[darkblue!30] (\y,0,7) -- (\y,5,7) -- (\y,5,8) -- (\y,0,8) -- (\y,0,7) ;
}

\fill[darkblue!70] (0,0,8) -- (3,0,8) -- (3,5,8) -- (0,5,8) -- (0,0,8);

\fill[darkblue!50] (0,5,0) -- (3,5,0) -- (3,5,2) -- (0,5,2) -- (0,5,0) ;
\fill[MyGreen!60] (0,5,2) -- (3,5,2) -- (3,5,3) -- (0,5,3) -- (0,5,2) ;
\fill[darkblue!50] (0,5,3) -- (3,5,3) -- (3,5,4) -- (0,5,4) -- (0,5,3) ;
\fill[MyGreen!60] (0,5,4) -- (3,5,4) -- (3,5,7) -- (0,5,7) -- (0,5,4) ;
\fill[darkblue!50] (0,5,7) -- (3,5,7) -- (3,5,8) -- (0,5,8) -- (0,5,7) ;

\fill[darkblue!30] (3,0,0) -- (3,5,0) -- (3,5,2) -- (3,0,2) -- (3,0,0) ;
\fill[MyGreen!40] (3,0,2) -- (3,5,2) -- (3,5,3) -- (3,0,3) -- (3,0,2) ;
\fill[darkblue!30] (3,0,3) -- (3,5,3) -- (3,5,4) -- (3,0,4) -- (3,0,3) ;
\fill[MyGreen!40] (3,0,4) -- (3,5,4) -- (3,5,7) -- (3,0,7) -- (3,0,4) ;
\fill[darkblue!30] (3,0,7) -- (3,5,7) -- (3,5,8) -- (3,0,8) -- (3,0,7) ;

\draw[darkblue,line cap=round] (0,5,0) -- (3,5,0);
\draw[darkblue,line cap=round] (3,0,0) -- (3,5,0);
\foreach \y in {0,1,...,5}{
\draw[darkblue,line cap=round] (0,\y,8) -- (3,\y,8);
\draw[darkblue,line cap=round] (3,\y,0) -- (3,\y,2);
\draw[MyGreen,line cap=round] (3,\y,2) -- (3,\y,3);
\draw[darkblue,line cap=round] (3,\y,3) -- (3,\y,4);
\draw[MyGreen,line cap=round] (3,\y,4) -- (3,\y,7);
\draw[darkblue,line cap=round] (3,\y,7) -- (3,\y,8);
}
\foreach \y in {0,1,...,3}{
\draw[darkblue,line cap=round] (\y,0,8) -- (\y,5,8);
\draw[darkblue,line cap=round] (\y,5,0) -- (\y,5,2);
\draw[MyGreen,line cap=round] (\y,5,2) -- (\y,5,3);
\draw[darkblue,line cap=round] (\y,5,3) -- (\y,5,4);
\draw[MyGreen,line cap=round] (\y,5,4) -- (\y,5,7);
\draw[darkblue,line cap=round] (\y,5,7) -- (\y,5,8);
}
\foreach \y in {0,1}{
\draw[darkblue,line cap=round] (0,5,\y) -- (3,5,\y);
\draw[darkblue,line cap=round] (3,0,\y) -- (3,5,\y);
}
\foreach \y in {2,3,...,7}{
\draw[MyGreen,line cap=round] (0,5,\y) -- (3,5,\y);
\draw[MyGreen,line cap=round] (3,0,\y) -- (3,5,\y);
}
\draw[darkblue,line cap=round] (0,5,8) -- (3,5,8);
\draw[darkblue,line cap=round] (3,0,8) -- (3,5,8);
\pgfsetfillopacity{1}
\node (d) at (3.7,-1.1,8.2) {$T'$};
\end{scope}
\end{tikzpicture}
\caption{\label{fig:modification} Graphical representation of the modification from \Cref{exam:modification}. Considering the canonical bases, the tensor $T'$ is obtained by adding to the tensor $T$, represented in blue, tensor products of the canonical basis of $\bbC^8$ and the generators of $V$, which are $e_1\otimes e_1$ and $e_2\otimes e_2$ colored in yellow. The \textit{modified} parts of $T'$ are shown in green.} 
\end{figure}
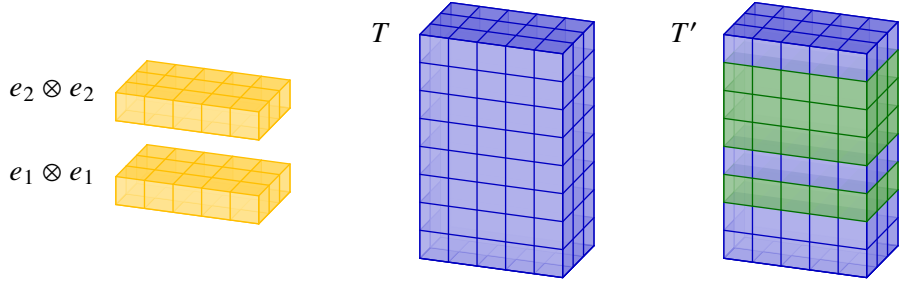
\subsection{Clones and augmented tensors}
Another key tool for the strategy we are going to use is the clone of a tensor, which we define here for a tensor of arbitrary order.
\begin{definition}\label{defn:sigma_clone}
Let $V$ be a vector space, let $\Sigma\coloneqq\{v_1,\dots,v_{\bfv}\}$ be a basis of $V$, and let
$S_{\Sigma}\coloneqq v_1+\cdots+v_{\bfv}$.
The \textit{$\Sigma$-clone} of a tensor $T\in A_1\otimes\cdots\otimes A_n$ is the tensor
\[
T^{\Sigma}\coloneqq T\boxtimes S_{\Sigma}^{\otimes n}\in(A_1\otimes V)\otimes\cdots\otimes(A_n\otimes V).
\]
\end{definition}
\begin{example}\label{exam:clone}
Let $T=e_1\otimes e_2\otimes e_1\in\bbC^2\otimes\bbC^2\otimes\bbC^2$.
If $\Sigma=\{e_1,e_2,e_3\}$ is the canonical basis of $\bbC^3$, then the $\Sigma$-clone of $T$ is a tensor in $\bigl(\bbC^2\otimes\bbC^3\bigl)^{\otimes 3}$ given by
\begin{align*}
T^{\Sigma}&=(e_1\otimes e_2\otimes e_1)\boxtimes (e_1+e_2+e_3)^{\otimes 3}=e_1\otimes(e_1+e_2+e_3)\otimes e_2\otimes (e_1+e_2+e_3)\otimes e_1\otimes (e_1+e_2+e_3)\\
&=\sum_{i,j,k=1}^3e_1\otimes e_i\otimes e_2\otimes e_j\otimes e_1\otimes e_k.
\end{align*}
It can be viewed as a block of several copies of the initial tensor $T$ along the three directions given by the elements of the basis $\Sigma$ (see \Cref{fig:clone}).
\end{example}
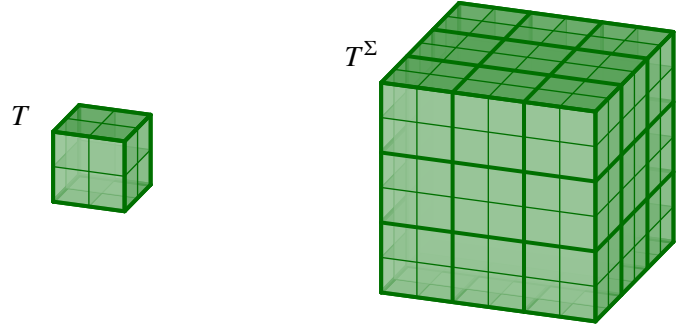
\begin{figure}[hbt!]
\centering
\tdplotsetmaincoords{68}{110}
\begin{tikzpicture}[tdplot_main_coords,line join=round,scale=0.5]
\begin{scope}
\pgfsetfillopacity{0.6}
\pgfsetlinewidth{1.5pt}
\foreach \y in {0,2}{
\draw[MyGreen,line cap=round] (0,\y,0) -- (2,\y,0);
\draw[MyGreen,line cap=round] (0,\y,0) -- (0,\y,2);
\draw[MyGreen,line cap=round] (\y,0,0) -- (\y,2,0);
\draw[MyGreen,line cap=round] (\y,0,0) -- (\y,0,2);
\draw[MyGreen,line cap=round] (0,0,\y) -- (2,0,\y);
\draw[MyGreen,line cap=round] (0,0,\y) -- (0,2,\y);
}
\pgfsetlinewidth{0.5pt}
\draw[MyGreen,line cap=round] (0,1,0) -- (2,1,0);
\draw[MyGreen,line cap=round] (0,1,0) -- (0,1,2);
\draw[MyGreen,line cap=round] (1,0,0) -- (1,2,0);
\draw[MyGreen,line cap=round] (1,0,0) -- (1,0,2);
\draw[MyGreen,line cap=round] (0,0,1) -- (2,0,1);
\draw[MyGreen,line cap=round] (0,0,1) -- (0,2,1);

\foreach \y in {0,1}{
\fill[MyGreen!80] (0,0,\y) -- (2,0,\y) -- (2,2,\y) -- (0,2,\y) -- (0,0,\y);
}
\foreach \y in {0,1}{
\fill[MyGreen!60] (0,\y,0) -- (2,\y,0) -- (2,\y,2) -- (0,\y,2) -- (0,\y,0) ;
}
\foreach \y in {0,1}{
\fill[MyGreen!40] (\y,0,0) -- (\y,2,0) -- (\y,2,2) -- (\y,0,2) -- (\y,0,0) ;
}

\fill[MyGreen!80] (0,0,2) -- (2,0,2) -- (2,2,2) -- (0,2,2) -- (0,0,2);
\fill[MyGreen!60] (0,2,0) -- (2,2,0) -- (2,2,2) -- (0,2,2) -- (0,2,0) ;
\fill[MyGreen!40] (2,0,0) -- (2,2,0) -- (2,2,2) -- (2,0,2) -- (2,0,0) ;

\pgfsetlinewidth{1.5pt}
\foreach \x in {0,2}{
\draw[MyGreen,line cap=round] (2,0,\x) -- (2,2,\x);
\draw[MyGreen,line cap=round] (2,\x,0) -- (2,\x,2);
\draw[MyGreen,line cap=round] (\x,2,0) -- (\x,2,2);
\draw[MyGreen,line cap=round] (0,2,\x) -- (2,2,\x);
\draw[MyGreen,line cap=round] (\x,0,2) -- (\x,2,2);
\draw[MyGreen,line cap=round] (0,\x,2) -- (2,\x,2);
}
\pgfsetlinewidth{0.5pt}
\draw[MyGreen,line cap=round] (0,1,2) -- (2,1,2);
\draw[MyGreen,line cap=round] (2,1,0) -- (2,1,2);
\draw[MyGreen,line cap=round] (1,0,2) -- (1,2,2);
\draw[MyGreen,line cap=round] (1,2,0) -- (1,2,2);
\draw[MyGreen,line cap=round] (0,2,1) -- (2,2,1);
\draw[MyGreen,line cap=round] (2,0,1) -- (2,2,1);

\pgfsetfillopacity{1}
\node (d) at (1.7,-1,2.2) {$T$};

\end{scope}

\begin{scope}[xshift=10cm,yshift=-1cm]
\pgfsetfillopacity{0.6}
\pgfsetlinewidth{1.5pt}
\foreach \y in {0,2,4,6}{
\draw[MyGreen,line cap=round] (0,\y,0) -- (6,\y,0);
\draw[MyGreen,line cap=round] (0,\y,0) -- (0,\y,6);
\draw[MyGreen,line cap=round] (\y,0,0) -- (\y,6,0);
\draw[MyGreen,line cap=round] (\y,0,0) -- (\y,0,6);
\draw[MyGreen,line cap=round] (0,0,\y) -- (6,0,\y);
\draw[MyGreen,line cap=round] (0,0,\y) -- (0,6,\y);
}
\pgfsetlinewidth{0.5pt}
\foreach \y in {1,3,5}{
\draw[MyGreen,line cap=round] (0,\y,0) -- (6,\y,0);
\draw[MyGreen,line cap=round] (0,\y,0) -- (0,\y,6);
\draw[MyGreen,line cap=round] (\y,0,0) -- (\y,6,0);
\draw[MyGreen,line cap=round] (\y,0,0) -- (\y,0,6);
\draw[MyGreen,line cap=round] (0,0,\y) -- (6,0,\y);
\draw[MyGreen,line cap=round] (0,0,\y) -- (0,6,\y);
}

\foreach \y in {0,1,...,5}{
\fill[MyGreen!80] (0,0,\y) -- (6,0,\y) -- (6,6,\y) -- (0,6,\y) -- (0,0,\y);
}
\foreach \y in {0,1,...,5}{
\fill[MyGreen!60] (0,\y,0) -- (6,\y,0) -- (6,\y,6) -- (0,\y,6) -- (0,\y,0) ;
}
\foreach \y in {0,1,...,5}{
\fill[MyGreen!40] (\y,0,0) -- (\y,6,0) -- (\y,6,6) -- (\y,0,6) -- (\y,0,0) ;
}

\fill[MyGreen!80] (0,0,6) -- (6,0,6) -- (6,6,6) -- (0,6,6) -- (0,0,6);
\fill[MyGreen!60] (0,6,0) -- (6,6,0) -- (6,6,6) -- (0,6,6) -- (0,6,0) ;
\fill[MyGreen!40] (6,0,0) -- (6,6,0) -- (6,6,6) -- (6,0,6) -- (6,0,0) ;

\pgfsetlinewidth{1.5pt}
\foreach \x in {0,2,4,6}{
\draw[MyGreen,line cap=round] (6,0,\x) -- (6,6,\x);
\draw[MyGreen,line cap=round] (6,\x,0) -- (6,\x,6);
\draw[MyGreen,line cap=round] (\x,6,0) -- (\x,6,6);
\draw[MyGreen,line cap=round] (0,6,\x) -- (6,6,\x);
\draw[MyGreen,line cap=round] (\x,0,6) -- (\x,6,6);
\draw[MyGreen,line cap=round] (0,\x,6) -- (6,\x,6);
}

\pgfsetlinewidth{0.5pt}
\foreach \y in {1,3,5}{
\draw[MyGreen,line cap=round] (0,\y,6) -- (6,\y,6);
\draw[MyGreen,line cap=round] (6,\y,0) -- (6,\y,6);
\draw[MyGreen,line cap=round] (\y,0,6) -- (\y,6,6);
\draw[MyGreen,line cap=round] (\y,6,0) -- (\y,6,6);
\draw[MyGreen,line cap=round] (0,6,\y) -- (6,6,\y);
\draw[MyGreen,line cap=round] (6,0,\y) -- (6,6,\y);
}

\pgfsetfillopacity{1}
\node (d) at (4.7,-1,6.2) {$T^{\Sigma}$};
\end{scope}
\end{tikzpicture}
\caption{\label{fig:clone} Graphical representation of the $\Sigma$-clone from \Cref{exam:clone}. The tensor $T^{\Sigma}$ can be viewed as a tensor obtained by taking $27$ copies of the initial tensor $T=e_1\otimes e_2\otimes e_1\in (\bbC^2)^{\otimes 3}$.  }
\end{figure}
Given a subspace $U\subset B\otimes C$, which is generated by some matrices in $B\otimes C$, it is quite natural to define an associated tensor by \textit{gluing} them together to form an order $3$ tensor, that is, a tensor in $\bbC^{\bfu}\otimes B\otimes C$, where $\bfu\coloneqq \dim U$.
\begin{definition}
    Let $U\subset B\otimes C$ and $\calB_U\coloneqq\{u_1,\dots,u_{\bfu}\}$ be a basis of $U$. The \textit{corresponding tensor} of $U$ is the tensor
    \[
    T_U\coloneqq \sum_{j=1}^{\bfu}e_j\otimes u_j\in \bbC^{\bfu}\otimes B\otimes C.
    \]
    Moreover, 
    for any vector space $V$ and a basis $\Sigma=\{v_1,\dots,v_{\bfv}\}$ the \textit{$\Sigma$-clone} of $U$ is the space
\[
U^{\Sigma}\coloneqq T_U^{\Sigma}\bigl((\bbC^{\bfu}\otimes V)^*\bigr)\subseteq(B\otimes V)\otimes(C\otimes V).
\] 
\end{definition}
\begin{figure}[hbt!]
\centering
\tdplotsetmaincoords{68}{110}
\begin{tikzpicture}[tdplot_main_coords,line join=round,scale=0.4]
\begin{scope}

\foreach \z in {0,3,6,9}{
\pgfsetfillopacity{0.6}
\pgfsetlinewidth{0.5pt}
\draw[darkblue,line cap=round] (0,0,\z) -- (5,0,\z);
\draw[darkblue,line cap=round] (0,0,\z) -- (0,5,\z);
\foreach \y in {0,1,...,5}{
\draw[darkblue,line cap=round] (0,\y,\z) -- (5,\y,\z);
\draw[darkblue,line cap=round] (0,\y,\z) -- (0,\y,\z+1);
}
\foreach \x in {0,1,...,5}{
\draw[darkblue,line cap=round] (\x,0,\z) -- (\x,5,\z);
\draw[darkblue,line cap=round] (\x,0,\z) -- (\x,0,\z+1);
}

\fill[darkblue!70] (0,0,\z) -- (5,0,\z) -- (5,5,\z) -- (0,5,\z) -- (0,0,\z);
\foreach \y in {0,1,...,4}{
\fill[darkblue!50] (0,\y,\z) -- (5,\y,\z) -- (5,\y,\z+1) -- (0,\y,\z+1) -- (0,\y,\z) ;
}
\foreach \y in {0,1,...,4}{
\fill[darkblue!30] (\y,0,\z) -- (\y,5,\z) -- (\y,5,\z+1) -- (\y,0,\z+1) -- (\y,0,\z) ;
}

\fill[darkblue!70] (0,0,\z+1) -- (5,0,\z+1) -- (5,5,\z+1) -- (0,5,\z+1) -- (0,0,\z+1);
\fill[darkblue!50] (0,5,\z) -- (5,5,\z) -- (5,5,\z+1) -- (0,5,\z+1) -- (0,5,\z) ;
\fill[darkblue!30] (5,0,\z) -- (5,5,\z) -- (5,5,\z+1) -- (5,0,\z+1) -- (5,0,\z) ;

\draw[darkblue,line cap=round] (0,5,\z) -- (5,5,\z);
\draw[darkblue,line cap=round] (5,0,\z) -- (5,5,\z);
\foreach \y in {0,1,...,5}{
\draw[darkblue,line cap=round] (0,\y,\z+1) -- (5,\y,\z+1);
\draw[darkblue,line cap=round] (5,\y,\z) -- (5,\y,\z+1);
}
\foreach \x in {0,1,...,5}{
\draw[darkblue,line cap=round] (\x,0,\z+1) -- (\x,5,\z+1);
\draw[darkblue,line cap=round] (\x,5,\z) -- (\x,5,\z+1);
}
}
\pgfsetfillopacity{1}
\node (a) at (3.7,-1.1,1.2) {$f_{{1}}$};
\node (b) at (3.7,-1.1,4.2) {$f_{{2}}$};
\node (c) at (3.7,-1.1,7.2) {$f_{{3}}$};
\node (d) at (3.7,-1.1,10.2) {$f_{{4}}$};

\end{scope}
\begin{scope}[xshift=12cm,yshift=2.5cm]
\pgfsetfillopacity{0.5}
\pgfsetlinewidth{0.5pt}

\foreach \y in {0,1,...,4}{
\draw[darkblue,line cap=round] (0,\y,0) -- (5,\y,0);
\draw[darkblue,line cap=round] (0,\y,0) -- (0,\y,4);
\draw[darkblue,line cap=round] (\y,0,0) -- (\y,5,0);
\draw[darkblue,line cap=round] (\y,0,0) -- (\y,0,4);
\draw[darkblue,line cap=round] (0,0,\y) -- (5,0,\y);
\draw[darkblue,line cap=round] (0,0,\y) -- (0,5,\y);
}
\draw[darkblue,line cap=round] (0,5,0) -- (5,5,0);
\draw[darkblue,line cap=round] (0,5,0) -- (0,5,4);
\draw[darkblue,line cap=round] (5,0,0) -- (5,5,0);
\draw[darkblue,line cap=round] (5,0,0) -- (5,0,4);

\foreach \y in {0,1,...,3}{
\fill[darkblue!70] (0,0,\y) -- (5,0,\y) -- (5,5,\y) -- (0,5,\y) -- (0,0,\y);
}
\foreach \y in {0,1,...,4}{
\fill[darkblue!50] (0,\y,0) -- (5,\y,0) -- (5,\y,4) -- (0,\y,4) -- (0,\y,0) ;
}
\foreach \y in {0,1,...,4}{
\fill[darkblue!30] (\y,0,0) -- (\y,5,0) -- (\y,5,4) -- (\y,0,4) -- (\y,0,0) ;
}

\fill[darkblue!70] (0,0,4) -- (5,0,4) -- (5,5,4) -- (0,5,4) -- (0,0,4);
\fill[darkblue!50] (0,5,0) -- (5,5,0) -- (5,5,4) -- (0,5,4) -- (0,5,0) ;
\fill[darkblue!30] (5,0,0) -- (5,5,0) -- (5,5,4) -- (5,0,4) -- (5,0,0) ;

\foreach \y in {0,1,...,4}{
\draw[darkblue,line cap=round] (0,\y,4) -- (5,\y,4);
\draw[darkblue,line cap=round] (5,\y,0) -- (5,\y,4);
\draw[darkblue,line cap=round] (\y,0,4) -- (\y,5,4);
\draw[darkblue,line cap=round] (\y,5,0) -- (\y,5,4);
\draw[darkblue,line cap=round] (0,5,\y) -- (5,5,\y);
\draw[darkblue,line cap=round] (5,0,\y) -- (5,5,\y);
}
\draw[darkblue,line cap=round] (0,5,4) -- (5,5,4);
\draw[darkblue,line cap=round] (5,5,0) -- (5,5,4);
\draw[darkblue,line cap=round] (5,0,4) -- (5,5,4);
\draw[darkblue,line cap=round] (5,5,0) -- (5,5,4);

\pgfsetfillopacity{1}
\node (d) at (3.7,-1.1,4.2) {$T_U$};
\end{scope}
\end{tikzpicture}
\caption{\label{fig:corresponding_tensor} Graphical representation of the corresponding tensor $T_U$ from \Cref{exam:corresponding_tensor}. The space $U$, generated by the four matrices $f_1,f_2,f_3,f_4\in\bbC^5\otimes \bbC^5$, can be associated to the tensor~$T_U$ obtained, once fixed a basis, as the union of the initial matrices.}
\end{figure}
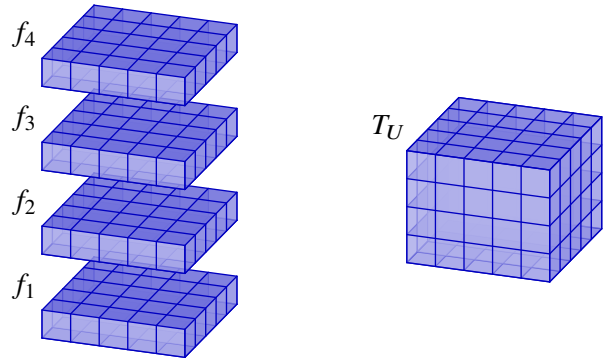
\begin{example}\label{exam:corresponding_tensor}
Let $U=\langle f_1,f_2,f_3,f_4\rangle$, where $f_1,f_2,f_3,f_4\in\bbC^5\otimes\bbC^5$ are four linearly independent matrices. The corresponding tensor of $U$, shown in \Cref{fig:corresponding_tensor}, is given  by the tensor
\[
T_U=e_1\otimes f_1+e_2\otimes f_2+e_3\otimes f_3+e_4\otimes f_4\in\bbC^4\otimes \bbC^5\otimes\bbC^5.\qedhere
\]
\end{example}
Considering how the corresponding tensors are defined, it is also natural to define the rank of a subspace of $B\otimes C$, and also its modifications.
\begin{definition}
The \textit{rank of a subspace} $U\subset B\otimes C$ is defined as the rank of its corresponding tensor, that is, as the value
$\rk(U)\coloneqq \rk(T_U)$.
\end{definition}
\begin{definition}
Let $U,M\subseteq B\otimes C$ be subspaces. Any sum
\[
T_U+\sum_{i=1}^{\bfu}e_i\otimes m_i,
\]
where $m_i\in M$ for every $i=1,\dots,\bfu$  is called a \textit{modification} of $U$ by $M$.
The \textit{rank of $U$ modulo $M$} is the smallest rank of a modification of $U$ by $M$, denoted by $\min\rk(T_U\bmod M)$.
\end{definition}

Observe that the definition of rank and modification is independent of the choice of basis used to define the corresponding tensor. 
\begin{definition}
Let $U_A\subseteq B\otimes C$, $U_B\subseteq A\otimes C$, and $U_C\subseteq A\otimes B$.
A tensor
\[
T+T_{U_A}+T_{U_B}+T_{U_C}\in (A\oplus\bbC^{\bfu_A})\otimes (B\oplus\bbC^{\bfu_B})\otimes (C\oplus\bbC^{\bfu_C})
\]
is said to be an \textit{augmented tensor} of $T$.
\end{definition}
\begin{example}\label{exam:augmented_tensor}
Let $A= \bbC^4$, $B= \bbC^5$, $C= \bbC^3$, and let $T= e_4\otimes e_5\otimes e_3\in \bbC^4\otimes \bbC^5\otimes \bbC^3$.
Consider 
\[
U_A=\langle f_1,f_2,f_3\rangle\subset \bbC^5\otimes\bbC^3,\quad U_B=\langle g_1,g_2,g_3\rangle\subset \bbC^4\otimes\bbC^3,\quad U_C=\langle h_1,h_2,h_3\rangle\subset \bbC^4\otimes\bbC^5
\]
generated by simple tensors
\begin{gather*}
f_1= e_5\otimes e_1,\quad f_2= e_5\otimes e_2,\quad f_3= e_5\otimes e_3,\quad g_1= e_4\otimes e_1,\quad g_2= e_4\otimes e_2,\\ g_3= e_4\otimes e_3,\quad
h_1= e_3\otimes e_1,\quad h_2= e_3\otimes e_2,\quad h_3= e_3\otimes e_3.
\end{gather*}
An augmented tensor
\begin{align*}
T+T_{U_A}+T_{U_B}+T_{U_C}&=e_4\otimes e_5\otimes e_3+e_5\otimes e_5\otimes e_1+e_6\otimes e_5\otimes e_2+e_7\otimes e_5\otimes e_3+e_4\otimes e_6\otimes e_1\\
&\hphantom{{}={}}+e_4\otimes e_7\otimes e_2+e_4\otimes e_8\otimes e_3+ e_3\otimes e_1\otimes e_4+ e_3\otimes e_2\otimes e_5+ e_3\otimes e_3\otimes e_6
\end{align*}
belongs to the space $(\bbC^4\oplus\bbC^{3})\otimes (\bbC^5\oplus\bbC^{3})\otimes (\bbC^3\oplus\bbC^{3})\simeq \bbC^7\otimes\bbC^8\otimes \bbC^6$.
\end{example}
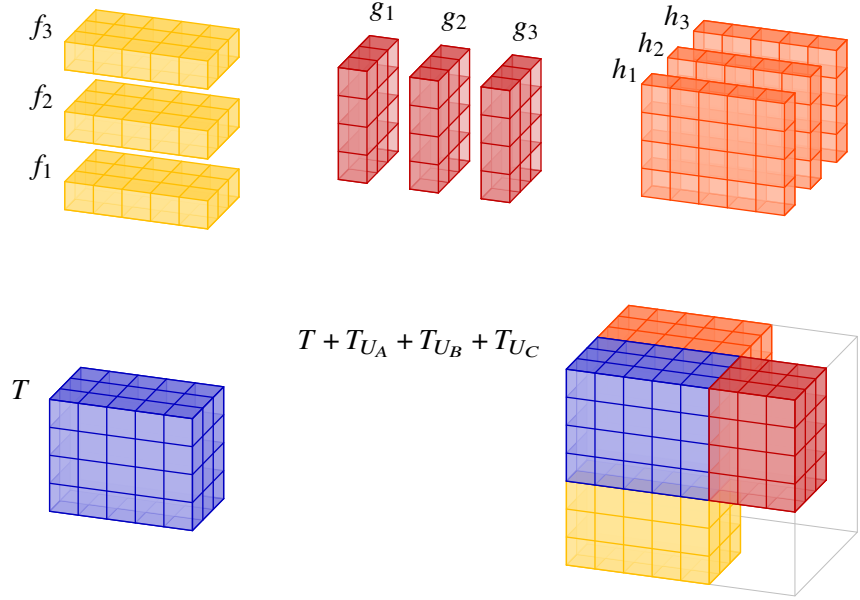
\begin{figure}[hbt!]
\centering
\tdplotsetmaincoords{68}{110}
\begin{tikzpicture}[tdplot_main_coords,line join=round,scale=0.4]
\begin{scope}[xshift=0.5cm]
\begin{scope}

\foreach \z in {0,2.5,5}{
\pgfsetfillopacity{0.6}
\pgfsetlinewidth{0.5pt}
\draw[amber,line cap=round] (0,0,\z) -- (3,0,\z);
\draw[amber,line cap=round] (0,0,\z) -- (0,5,\z);
\foreach \y in {0,1,...,5}{
\draw[amber,line cap=round] (0,\y,\z) -- (3,\y,\z);
\draw[amber,line cap=round] (0,\y,\z) -- (0,\y,\z+1);
}
\foreach \x in {0,1,...,3}{
\draw[amber,line cap=round] (\x,0,\z) -- (\x,5,\z);
\draw[amber,line cap=round] (\x,0,\z) -- (\x,0,\z+1);
}

\fill[amber!80] (0,0,\z) -- (3,0,\z) -- (3,5,\z) -- (0,5,\z) -- (0,0,\z);
\foreach \y in {0,1,...,4}{
\fill[amber!60] (0,\y,\z) -- (3,\y,\z) -- (3,\y,\z+1) -- (0,\y,\z+1) -- (0,\y,\z) ;
}
\foreach \y in {0,1,...,2}{
\fill[amber!40] (\y,0,\z) -- (\y,5,\z) -- (\y,5,\z+1) -- (\y,0,\z+1) -- (\y,0,\z) ;
}

\fill[amber!80] (0,0,\z+1) -- (3,0,\z+1) -- (3,5,\z+1) -- (0,5,\z+1) -- (0,0,\z+1);
\fill[amber!60] (0,5,\z) -- (3,5,\z) -- (3,5,\z+1) -- (0,5,\z+1) -- (0,5,\z) ;
\fill[amber!40] (3,0,\z) -- (3,5,\z) -- (3,5,\z+1) -- (3,0,\z+1) -- (3,0,\z) ;

\draw[amber,line cap=round] (0,5,\z) -- (3,5,\z);
\draw[amber,line cap=round] (3,0,\z) -- (3,5,\z);
\foreach \y in {0,1,...,5}{
\draw[amber,line cap=round] (0,\y,\z+1) -- (3,\y,\z+1);
\draw[amber,line cap=round] (3,\y,\z) -- (3,\y,\z+1);
}
\foreach \x in {0,1,...,3}{
\draw[amber,line cap=round] (\x,0,\z+1) -- (\x,5,\z+1);
\draw[amber,line cap=round] (\x,5,\z) -- (\x,5,\z+1);
}
}
\pgfsetfillopacity{1}
\node (a) at (2.5,-1,1.2) {$f_{{1}}$};
\node (b) at (2.5,-1,3.7) {$f_{{2}}$};
\node (c) at (2.5,-1,6.2) {$f_{{3}}$};

\end{scope}

\begin{scope}[xshift=20cm,yshift=1.5cm]

\foreach \z in {0,2.5,5}{
\pgfsetfillopacity{0.6}
\pgfsetlinewidth{0.5pt}

\foreach \x in {0,1,...,4}{
\draw[darkorange,line cap=round] (\z,0,\x) -- (\z,5,\x);
\draw[darkorange,line cap=round] (\z,0,\x) -- (\z+1,0,\x);
}
\foreach \y in {0,1,...,5}{
\draw[darkorange,line cap=round] (\z,\y,0) -- (\z,\y,4);
\draw[darkorange,line cap=round] (\z,\y,0) -- (\z+1,\y,0);
}
\draw[darkorange,line cap=round] (\z,0,0) -- (\z,0,4);
\draw[darkorange,line cap=round] (\z,0,0) -- (\z,5,0);

\foreach \y in {0,1,...,3}{
\fill[darkorange!80] (\z,0,\y) -- (\z,5,\y) -- (\z+1,5,\y) -- (\z+1,0,\y) -- (\z,0,\y) ;
}
\foreach \y in {0,1,...,4}{
\fill[darkorange!60] (\z,\y,0) -- (\z,\y,4) -- (\z+1,\y,4) -- (\z+1,\y,0) -- (\z,\y,0) ;
}
\fill[darkorange!40] (\z,0,0) -- (\z,0,4) -- (\z,5,4) -- (\z,5,0) -- (\z,0,0);

\fill[darkorange!80] (\z,0,4) -- (\z,5,4) -- (\z+1,5,4) -- (\z+1,0,4) -- (\z,0,4) ;
\fill[darkorange!60] (\z,5,0) -- (\z,5,4) -- (\z+1,5,4) -- (\z+1,5,0) -- (\z,5,0) ;
\fill[darkorange!40] (\z+1,0,0) -- (\z+1,0,4) -- (\z+1,5,4) -- (\z+1,5,0) -- (\z+1,0,0);

\draw[darkorange,line cap=round] (\z,5,0) -- (\z,5,4);
\draw[darkorange,line cap=round] (\z,0,4) -- (\z,5,4);
\foreach \y in {0,1,...,5}{
\draw[darkorange,line cap=round] (\z+1,\y,0) -- (\z+1,\y,4);
\draw[darkorange,line cap=round] (\z,\y,4) -- (\z+1,\y,4);
}
\foreach \x in {0,1,...,4}{
\draw[darkorange,line cap=round] (\z+1,0,\x) -- (\z+1,5,\x);
\draw[darkorange,line cap=round] (\z,5,\x) -- (\z+1,5,\x);
}
}
\pgfsetfillopacity{1}
\node (a) at (1.2,-0.5,4.5) {$h_{{3}}$};
\node (b) at (3.7,-0.5,4.5) {$h_{{2}}$};
\node (c) at (6.2,-0.5,4.5) {$h_{{1}}$};

\end{scope}

\begin{scope}[xshift=9cm,yshift=1cm]

\foreach \z in {0,2.5,5}{
\pgfsetfillopacity{0.6}
\pgfsetlinewidth{0.5pt}
\foreach \y in {0,1,...,4}{
\draw[darkred,line cap=round] (0,\z,\y) -- (3,\z,\y);
\draw[darkred,line cap=round] (0,\z,\y) -- (0,\z+1,\y);
}
\draw[darkred,line cap=round] (0,\z,0) -- (3,\z,0);
\draw[darkred,line cap=round] (0,\z,0) -- (0,\z,4);
\foreach \x in {0,1,...,3}{
\draw[darkred,line cap=round] (\x,\z,0) -- (\x,\z,4);
\draw[darkred,line cap=round] (\x,\z,0) -- (\x,\z+1,0);
}

\foreach \y in {0,1,...,3}{
\fill[darkred!80] (0,\z,\y) -- (3,\z,\y) -- (3,\z+1,\y) -- (0,\z+1,\y) -- (0,\z,\y) ;
}
\fill[darkred!60] (0,\z,0) -- (3,\z,0) -- (3,\z,4) -- (0,\z,4) -- (0,\z,0);
\foreach \y in {0,1,...,2}{
\fill[darkred!40] (\y,\z,0) -- (\y,\z,4) -- (\y,\z+1,4) -- (\y,\z+1,0) -- (\y,\z,0) ;
}

\fill[darkred!80] (0,\z,4) -- (3,\z,4) -- (3,\z+1,4) -- (0,\z+1,4) -- (0,\z,4) ;
\fill[darkred!60] (0,\z+1,0) -- (3,\z+1,0) -- (3,\z+1,0) -- (0,\z+1,0) -- (0,0,\z+1);
\fill[darkred!40] (3,\z,0) -- (3,\z,4) -- (3,\z+1,4) -- (3,\z+1,0) -- (3,\z,0) ;

\draw[darkred,line cap=round] (0,\z,4) -- (3,\z,4);
\draw[darkred,line cap=round] (3,\z,0) -- (3,\z,4);
\foreach \y in {0,1,...,4}{
\draw[darkred,line cap=round] (0,\z+1,\y) -- (3,\z+1,\y);
\draw[darkred,line cap=round] (3,\z,\y) -- (3,\z+1,\y);
}
\foreach \x in {0,1,...,3}{
\draw[darkred,line cap=round] (\x,\z+1,0) -- (\x,\z+1,4);
\draw[darkred,line cap=round] (\x,\z,4) -- (\x,\z+1,4);
}
}
\pgfsetfillopacity{1}
\node (a) at (-0.5,0.3,4.7) {$g_{{1}}$};
\node (b) at (-0.5,2.8,4.7) {$g_{{2}}$};
\node (c) at (-0.5,5.3,4.7) {$g_{{3}}$};

\end{scope}
\end{scope}

\begin{scope}[yshift=-10cm]
\begin{scope}[xshift=0]
\pgfsetfillopacity{0.6}
\pgfsetlinewidth{0.5pt}

\foreach \y in {0,1,...,4}{
\draw[darkblue,line cap=round] (0,0,\y) -- (3,0,\y);
\draw[darkblue,line cap=round] (0,0,\y) -- (0,5,\y);
}

\foreach \y in {0,1,...,5}{
\draw[darkblue,line cap=round] (0,\y,0) -- (3,\y,0);
\draw[darkblue,line cap=round] (0,\y,0) -- (0,\y,4);
}
\foreach \y in {0,1,...,3}{
\draw[darkblue,line cap=round] (\y,0,0) -- (\y,5,0);
\draw[darkblue,line cap=round] (\y,0,0) -- (\y,0,4);
}

\foreach \y in {0,1,...,3}{
\fill[darkblue!70] (0,0,\y) -- (3,0,\y) -- (3,5,\y) -- (0,5,\y) -- (0,0,\y);
}
\foreach \y in {0,1,...,4}{
\fill[darkblue!50] (0,\y,0) -- (3,\y,0) -- (3,\y,4) -- (0,\y,4) -- (0,\y,0) ;
}
\foreach \y in {0,1,...,2}{
\fill[darkblue!30] (\y,0,0) -- (\y,5,0) -- (\y,5,4) -- (\y,0,4) -- (\y,0,0) ;
}

\fill[darkblue!70] (0,0,4) -- (3,0,4) -- (3,5,4) -- (0,5,4) -- (0,0,4);
\fill[darkblue!50] (0,5,0) -- (3,5,0) -- (3,5,4) -- (0,5,4) -- (0,5,0);
\fill[darkblue!30] (3,0,0) -- (3,5,0) -- (3,5,4) -- (3,0,4) -- (3,0,0);

\foreach \y in {0,1,...,5}{
\draw[darkblue,line cap=round] (0,\y,4) -- (3,\y,4);
\draw[darkblue,line cap=round] (3,\y,0) -- (3,\y,4);
}
\foreach \y in {0,1,...,4}{
\draw[darkblue,line cap=round] (0,5,\y) -- (3,5,\y);
\draw[darkblue,line cap=round] (3,0,\y) -- (3,5,\y);
}
\foreach \y in {0,1,...,3}{
\draw[darkblue,line cap=round] (\y,0,4) -- (\y,5,4);
\draw[darkblue,line cap=round] (\y,5,0) -- (\y,5,4);
}

\pgfsetfillopacity{1}
\node (d) at (3,-1,4.2) {$T$};

\end{scope}

\begin{scope}[xshift=17cm,yshift=1cm,shift={(-3,0,0)}]

\pgfsetfillopacity{0.6}
\pgfsetlinewidth{0.5pt}
\foreach \y in {0,1,...,4}{
\draw[darkorange,line cap=round] (0,0,\y) -- (3,0,\y);
\draw[darkorange,line cap=round] (0,0,\y) -- (0,5,\y);
}

\foreach \y in {0,1,...,5}{
\draw[darkorange,line cap=round] (0,\y,0) -- (3,\y,0);
\draw[darkorange,line cap=round] (0,\y,0) -- (0,\y,4);
}
\foreach \y in {0,1,...,3}{
\draw[darkorange,line cap=round] (\y,0,0) -- (\y,5,0);
\draw[darkorange,line cap=round] (\y,0,0) -- (\y,0,4);
}

\foreach \y in {0,1,...,3}{
\fill[darkorange!80] (0,0,\y) -- (3,0,\y) -- (3,5,\y) -- (0,5,\y) -- (0,0,\y);
}
\foreach \y in {0,1,...,4}{
\fill[darkorange!60] (0,\y,0) -- (3,\y,0) -- (3,\y,4) -- (0,\y,4) -- (0,\y,0) ;
}
\foreach \y in {0,1,...,2}{
\fill[darkorange!40] (\y,0,0) -- (\y,5,0) -- (\y,5,4) -- (\y,0,4) -- (\y,0,0) ;
}

\fill[darkorange!80] (0,0,4) -- (3,0,4) -- (3,5,4) -- (0,5,4) -- (0,0,4);
\fill[darkorange!60] (0,5,0) -- (3,5,0) -- (3,5,4) -- (0,5,4) -- (0,5,0);
\fill[darkorange!40] (3,0,0) -- (3,5,0) -- (3,5,4) -- (3,0,4) -- (3,0,0);

\foreach \y in {0,1,...,5}{
\draw[darkorange,line cap=round] (0,\y,4) -- (3,\y,4);
\draw[darkorange,line cap=round] (3,\y,0) -- (3,\y,4);
}
\foreach \y in {0,1,...,4}{
\draw[darkorange,line cap=round] (0,5,\y) -- (3,5,\y);
\draw[darkorange,line cap=round] (3,0,\y) -- (3,5,\y);
}
\foreach \y in {0,1,...,3}{
\draw[darkorange,line cap=round] (\y,0,4) -- (\y,5,4);
\draw[darkorange,line cap=round] (\y,5,0) -- (\y,5,4);
}

\end{scope}

\begin{scope}[xshift=17cm,yshift=1cm,shift={(0,0,-3)}]

\pgfsetfillopacity{0.6}
\pgfsetlinewidth{0.3pt}
\draw[gray!50,line cap=round] (3,5,0) -- (3,8,0) -- (-3,8,0) -- (-3,0,0);
\draw[gray!50,line cap=round] (3,8,0) -- (3,8,3);
\draw[gray!50,line cap=round] (-3,8,0) -- (-3,8,7);
\draw[gray!50,line cap=round] (-3,5,7) -- (-3,8,7) -- (0,8,7);
\pgfsetfillopacity{0.6}
\pgfsetlinewidth{0.5pt}
\foreach \y in {0,1,...,3}{
\draw[amber,line cap=round] (0,0,\y) -- (3,0,\y);
\draw[amber,line cap=round] (0,0,\y) -- (0,5,\y);
}

\foreach \y in {0,1,...,5}{
\draw[amber,line cap=round] (0,\y,0) -- (3,\y,0);
\draw[amber,line cap=round] (0,\y,0) -- (0,\y,3);
}
\foreach \y in {0,1,...,3}{
\draw[amber,line cap=round] (\y,0,0) -- (\y,5,0);
\draw[amber,line cap=round] (\y,0,0) -- (\y,0,3);
}

\foreach \y in {0,1,...,2}{
\fill[amber!80] (0,0,\y) -- (3,0,\y) -- (3,5,\y) -- (0,5,\y) -- (0,0,\y);
}
\foreach \y in {0,1,...,4}{
\fill[amber!60] (0,\y,0) -- (3,\y,0) -- (3,\y,3) -- (0,\y,3) -- (0,\y,0) ;
}
\foreach \y in {0,1,...,2}{
\fill[amber!40] (\y,0,0) -- (\y,5,0) -- (\y,5,3) -- (\y,0,3) -- (\y,0,0) ;
}

\fill[amber!80] (0,0,3) -- (3,0,3) -- (3,5,3) -- (0,5,3) -- (0,0,3);
\fill[amber!60] (0,5,0) -- (3,5,0) -- (3,5,3) -- (0,5,3) -- (0,5,0);
\fill[amber!40] (3,0,0) -- (3,5,0) -- (3,5,3) -- (3,0,3) -- (3,0,0);

\foreach \y in {0,1,...,5}{
\draw[amber,line cap=round] (0,\y,3) -- (3,\y,3);
\draw[amber,line cap=round] (3,\y,0) -- (3,\y,3);
}
\foreach \y in {0,1,...,3}{
\draw[amber,line cap=round] (0,5,\y) -- (3,5,\y);
\draw[amber,line cap=round] (3,0,\y) -- (3,5,\y);
}
\foreach \y in {0,1,...,3}{
\draw[amber,line cap=round] (\y,0,3) -- (\y,5,3);
\draw[amber,line cap=round] (\y,5,0) -- (\y,5,3);
}

\end{scope}

\begin{scope}[xshift=17cm,yshift=1cm,shift={(0,0,0)}]
\pgfsetfillopacity{0.6}
\pgfsetlinewidth{0.5pt}

\foreach \y in {0,1,...,4}{
\draw[darkblue,line cap=round] (0,0,\y) -- (3,0,\y);
\draw[darkblue,line cap=round] (0,0,\y) -- (0,5,\y);
}

\foreach \y in {0,1,...,5}{
\draw[darkblue,line cap=round] (0,\y,0) -- (3,\y,0);
\draw[darkblue,line cap=round] (0,\y,0) -- (0,\y,4);
}
\foreach \y in {0,1,...,3}{
\draw[darkblue,line cap=round] (\y,0,0) -- (\y,5,0);
\draw[darkblue,line cap=round] (\y,0,0) -- (\y,0,4);
}

\foreach \y in {0,1,...,3}{
\fill[darkblue!70] (0,0,\y) -- (3,0,\y) -- (3,5,\y) -- (0,5,\y) -- (0,0,\y);
}
\foreach \y in {0,1,...,4}{
\fill[darkblue!50] (0,\y,0) -- (3,\y,0) -- (3,\y,4) -- (0,\y,4) -- (0,\y,0) ;
}
\foreach \y in {0,1,...,2}{
\fill[darkblue!30] (\y,0,0) -- (\y,5,0) -- (\y,5,4) -- (\y,0,4) -- (\y,0,0) ;
}

\fill[darkblue!70] (0,0,4) -- (3,0,4) -- (3,5,4) -- (0,5,4) -- (0,0,4);
\fill[darkblue!50] (0,5,0) -- (3,5,0) -- (3,5,4) -- (0,5,4) -- (0,5,0);
\fill[darkblue!30] (3,0,0) -- (3,5,0) -- (3,5,4) -- (3,0,4) -- (3,0,0);

\foreach \y in {0,1,...,5}{
\draw[darkblue,line cap=round] (0,\y,4) -- (3,\y,4);
\draw[darkblue,line cap=round] (3,\y,0) -- (3,\y,4);
}
\foreach \y in {0,1,...,4}{
\draw[darkblue,line cap=round] (0,5,\y) -- (3,5,\y);
\draw[darkblue,line cap=round] (3,0,\y) -- (3,5,\y);
}
\foreach \y in {0,1,...,3}{
\draw[darkblue,line cap=round] (\y,0,4) -- (\y,5,4);
\draw[darkblue,line cap=round] (\y,5,0) -- (\y,5,4);
}

\pgfsetfillopacity{1}
\node (d) at (3,-5.1,4.2) {$T+T_{U_A}+T_{U_B}+T_{U_C}$};

\end{scope}

\begin{scope}[xshift=17cm,yshift=1cm,shift={(0,5,0)}]
\pgfsetfillopacity{0.6}
\pgfsetlinewidth{0.5pt}

\foreach \y in {0,1,...,4}{
\draw[darkred,line cap=round] (0,0,\y) -- (3,0,\y);
\draw[darkred,line cap=round] (0,0,\y) -- (0,3,\y);
}

\foreach \y in {0,1,...,3}{
\draw[darkred,line cap=round] (0,\y,0) -- (3,\y,0);
\draw[darkred,line cap=round] (0,\y,0) -- (0,\y,4);
}
\foreach \y in {0,1,...,3}{
\draw[darkred,line cap=round] (\y,0,0) -- (\y,3,0);
\draw[darkred,line cap=round] (\y,0,0) -- (\y,0,4);
}

\foreach \y in {0,1,...,3}{
\fill[darkred!80] (0,0,\y) -- (3,0,\y) -- (3,3,\y) -- (0,3,\y) -- (0,0,\y);
}
\foreach \y in {0,1,...,2}{
\fill[darkred!60] (0,\y,0) -- (3,\y,0) -- (3,\y,4) -- (0,\y,4) -- (0,\y,0) ;
}
\foreach \y in {0,1,...,2}{
\fill[darkred!40] (\y,0,0) -- (\y,3,0) -- (\y,3,4) -- (\y,0,4) -- (\y,0,0) ;
}

\fill[darkred!80] (0,0,4) -- (3,0,4) -- (3,3,4) -- (0,3,4) -- (0,0,4);
\fill[darkred!60] (0,3,0) -- (3,3,0) -- (3,3,4) -- (0,3,4) -- (0,3,0);
\fill[darkred!40] (3,0,0) -- (3,3,0) -- (3,3,4) -- (3,0,4) -- (3,0,0);

\foreach \y in {0,1,...,3}{
\draw[darkred,line cap=round] (0,\y,4) -- (3,\y,4);
\draw[darkred,line cap=round] (3,\y,0) -- (3,\y,4);
}
\foreach \y in {0,1,...,4}{
\draw[darkred,line cap=round] (0,3,\y) -- (3,3,\y);
\draw[darkred,line cap=round] (3,0,\y) -- (3,3,\y);
}
\foreach \y in {0,1,...,3}{
\draw[darkred,line cap=round] (\y,0,4) -- (\y,3,4);
\draw[darkred,line cap=round] (\y,3,0) -- (\y,3,4);
}

\end{scope}
\end{scope}
\end{tikzpicture}
\caption{\label{fig:augmented_tensor} Graphical representation of the augmented tensor from \Cref{exam:augmented_tensor}. The generators of the spaces $U_A$, $U_B$, and $U_C$ are represented in yellow, orange, and red, respectively. Considering the canonical bases, the augmented tensor can be viewed as a tensor obtained by adding  blocks of the corresponding tensors of the three spaces $U_A$, $U_B$, and $U_C$.} 
\end{figure}
A typical procedure to obtain a lower bound is known as the \textit{substitution method}. Its proof relies on the following lemma, which is nothing but a simpler version of the same result (see also \cite[Proposition 2.3.3.2]{Lan19}).
\begin{lemma}
\label{lem:rank_mod_(00V)}
Let $T\in A\otimes B\otimes C$ and let $W_C \subset T(C^*)$. Then
\[\rk(T) \geq \min \rk\bigl(T \bmod(0, 0,W_C)\bigr) + \bfw_C.\]
Moreover, if $W_C$ is spanned by rank one elements, then equality holds.
\end{lemma}
\Cref{lem:rank_mod_(00V)} is obtained by applying the following proposition, which appeared for the first time in \cite{AFT11}, and is also presented in \cite[Proposition 3.1]{LM17}.
\begin{proposition}
\label{prop:LM17_result}
Let $T\in A\otimes B\otimes C$ such that $\rk (T)=r$, let $\{c_1,\dots,c_{\bfc}\}$ be a basis of $C$, and let 
\[ T=\sum_{i=1}^{\bfc}M_i\otimes c_i\] for some $M_1,\dots,M_{\bfc}\in A\otimes B$, with $M_1\neq 0$.
Then there exist constants $\lambda_2,\dots,\lambda_{\bfc}\in\bbC$ such that a tensor
\[
\hat{T}\coloneqq \sum_{j=2}^{\bfc}(M_j- \lambda_j M_1)\otimes c_j\in A\otimes B\otimes (c_1^\perp)^*
\]
has rank at most $r-1$. In particular, if $\rk( M_1)=1$, then $\rk(\hat{T})=r-1$.
\end{proposition}
The substitution method, first appeared in \cite{Pan66}, is described in detail in \cite[Section 3]{LM17} and it is an improvement of \Cref{prop:LM17_result},
whose proof relies 
on the following proposition. It appears for the first time in \cite[Proposition 14.45]{BCS97} in a weaker version as an inequality, while the version with the equality was given for the first time in \cite[Theorem 2.4]{Fri12} for tensors of order $3$, see also \cite[Theorem 3.1.1.1]{Lan12} and \cite[Proposition 2.1]{LM17}. It can be easily generalized to the tensor product of an arbitrary number of vector spaces, see \cite[Exercise 3.1.1.2]{Lan12}. 
\begin{proposition}
\label{prop:T(A)_contained_BC}
There exist $r$ rank one elements of $B\otimes C$ such that $T(A^*)$ is contained in their span if and only if $\rk (T)\leq r$. In particular, $\rk(T)$ equals the minimum number of rank one elements spanning~$T(A^*)$.
\end{proposition}
\Cref{lem:rank_mod_(00V)} is obtained by applying \Cref{prop:LM17_result} for a number of times equal to $\bfw_C$. 
The next corollary, which is essential later on, is a direct consequence of \Cref{lem:rank_mod_(00V)}. 
\begin{corollary}\label{cor2}
Let $U_A \subset B\otimes C$, $U_B \subset A\otimes C$, $U_C \subset A\otimes B$. Then, the augmented tensor \[T+T_{U_A} +T_{U_B} +T_{U_C} \in (A\oplus \C^{{\bf u}_A})\otimes(B\oplus \C^{{\bf u}_B})\otimes (C\oplus \C^{{\bf u}_C})\] satisfies the inequality
    \[\mathbf{R}(T+T_{U_A} +T_{U_B} +T_{U_C})\geq \min \mathbf{R}(T \text{ \rm mod }(U_A,U_B,U_C))+{\bf u}_A+{\bf u}_B+{\bf u}_C.\]
    Moreover, equality holds if $U_A,U_B,U_C$ are each spanned by rank one elements.
\end{corollary}
\section{Construction of the counterexample}\label{Sec3}
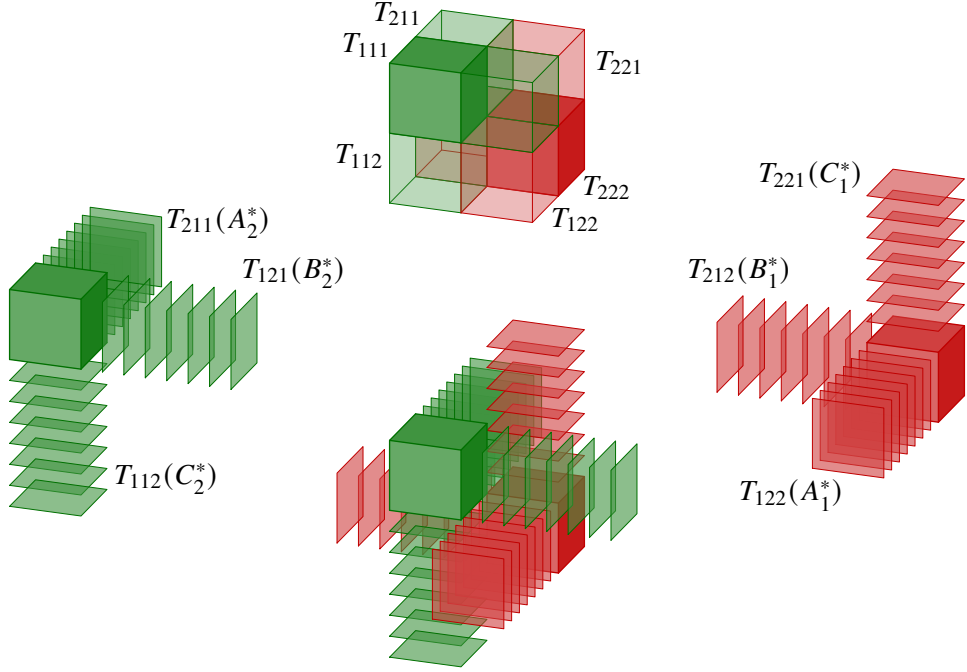
\begin{figure}[hbt!]
\centering
\tdplotsetmaincoords{68}{110}
\begin{tikzpicture}[tdplot_main_coords,line join=round,line cap=round,scale=1]
\begin{scope}
\begin{scope}[shift={(-1,0,-1)}]
\pgfsetfillopacity{0.2}
\pgfsetlinewidth{0.3pt}

\draw[color=darkred] (0,0,0) -- (1,0,0);
\draw[color=darkred] (0,0,0) -- (0,1,0);
\draw[color=darkred] (0,0,0) -- (0,0,1);

\fill[color=darkred!90] (0,0,0) -- (1,0,0) -- (1,0,1) -- (0,0,1);
\fill[color=darkred!90] (0,0,0) -- (1,0,0) -- (1,1,0) -- (0,1,0);
\fill[color=darkred!90] (0,0,0) -- (0,1,0) -- (0,1,1) -- (0,0,1);
\fill[color=darkred!60] (1,0,0) -- (1,1,0) -- (1,1,1) -- (1,0,1);
\fill[color=darkred!75] (0,0,1) -- (1,0,1) -- (1,1,1) -- (0,1,1);
\fill[color=darkred!90] (0,1,0) -- (1,1,0) -- (1,1,1) -- (0,1,1);

\draw[color=darkred] (1,0,0) -- (1,1,0);
\draw[color=darkred] (1,0,1) -- (1,1,1);
\draw[color=darkred] (0,0,1) -- (0,1,1);

\draw[color=darkred] (1,1,0) -- (0,1,0);
\draw[color=darkred] (1,1,1) -- (0,1,1);
\draw[color=darkred] (1,0,1) -- (0,0,1);

\draw[color=darkred] (1,0,0) -- (1,0,1);
\draw[color=darkred] (1,1,0) -- (1,1,1);
\draw[color=darkred] (0,1,0) -- (0,1,1);

\end{scope}

\begin{scope}[shift={(-1,1,-1)}]
\pgfsetfillopacity{1}
\pgfsetlinewidth{0.5pt}

\draw[line width=0.5pt,color=darkred] (0,0,0) -- (1,0,0);
\draw[line width=0.5pt,color=darkred] (0,0,0) -- (0,1,0);
\draw[line width=0.5pt,color=darkred] (0,0,0) -- (0,0,1);

\fill[color=darkred!90] (0,0,0) -- (1,0,0) -- (1,0,1) -- (0,0,1);
\fill[color=darkred!90] (0,0,0) -- (1,0,0) -- (1,1,0) -- (0,1,0);
\fill[color=darkred!90] (0,0,0) -- (0,1,0) -- (0,1,1) -- (0,0,1);
\fill[color=darkred!60] (1,0,0) -- (1,1,0) -- (1,1,1) -- (1,0,1);
\fill[color=darkred!75] (0,0,1) -- (1,0,1) -- (1,1,1) -- (0,1,1);
\fill[color=darkred!90] (0,1,0) -- (1,1,0) -- (1,1,1) -- (0,1,1);

\draw[color=darkred] (1,0,0) -- (1,1,0);
\draw[color=darkred] (1,0,1) -- (1,1,1);
\draw[color=darkred] (0,0,1) -- (0,1,1);

\draw[color=darkred] (1,1,0) -- (0,1,0);
\draw[color=darkred] (1,1,1) -- (0,1,1);
\draw[color=darkred] (1,0,1) -- (0,0,1);

\draw[color=darkred] (1,0,0) -- (1,0,1);
\draw[color=darkred] (1,1,0) -- (1,1,1);
\draw[color=darkred] (0,1,0) -- (0,1,1);

\node (d) at (0.5,1.5,0) {$T_{222}$};

\end{scope}

\begin{scope}[shift={(-1,1,0)}]
\pgfsetfillopacity{0.2}
\pgfsetlinewidth{0.3pt}

\draw[color=darkred] (0,0,0) -- (1,0,0);
\draw[color=darkred] (0,0,0) -- (0,1,0);
\draw[color=darkred] (0,0,0) -- (0,0,1);

\fill[color=darkred!90] (0,0,0) -- (1,0,0) -- (1,0,1) -- (0,0,1);
\fill[color=darkred!90] (0,0,0) -- (1,0,0) -- (1,1,0) -- (0,1,0);
\fill[color=darkred!90] (0,0,0) -- (0,1,0) -- (0,1,1) -- (0,0,1);
\fill[color=darkred!60] (1,0,0) -- (1,1,0) -- (1,1,1) -- (1,0,1);
\fill[color=darkred!75] (0,0,1) -- (1,0,1) -- (1,1,1) -- (0,1,1);
\fill[color=darkred!90] (0,1,0) -- (1,1,0) -- (1,1,1) -- (0,1,1);

\draw[color=darkred] (1,0,0) -- (1,1,0);
\draw[color=darkred] (1,0,1) -- (1,1,1);
\draw[color=darkred] (0,0,1) -- (0,1,1);

\draw[color=darkred] (1,1,0) -- (0,1,0);
\draw[color=darkred] (1,1,1) -- (0,1,1);
\draw[color=darkred] (1,0,1) -- (0,0,1);

\draw[color=darkred] (1,0,0) -- (1,0,1);
\draw[color=darkred] (1,1,0) -- (1,1,1);
\draw[color=darkred] (0,1,0) -- (0,1,1);

\pgfsetfillopacity{1}
\node (d) at (0,1.5,0.6) {$T_{221}$};
\end{scope}

\begin{scope}[shift={(-1,0,0)}]
\pgfsetfillopacity{0.2}
\pgfsetlinewidth{0.3pt}

\draw[color=MyGreen] (0,0,0) -- (1,0,0);
\draw[color=MyGreen] (0,0,0) -- (0,1,0);
\draw[color=MyGreen] (0,0,0) -- (0,0,1);

\fill[color=MyGreen!90] (0,0,0) -- (1,0,0) -- (1,0,1) -- (0,0,1);
\fill[color=MyGreen!90] (0,0,0) -- (1,0,0) -- (1,1,0) -- (0,1,0);
\fill[color=MyGreen!90] (0,0,0) -- (0,1,0) -- (0,1,1) -- (0,0,1);
\fill[color=MyGreen!60] (1,0,0) -- (1,1,0) -- (1,1,1) -- (1,0,1);
\fill[color=MyGreen!75] (0,0,1) -- (1,0,1) -- (1,1,1) -- (0,1,1);
\fill[color=MyGreen!90] (0,1,0) -- (1,1,0) -- (1,1,1) -- (0,1,1);

\draw[color=MyGreen] (1,0,0) -- (1,1,0);
\draw[color=MyGreen] (1,0,1) -- (1,1,1);
\draw[color=MyGreen] (0,0,1) -- (0,1,1);

\draw[color=MyGreen] (1,1,0) -- (0,1,0);
\draw[color=MyGreen] (1,1,1) -- (0,1,1);
\draw[color=MyGreen] (1,0,1) -- (0,0,1);

\draw[color=MyGreen] (1,0,0) -- (1,0,1);
\draw[color=MyGreen] (1,1,0) -- (1,1,1);
\draw[color=MyGreen] (0,1,0) -- (0,1,1);

\pgfsetfillopacity{1pt}
\node (d) at (0.3,-0.5,1) {$T_{211}$};
\end{scope}

\begin{scope}[shift={(0,0,-1)}]
\pgfsetfillopacity{0.2}
\pgfsetlinewidth{0.3pt}

\draw[color=MyGreen] (0,0,0) -- (1,0,0);
\draw[color=MyGreen] (0,0,0) -- (0,1,0);
\draw[color=MyGreen] (0,0,0) -- (0,0,1);

\fill[color=MyGreen!90] (0,0,0) -- (1,0,0) -- (1,0,1) -- (0,0,1);
\fill[color=MyGreen!90] (0,0,0) -- (1,0,0) -- (1,1,0) -- (0,1,0);
\fill[color=MyGreen!90] (0,0,0) -- (0,1,0) -- (0,1,1) -- (0,0,1);
\fill[color=MyGreen!60] (1,0,0) -- (1,1,0) -- (1,1,1) -- (1,0,1);
\fill[color=MyGreen!75] (0,0,1) -- (1,0,1) -- (1,1,1) -- (0,1,1);
\fill[color=MyGreen!90] (0,1,0) -- (1,1,0) -- (1,1,1) -- (0,1,1);

\draw[color=MyGreen] (1,0,0) -- (1,1,0);
\draw[color=MyGreen] (1,0,1) -- (1,1,1);
\draw[color=MyGreen] (0,0,1) -- (0,1,1);

\draw[color=MyGreen] (1,1,0) -- (0,1,0);
\draw[color=MyGreen] (1,1,1) -- (0,1,1);
\draw[color=MyGreen] (1,0,1) -- (0,0,1);

\draw[color=MyGreen] (1,0,0) -- (1,0,1);
\draw[color=MyGreen] (1,1,0) -- (1,1,1);
\draw[color=MyGreen] (0,1,0) -- (0,1,1);

\pgfsetfillopacity{1}
\node (d) at (0.5,-0.6,0.4) {$T_{112}$};
\end{scope}

\begin{scope}[shift={(0,1,-1)}]
\pgfsetfillopacity{0.2}
\pgfsetlinewidth{0.3pt}

\draw[color=darkred] (0,0,0) -- (1,0,0);
\draw[color=darkred] (0,0,0) -- (0,1,0);
\draw[color=darkred] (0,0,0) -- (0,0,1);

\fill[color=darkred!90] (0,0,0) -- (1,0,0) -- (1,0,1) -- (0,0,1);
\fill[color=darkred!90] (0,0,0) -- (1,0,0) -- (1,1,0) -- (0,1,0);
\fill[color=darkred!90] (0,0,0) -- (0,1,0) -- (0,1,1) -- (0,0,1);
\fill[color=darkred!60] (1,0,0) -- (1,1,0) -- (1,1,1) -- (1,0,1);
\fill[color=darkred!75] (0,0,1) -- (1,0,1) -- (1,1,1) -- (0,1,1);
\fill[color=darkred!90] (0,1,0) -- (1,1,0) -- (1,1,1) -- (0,1,1);

\draw[color=darkred] (1,0,0) -- (1,1,0);
\draw[color=darkred] (1,0,1) -- (1,1,1);
\draw[color=darkred] (0,0,1) -- (0,1,1);

\draw[color=darkred] (1,1,0) -- (0,1,0);
\draw[color=darkred] (1,1,1) -- (0,1,1);
\draw[color=darkred] (1,0,1) -- (0,0,1);

\draw[color=darkred] (1,0,0) -- (1,0,1);
\draw[color=darkred] (1,1,0) -- (1,1,1);
\draw[color=darkred] (0,1,0) -- (0,1,1);

\pgfsetfillopacity{1}
\node (d) at (0.7,1.5,0) {$T_{122}$};
\end{scope}

\begin{scope}
\pgfsetfillopacity{1}
\pgfsetlinewidth{0.5pt}

\draw[line width=0.5pt,color=MyGreen] (0,0,0) -- (1,0,0);
\draw[line width=0.5pt,color=MyGreen] (0,0,0) -- (0,1,0);
\draw[line width=0.5pt,color=MyGreen] (0,0,0) -- (0,0,1);

\fill[color=MyGreen!90] (0,0,0) -- (1,0,0) -- (1,0,1) -- (0,0,1);
\fill[color=MyGreen!90] (0,0,0) -- (1,0,0) -- (1,1,0) -- (0,1,0);
\fill[color=MyGreen!90] (0,0,0) -- (0,1,0) -- (0,1,1) -- (0,0,1);
\fill[color=MyGreen!60] (1,0,0) -- (1,1,0) -- (1,1,1) -- (1,0,1);
\fill[color=MyGreen!75] (0,0,1) -- (1,0,1) -- (1,1,1) -- (0,1,1);
\fill[color=MyGreen!90] (0,1,0) -- (1,1,0) -- (1,1,1) -- (0,1,1);

\draw[color=MyGreen] (1,0,0) -- (1,1,0);
\draw[color=MyGreen] (1,0,1) -- (1,1,1);
\draw[color=MyGreen] (0,0,1) -- (0,1,1);

\draw[color=MyGreen] (1,1,0) -- (0,1,0);
\draw[color=MyGreen] (1,1,1) -- (0,1,1);
\draw[color=MyGreen] (1,0,1) -- (0,0,1);

\draw[color=MyGreen] (1,0,0) -- (1,0,1);
\draw[color=MyGreen] (1,1,0) -- (1,1,1);
\draw[color=MyGreen] (0,1,0) -- (0,1,1);

\node (d) at (0.5,-0.5,1) {$T_{111}$};
\end{scope}

\begin{scope}[shift={(0,1,0)}]
\pgfsetfillopacity{0.2}
\pgfsetlinewidth{0.3pt}

\draw[color=MyGreen] (0,0,0) -- (1,0,0);
\draw[color=MyGreen] (0,0,0) -- (0,1,0);
\draw[color=MyGreen] (0,0,0) -- (0,0,1);

\fill[color=MyGreen!90] (0,0,0) -- (1,0,0) -- (1,0,1) -- (0,0,1);
\fill[color=MyGreen!90] (0,0,0) -- (1,0,0) -- (1,1,0) -- (0,1,0);
\fill[color=MyGreen!90] (0,0,0) -- (0,1,0) -- (0,1,1) -- (0,0,1);
\fill[color=MyGreen!60] (1,0,0) -- (1,1,0) -- (1,1,1) -- (1,0,1);
\fill[color=MyGreen!75] (0,0,1) -- (1,0,1) -- (1,1,1) -- (0,1,1);
\fill[color=MyGreen!90] (0,1,0) -- (1,1,0) -- (1,1,1) -- (0,1,1);

\draw[color=MyGreen] (1,0,0) -- (1,1,0);
\draw[color=MyGreen] (1,0,1) -- (1,1,1);
\draw[color=MyGreen] (0,0,1) -- (0,1,1);

\draw[color=MyGreen] (1,1,0) -- (0,1,0);
\draw[color=MyGreen] (1,1,1) -- (0,1,1);
\draw[color=MyGreen] (1,0,1) -- (0,0,1);

\draw[color=MyGreen] (1,0,0) -- (1,0,1);
\draw[color=MyGreen] (1,1,0) -- (1,1,1);
\draw[color=MyGreen] (0,1,0) -- (0,1,1);
\end{scope}
\end{scope}

\begin{scope}[yshift=-5cm]
\begin{scope}
\draw[line width=0.5pt,dash pattern=on 2pt off 1pt,color=darkred] (-1,1,-1) -- (0,1,-1);
\draw[line width=0.5pt,dash pattern=on 2pt off 1pt,color=darkred] (-1,1,-1) -- (-1,2,-1);
\draw[line width=0.5pt,dash pattern=on 2pt off 1pt,color=darkred] (-1,1,-1) -- (-1,1,0);

\foreach \x in {0.3,0.6,...,2.4}{
\pgfsetlinewidth{0.3pt}
\pgfsetfillopacity{0.6pt}
\fill[color=darkred!75] (0,1-\x,-1) -- (-1,1-\x,-1) -- (-1,1-\x,0) -- (0,1-\x,0) -- (0,1-\x,-1);
\draw[color=darkred] (0,1-\x,-1) -- (-1,1-\x,-1) -- (-1,1-\x,0) -- (0,1-\x,0) -- (0,1-\x,-1);
}

\begin{scope}[shift={(-1,1,-1)}]
\pgfsetfillopacity{1}
\pgfsetlinewidth{0.5pt}

\draw[line width=0.5pt,dash pattern=on 2pt off 1pt,color=darkred] (0,0,0) -- (1,0,0);
\draw[line width=0.5pt,dash pattern=on 2pt off 1pt,color=darkred] (0,0,0) -- (0,1,0);
\draw[line width=0.5pt,dash pattern=on 2pt off 1pt,color=darkred] (0,0,0) -- (0,0,1);

\fill[color=darkred!90] (0,0,0) -- (1,0,0) -- (1,0,1) -- (0,0,1);
\fill[color=darkred!90] (0,0,0) -- (1,0,0) -- (1,1,0) -- (0,1,0);
\fill[color=darkred!90] (0,0,0) -- (0,1,0) -- (0,1,1) -- (0,0,1);
\fill[color=darkred!60] (1,0,0) -- (1,1,0) -- (1,1,1) -- (1,0,1);
\fill[color=darkred!75] (0,0,1) -- (1,0,1) -- (1,1,1) -- (0,1,1);
\fill[color=darkred!90] (0,1,0) -- (1,1,0) -- (1,1,1) -- (0,1,1);

\draw[color=darkred] (1,0,0) -- (1,1,0);
\draw[color=darkred] (1,0,1) -- (1,1,1);
\draw[color=darkred] (0,0,1) -- (0,1,1);

\draw[color=darkred] (1,1,0) -- (0,1,0);
\draw[color=darkred] (1,1,1) -- (0,1,1);
\draw[color=darkred] (1,0,1) -- (0,0,1);

\draw[color=darkred] (1,0,0) -- (1,0,1);
\draw[color=darkred] (1,1,0) -- (1,1,1);
\draw[color=darkred] (0,1,0) -- (0,1,1);
\end{scope}

\foreach \x in {0.3,0.6,...,2.4}{
\pgfsetlinewidth{0.3pt}
\pgfsetfillopacity{0.6pt}
\fill[color=MyGreen!75] (-\x,0,0) -- (-\x,1,0) -- (-\x,1,1) -- (-\x,0,1) -- (-\x,0,0);
\draw[color=MyGreen] (-\x,0,0) -- (-\x,1,0) -- (-\x,1,1) -- (-\x,0,1) -- (-\x,0,0);
}

\foreach \x in {0.3,0.6,...,2.4}{
\pgfsetlinewidth{0.3pt}
\pgfsetfillopacity{0.6pt}
\fill[color=MyGreen!75] (1,0,-\x) -- (0,0,-\x) -- (0,1,-\x) -- (1,1,-\x) -- (1,0,-\x);
\draw[color=MyGreen] (1,0,-\x) -- (0,0,-\x) -- (0,1,-\x) -- (1,1,-\x) -- (1,0,-\x);
}

\foreach \x in {0.3,0.6,...,2.4}{
\pgfsetlinewidth{0.3pt}
\pgfsetfillopacity{0.6pt}
\fill[color=darkred!75] (0,1,\x) -- (0,2,\x) -- (-1,2,\x) -- (-1,1,\x) -- (0,1,\x);
\draw[color=darkred] (0,1,\x) -- (0,2,\x) -- (-1,2,\x) -- (-1,1,\x) -- (0,1,\x);
}

\foreach \x in {0.3,0.6,...,2.4}{
\pgfsetlinewidth{0.3pt}
\pgfsetfillopacity{0.6pt}
\fill[color=darkred!75] (\x,1,-1) -- (\x,2,-1) -- (\x,2,0) -- (\x,1,0) -- (\x,1,-1);
\draw[color=darkred] (\x,1,-1) -- (\x,2,-1) -- (\x,2,0) -- (\x,1,0) -- (\x,1,-1);
}

\begin{scope}
\pgfsetfillopacity{1}
\pgfsetlinewidth{0.5pt}

\draw[line width=0.5pt,dash pattern=on 2pt off 1pt,color=MyGreen] (0,0,0) -- (1,0,0);
\draw[line width=0.5pt,dash pattern=on 2pt off 1pt,color=MyGreen] (0,0,0) -- (0,1,0);
\draw[line width=0.5pt,dash pattern=on 2pt off 1pt,color=MyGreen] (0,0,0) -- (0,0,1);

\fill[color=MyGreen!90] (0,0,0) -- (1,0,0) -- (1,0,1) -- (0,0,1);
\fill[color=MyGreen!90] (0,0,0) -- (1,0,0) -- (1,1,0) -- (0,1,0);
\fill[color=MyGreen!90] (0,0,0) -- (0,1,0) -- (0,1,1) -- (0,0,1);
\fill[color=MyGreen!60] (1,0,0) -- (1,1,0) -- (1,1,1) -- (1,0,1);
\fill[color=MyGreen!75] (0,0,1) -- (1,0,1) -- (1,1,1) -- (0,1,1);
\fill[color=MyGreen!90] (0,1,0) -- (1,1,0) -- (1,1,1) -- (0,1,1);

\draw[color=MyGreen] (1,0,0) -- (1,1,0);
\draw[color=MyGreen] (1,0,1) -- (1,1,1);
\draw[color=MyGreen] (0,0,1) -- (0,1,1);

\draw[color=MyGreen] (1,1,0) -- (0,1,0);
\draw[color=MyGreen] (1,1,1) -- (0,1,1);
\draw[color=MyGreen] (1,0,1) -- (0,0,1);

\draw[color=MyGreen] (1,0,0) -- (1,0,1);
\draw[color=MyGreen] (1,1,0) -- (1,1,1);
\draw[color=MyGreen] (0,1,0) -- (0,1,1);
\end{scope}

\foreach \x in {0.3,0.6,...,2.4}{
\pgfsetlinewidth{0.3pt}
\pgfsetfillopacity{0.6}
\fill[color=MyGreen!75] (1,1+\x,0) -- (0,1+\x,0) -- (0,1+\x,1) -- (1,1+\x,1) -- (1,1+\x,0);
\draw[color=MyGreen] (1,1+\x,0) -- (0,1+\x,0) -- (0,1+\x,1) -- (1,1+\x,1) -- (1,1+\x,0);
}
\end{scope}
\begin{scope}[xshift=-5cm,yshift=2cm]
\foreach \x in {0.3,0.6,...,2.4}{
\pgfsetlinewidth{0.3pt}
\pgfsetfillopacity{0.6}
\fill[color=MyGreen!75] (-\x,0,0) -- (-\x,1,0) -- (-\x,1,1) -- (-\x,0,1) -- (-\x,0,0);
\draw[color=MyGreen] (-\x,0,0) -- (-\x,1,0) -- (-\x,1,1) -- (-\x,0,1) -- (-\x,0,0);
}

\pgfsetfillopacity{1}
\node (d) at (-1,2.2,1.5) {$T_{211}(A_2^*)$};

\foreach \x in {0.3,0.6,...,2.4}{
\pgfsetlinewidth{0.3pt}
\pgfsetfillopacity{0.6pt}
\fill[color=MyGreen!75] (1,0,-\x) -- (0,0,-\x) -- (0,1,-\x) -- (1,1,-\x) -- (1,0,-\x);
\draw[color=MyGreen] (1,0,-\x) -- (0,0,-\x) -- (0,1,-\x) -- (1,1,-\x) -- (1,0,-\x);
}
\pgfsetfillopacity{1}
\node (d) at (1,2.2,-1.4) {$T_{112}(C_2^*)$};

\begin{scope}
\pgfsetfillopacity{1}
\pgfsetlinewidth{0.5pt}

\draw[line width=0.5pt,dash pattern=on 2pt off 1pt,color=MyGreen] (0,0,0) -- (1,0,0);
\draw[line width=0.5pt,dash pattern=on 2pt off 1pt,color=MyGreen] (0,0,0) -- (0,1,0);
\draw[line width=0.5pt,dash pattern=on 2pt off 1pt,color=MyGreen] (0,0,0) -- (0,0,1);

\fill[color=MyGreen!90] (0,0,0) -- (1,0,0) -- (1,0,1) -- (0,0,1);
\fill[color=MyGreen!90] (0,0,0) -- (1,0,0) -- (1,1,0) -- (0,1,0);
\fill[color=MyGreen!90] (0,0,0) -- (0,1,0) -- (0,1,1) -- (0,0,1);
\fill[color=MyGreen!60] (1,0,0) -- (1,1,0) -- (1,1,1) -- (1,0,1);
\fill[color=MyGreen!75] (0,0,1) -- (1,0,1) -- (1,1,1) -- (0,1,1);
\fill[color=MyGreen!90] (0,1,0) -- (1,1,0) -- (1,1,1) -- (0,1,1);

\draw[color=MyGreen] (1,0,0) -- (1,1,0);
\draw[color=MyGreen] (1,0,1) -- (1,1,1);
\draw[color=MyGreen] (0,0,1) -- (0,1,1);

\draw[color=MyGreen] (1,1,0) -- (0,1,0);
\draw[color=MyGreen] (1,1,1) -- (0,1,1);
\draw[color=MyGreen] (1,0,1) -- (0,0,1);

\draw[color=MyGreen] (1,0,0) -- (1,0,1);
\draw[color=MyGreen] (1,1,0) -- (1,1,1);
\draw[color=MyGreen] (0,1,0) -- (0,1,1);
\end{scope}

\foreach \x in {0.3,0.6,...,2.4}{
\pgfsetlinewidth{0.3pt}
\pgfsetfillopacity{0.6pt}
\fill[color=MyGreen!75] (1,1+\x,0) -- (0,1+\x,0) -- (0,1+\x,1) -- (1,1+\x,1) -- (1,1+\x,0);
\draw[color=MyGreen] (1,1+\x,0) -- (0,1+\x,0) -- (0,1+\x,1) -- (1,1+\x,1) -- (1,1+\x,0);
}

\pgfsetfillopacity{1}
\node (d) at (-0.6,3.4,1.1) {$T_{121}(B_2^*)$};

\end{scope}
\begin{scope}[xshift=5cm,yshift=2cm]
\draw[line width=0.5pt,dash pattern=on 2pt off 1pt,color=darkred] (-1,1,-1) -- (0,1,-1);
\draw[line width=0.5pt,dash pattern=on 2pt off 1pt,color=darkred] (-1,1,-1) -- (-1,2,-1);
\draw[line width=0.5pt,dash pattern=on 2pt off 1pt,color=darkred] (-1,1,-1) -- (-1,1,0);

\foreach \x in {0.3,0.6,...,2.4}{
\pgfsetlinewidth{0.3pt}
\pgfsetfillopacity{0.6pt}
\fill[color=darkred!75] (0,1-\x,-1) -- (-1,1-\x,-1) -- (-1,1-\x,0) -- (0,1-\x,0) -- (0,1-\x,-1);
\draw[color=darkred] (0,1-\x,-1) -- (-1,1-\x,-1) -- (-1,1-\x,0) -- (0,1-\x,0) -- (0,1-\x,-1);
}
\pgfsetfillopacity{1}
\node (d) at (-0.6,-1,0.5) {$T_{212}(B_1^*)$};

\begin{scope}[shift={(-1,1,-1)}]
\pgfsetfillopacity{1}
\pgfsetlinewidth{0.5pt}

\draw[line width=0.5pt,dash pattern=on 2pt off 1pt,color=darkred] (0,0,0) -- (1,0,0);
\draw[line width=0.5pt,dash pattern=on 2pt off 1pt,color=darkred] (0,0,0) -- (0,1,0);
\draw[line width=0.5pt,dash pattern=on 2pt off 1pt,color=darkred] (0,0,0) -- (0,0,1);

\fill[color=darkred!90] (0,0,0) -- (1,0,0) -- (1,0,1) -- (0,0,1);
\fill[color=darkred!90] (0,0,0) -- (1,0,0) -- (1,1,0) -- (0,1,0);
\fill[color=darkred!90] (0,0,0) -- (0,1,0) -- (0,1,1) -- (0,0,1);
\fill[color=darkred!60] (1,0,0) -- (1,1,0) -- (1,1,1) -- (1,0,1);
\fill[color=darkred!75] (0,0,1) -- (1,0,1) -- (1,1,1) -- (0,1,1);
\fill[color=darkred!90] (0,1,0) -- (1,1,0) -- (1,1,1) -- (0,1,1);

\draw[color=darkred] (1,0,0) -- (1,1,0);
\draw[color=darkred] (1,0,1) -- (1,1,1);
\draw[color=darkred] (0,0,1) -- (0,1,1);

\draw[color=darkred] (1,1,0) -- (0,1,0);
\draw[color=darkred] (1,1,1) -- (0,1,1);
\draw[color=darkred] (1,0,1) -- (0,0,1);

\draw[color=darkred] (1,0,0) -- (1,0,1);
\draw[color=darkred] (1,1,0) -- (1,1,1);
\draw[color=darkred] (0,1,0) -- (0,1,1);
\end{scope}

\foreach \x in {0.3,0.6,...,2.4}{
\pgfsetlinewidth{0.3pt}
\pgfsetfillopacity{0.6pt}
\fill[color=darkred!75] (0,1,\x) -- (0,2,\x) -- (-1,2,\x) -- (-1,1,\x) -- (0,1,\x);
\draw[color=darkred] (0,1,\x) -- (0,2,\x) -- (-1,2,\x) -- (-1,1,\x) -- (0,1,\x);
}
\pgfsetfillopacity{1}
\node (d) at (-0.6,0,2) {$T_{221}(C_1^*)$};

\foreach \x in {0.3,0.6,...,2.4}{
\pgfsetlinewidth{0.3pt}
\pgfsetfillopacity{0.6pt}
\fill[color=darkred!75] (\x,1,-1) -- (\x,2,-1) -- (\x,2,0) -- (\x,1,0) -- (\x,1,-1);
\draw[color=darkred] (\x,1,-1) -- (\x,2,-1) -- (\x,2,0) -- (\x,1,0) -- (\x,1,-1);
}

\pgfsetfillopacity{1}
\node (d) at (-0.4,-0.2,-2.5) {$T_{122}(A_1^*)$};
\end{scope}
\end{scope}
\end{tikzpicture}
\caption{A graphical representation of the 8 tensors appearing in the statement of \Cref{lemma2}. The tensors \( T_{111} \) and \( T_{222} \), colored green and red, respectively, are modified by the generators of the spaces \( T_{211}(A_2^*) \), \( T_{121}(B_2^*) \), and \( T_{112}(C_2^*) \), all in green, and \( T_{122}(A_1^*) \), \( T_{212}(B_1^*) \), and \( T_{221}(C_1^*) \), all in red.
}
\label{fig:Lemma_2}
\end{figure}

The refutation of Strassen's rank additivity conjecture is a consequence of the following two lemmas, which are proved in the following sections. The first of these results concerns the construction of three specific subspaces of certain clones of $B\otimes C$, $A\otimes C$, and $A\otimes B$, respectively.
\begin{lemma}[{\cite[Claim 5]{shitov}}] \label{lemma1}
Let 
$W_A \subset B\otimes C$, $W_B \subset A\otimes C$, $W_C \subset A\otimes B$. Then there exists a vector space $V$ with basis $\Sigma$ and three spaces
\[
    \mathcal{M}_A \subset (B\otimes V) \otimes (C\otimes V), \quad
    \mathcal{M}_B \subset (A\otimes V) \otimes (C\otimes V), \quad
    \mathcal{M}_C \subset (A\otimes V) \otimes (B\otimes V),
\]
 such that:
 \begin{enumerate}[label=(\arabic*), left= 0pt, widest=*,nosep]
     \item $\mathcal{M}_A$, $\mathcal{M}_B$, and $\mathcal{M}_C$ are spanned by rank one tensors;
     \item $W_A^{\Sigma} \subset \mathcal{M}_A$, $W_B^{\Sigma} \subset \mathcal{M}_B$, $W_C^{\Sigma} \subset \mathcal{M}_C$;
     \item for any tensor $T\in A\otimes B\otimes C$, the $\Sigma$-clone $T^{\Sigma}$
     satisfies the equation 
      \[
     \min \rk\bigl(T^{\Sigma}  \bmod (\calM_A,\calM_B,\calM_C)\bigr)=\min \rk\bigl(T  \bmod (W_A,W_B,W_C)\bigr),
     \]
     and, in particular, the rank of the augmented tensor 
     $$T^{\Sigma} + T_{\mathcal{M}_A} + T_{\mathcal{M}_B} + T_{\mathcal{M}_C} \in \bigl((A\otimes {V}) \oplus \C^{{\bf m}_A}\bigr)\otimes\bigl((B\otimes {V}) \oplus \C^{{\bf m}_B}\bigr)\otimes\bigl((C\otimes{V}) \oplus \C^{{\bf m}_C}\bigr)$$
     is given by  
     \[\rk(T^{\Sigma} + T_{\mathcal{M}_A} +  T_{\mathcal{M}_B} + T_{\mathcal{M}_C}) =\min \rk\bigl(T  \bmod (W_A,W_B,W_C)\bigr)+{\bf m}_A+{\bf m}_B+{\bf m}_C.\]
 \end{enumerate}
\end{lemma}
We want to highlight the fact that the construction of the spaces $\calM_A$, $\calM_B$, and $\calM_C$ in \Cref{lemma1} can be considered for an arbitrary tensor $T\in A\otimes B\otimes C$.
\begin{lemma}[{\cite[Claim 6]{shitov}}] \label{lemma2}
    There exist $6$ vector spaces $A_1, A_2, B_1, B_2, C_1, C_2$ and 
    a tensor \[T\in(A_1 \oplus A_2) \otimes (B_1 \oplus B_2) \otimes (C_1 \oplus C_2),\] such that,
    \begin{align*}
        \rk(T) &<  \min \rk\Bigl(T_{111}\bmod\bigl(T_{211}(A_2^*),T_{121}(B_2^*),T_{112}(C_2^*)\bigr)\Bigr) \\[.5ex]
    &\hphantom{{}={}}    +  \min \rk\Bigl(T_{222}\bmod \bigl(T_{122}(A_1^*),T_{212}(B_1^*),T_{221}(C_1^*)\bigr)\Bigr),
    \end{align*}    
    where $T_{ijk}\coloneqq T|_{A_i^*\otimes B_j^*\otimes C_k^*}\in A_i\otimes B_j\otimes C_k$ for every
    $i,j,k=1,2$.
\end{lemma}
A graphical representation of the eight blocks of the tensor $T$ of \Cref{lemma2} is given in \Cref{fig:Lemma_2}.
Now we give the proof of \Cref{maintheorem}, which is based on \Cref{lemma1,lemma2}. Notice that the latter is an existence statement. The proof gives a recipe for constructing a counterexample, but not an explicit~one.

\begin{proof}[Proof of \Cref{maintheorem}]
 Let $T \in (A_1 \oplus A_2) \otimes (B_1 \oplus B_2) \otimes (C_1 \oplus C_2)$ be a tensor as in \Cref{lemma2}. We apply \Cref{lemma1} twice. The first time to the tensor $T_{111} \in A_1 \otimes B_1 \otimes C_1$ and linear spaces $T_{211}(A_2^{*})$,  $T_{121}(B_2^{*})$, and $T_{112}(C_2^{*})$, and the second time to the tensor $T_{222} \in A_2 \otimes B_2 \otimes C_2$ and  linear spaces $T_{122}(A_1^{*})$, $T_{212}(B_1^{*})$, and $T_{221}(C_1^{*})$. We get two spaces $V_1$ and $V_2$ together with two natural values $\sigma_1$, $\sigma_2\in\mathbb{Z}_{>0}$, such that 
 $V_1\simeq \bbC^{\sigma_1}$ and $V_2\simeq\bbC^{\sigma_2}$,
and spaces $ \mathcal{M}_{I_j}$, for each $I=A,B,C$ and every $j = 1,2$, spanned by rank one elements and satisfying the conditions of \Cref{lemma1}:
 \begin{align*}
      \mathcal{M}_{A_1} & \subset (B_1\otimes \C^{\sigma_1}) \otimes (C_1\otimes \C^{\sigma_1}), & \mathcal{M}_{A_2} \subset (B_2\otimes \C^{\sigma_2}) \otimes (C_2\otimes \C^{\sigma_2}),\\
    \mathcal{M}_{B_1} & \subset (A_1\otimes \C^{\sigma_1}) \otimes (C_1\otimes \C^{\sigma_1}), & \mathcal{M}_{B_2}  \subset (A_2\otimes \C^{\sigma_2}) \otimes (C_2\otimes \C^{\sigma_2}),\\
    \mathcal{M}_{C_1} & \subset (A_1\otimes \C^{\sigma_1}) \otimes (B_1\otimes \C^{\sigma_1}), &  \mathcal{M}_{C_2}  \subset (A_2\otimes \C^{\sigma_2}) \otimes (B_2\otimes \C^{\sigma_2}).
 \end{align*}
 In particular, by point (2) of \Cref{lemma1} we have
 \begin{align}\label{subspacesM}
\nonumber         T_{211}^{\Sigma}(A_2^{*}\otimes \C^{\sigma_1}) & \subset \mathcal{M}_{A_1},  \qquad T_{122}^{\Sigma}(A_1^{*}\otimes \C^{\sigma_2})  \subset \mathcal{M}_{A_2}, \qquad
            T_{121}^{\Sigma}(B_2^{*}\otimes \C^{\sigma_1})  \subset \mathcal{M}_{B_1},\\  
 T_{212}^{\Sigma}(B_1^{*}\otimes \C^{\sigma_2})  &\subset \mathcal{M}_{B_2},\qquad 
               T_{112}^{\Sigma}(C_2^{*}\otimes \C^{\sigma_1})  \subset \mathcal{M}_{C_1},  \qquad T_{221}^{\Sigma}(C_1^{*}\otimes \C^{\sigma_2})  \subset \mathcal{M}_{C_2}.
 \end{align}    
Up to considering equivalent tensors (see \Cref{defn:support}), we 
can assume that $\sigma_1 = \sigma_2 = \sigma$, just defining $\sigma\coloneqq\max\{\sigma_1,\sigma_2\}$. 
Let $\mathcal{T}_1$ and $\mathcal{T}_2$ be the augmented tensors
 \[
     \mathcal{T}_1   =  T_{111}^{\Sigma}  +T_{\mathcal{M}_{A_1}}  +T_{\mathcal{M}_{B_1}}  +T_{\mathcal{M}_{C_1}}, \qquad
     \mathcal{T}_2  =T_{222}^{\Sigma} +T_{\mathcal{M}_{A_2}} +T_{\mathcal{M}_{B_2}}  +T_{\mathcal{M}_{C_2}},
     \]
 where, in particular, setting $\mathbf{m}_{I_j}\coloneqq \dim \mathcal{M}_{I_j}$ for every $I = A, B, C$ and $j = 1,2$,
     \begin{align*}
      \mathcal{T}_1  & \in  \bigl((A_1 \otimes \C^{\sigma})\oplus \C^{{\bf m}_{A_1}}\bigr)\otimes\bigl((B_1 \otimes \C^{\sigma})\oplus \C^{{\bf m}_{B_1}}\bigr)\otimes\bigl((C_1 \otimes \C^{\sigma})\oplus \C^{{\bf m}_{C_1}}\bigr),\\
    \mathcal{T}_2 & \in  \bigl((A_2 \otimes \C^{\sigma})\oplus \C^{{\bf m}_{A_2}}\bigr)\otimes\bigl((B_2 \otimes \C^{\sigma})\oplus \C^{{\bf m}_{B_2}}\bigr)\otimes\bigl((C_2 \otimes \C^{\sigma})\oplus \C^{{\bf m}_{C_2}}\bigr).
 \end{align*}
 By point (3) of \Cref{lemma1}, the ranks of the augmented tensors $\mathcal{T}_1$ and $\mathcal{T}_2$ are
 \begin{align}\label{ranksT1_T2}
\nonumber      \rk(\mathcal{T}_1) & =\min \rk\Bigl(T_{111}\bmod\bigl(T_{211}(A_2^*),T_{121}(B_2^*),T_{112}(C_2^*)\bigr)\Bigr)+{\bf m}_{A_1}+{\bf m}_{B_1}+{\bf m}_{C_1},\\[0.5ex]
       \rk(\mathcal{T}_2) & =\min \rk\Bigl(T_{222}\bmod \bigl(T_{122}(A_1^*),T_{212}(B_1^*),T_{221}(C_1^*)\bigr)\Bigr)+{\bf m}_{A_2}+{\bf m}_{B_2}+{\bf m}_{C_2}.
     \end{align}
Let us consider the spaces \[\mathcal{M}_A \coloneqq \bigl\langle  \mathcal{M}_{A_1}, \mathcal{M}_{A_2}\bigr\rangle,\quad \mathcal{M}_B \coloneqq \bigl\langle  \mathcal{M}_{B_1},  \mathcal{M}_{B_2}\bigr\rangle,\quad \mathcal{M}_C \coloneqq \bigl\langle \mathcal{M}_{C_1}, \mathcal{M}_{C_2}\bigr\rangle,\]
where, in particular,
 \begin{gather*}
     \mathcal{M}_A   \subset \bigl((B_1 \oplus B_2)\otimes \C^{\sigma}\bigr) \otimes \bigl((C_1 \oplus C_2)\otimes \C^{\sigma}\bigr),\qquad
     \mathcal{M}_B   \subset \bigl((A_1 \oplus A_2)\otimes \C^{\sigma}\bigr) \otimes \bigl((C_1 \oplus C_2)\otimes \C^{\sigma}\bigr),\\
     \mathcal{M}_C   \subset \bigl((A_1 \oplus A_2)\otimes \C^{\sigma}\bigr) \otimes \bigl((B_1 \oplus B_2)\otimes \C^{\sigma}\bigr),
 \end{gather*}
 of dimensions $\mathbf{m}_{A} = \mathbf{m}_{A_1}+\mathbf{m}_{A_2}$, $\mathbf{m}_{B} = \mathbf{m}_{B_1}+\mathbf{m}_{B_2}$, and $\mathbf{m}_{C} = \mathbf{m}_{C_1}+\mathbf{m}_{C_2}$, respectively. 
Then, $\mathcal{T}_1 \oplus \mathcal{T}_2$ is the augmented tensor  \[(T_{111}^{\Sigma} \oplus T_{222}^{\Sigma}) + T_{\mathcal{M}_A} + T_{\mathcal{M}_B} + T_{\mathcal{M}_C}.\] Since the subspaces $\mathcal{M}_A, \mathcal{M}_B, \mathcal{M}_C$ are spanned by rank one elements, by \Cref{cor2} we have
\begin{equation}\label{minR_equality}
  \rk(\mathcal{T}_1 \oplus \mathcal{T}_2) = \min \rk\bigl((T_{111}^{\Sigma} \oplus T_{222}^{\Sigma}) \bmod (\mathcal{M}_A, \mathcal{M}_B, \mathcal{M}_C)\bigr) + \mathbf{m}_{A} + \mathbf{m}_{B}  + \mathbf{m}_{C}.  
\end{equation}
Notice that \eqref{subspacesM} implies $T^{\Sigma} \in (T_{111}^{\Sigma} \oplus T_{222}^{\Sigma}) \bmod(\mathcal{M}_A, \mathcal{M}_B, \mathcal{M}_C)$. Thus,
\begin{equation}\label{main_inequality_1}
       \rk(T) = \rk(T^{\Sigma}) \geq \min \rk\bigl((T_{111}^{\Sigma} \oplus T_{222}^{\Sigma}) \bmod(\mathcal{M}_A, \mathcal{M}_B, \mathcal{M}_C)\bigr)= \rk(\mathcal{T}_1 \oplus \mathcal{T}_2) - \mathbf{m}_{A} - \mathbf{m}_{B}  - \mathbf{m}_{C},
    \end{equation}
where the last equality is obtained from \eqref{minR_equality}. On the other hand,  $T$ satisfies \Cref{lemma2}. Thus,
\begin{align}\label{main_inequality_2}
\nonumber        \rk(T) &<  \min \mathbf{R}\Bigl(T_{111}\bmod \bigl(T_{211}(A_2^*),T_{121}(B_2^*),T_{112}(C_2^*)\bigr)\Bigr) \\[1ex]
\nonumber &\hphantom{{}={}}   + \min \rk\Bigl(T_{222}\bmod\bigl(T_{122}(A_1^*),T_{212}(B_1^*),T_{221}(C_1^*)\bigr)\Bigr)\\[1ex]
\nonumber &= \rk(\mathcal{T}_1) - \mathbf{m}_{A_1} - \mathbf{m}_{B_1}  - \mathbf{m}_{C_1}
        +  \rk(\mathcal{T}_2) - \mathbf{m}_{A_2} - \mathbf{m}_{B_2}  - \mathbf{m}_{C_2} \\[1ex]
       & = \rk(\mathcal{T}_1) +  \rk(\mathcal{T}_2) - \mathbf{m}_{A} - \mathbf{m}_{B}  - \mathbf{m}_{C},
    \end{align}   
where the first equality comes from \eqref{ranksT1_T2}. Combining inequalities \eqref{main_inequality_1} and \eqref{main_inequality_2}, we have
\[\rk(\mathcal{T}_1 \oplus \mathcal{T}_2) - \mathbf{m}_{A} - \mathbf{m}_{B}  - \mathbf{m}_{C} \leq \rk(T) < \rk(\mathcal{T}_1) +  \rk(\mathcal{T}_2) - \mathbf{m}_{A} - \mathbf{m}_{B}  - \mathbf{m}_{C},\]
which implies \[\rk(\mathcal{T}_1 \oplus \mathcal{T}_2) < \rk(\mathcal{T}_1) +  \rk(\mathcal{T}_2),\] and completes the proof.
\end{proof}

\section{Rank of augmented tensors of clones}\label{Sec4}
In this section, we provide the proof of \Cref{lemma1}. All the following contents can be found in \cite[Section 2]{shitov}. We provide here a version that is more connected to the language of tensors.
We use induction on the sum of the dimensions of the three spaces $W_A$, $W_B$, $W_C$. The first case to analyze is the case where
 $W_A=W_B=0$ and $W_C=\langle w\rangle$ for some $w\in A\otimes B$. This basic 
 step is stated in the following proposition.
 \begin{proposition}\label{prop:induction_case1}
Let $w\in A\otimes B$. Then there exists a vector space $V$, with basis $\Sigma$, and a space 
\[
\calM\subset (A\otimes V)\otimes(B\otimes V)
\]
such that:   
\begin{enumerate}[label=(\arabic*), left= 0pt, widest=*,nosep]
     \item $\mathcal{M}$ is spanned by rank one elements;
     \item $w^{\Sigma} \in \mathcal{M}$;
     \item 
     for any tensor $T\in A\otimes B\otimes C$, the $\Sigma$-clone $T^{\Sigma}$
     satisfies the equation 
     \[
     \min \rk\bigl(T^{\Sigma}  \bmod (0,0,\calM)\bigr)=\min \rk\bigl(T  \bmod (0,0,\langle w\rangle)\bigr),
     \]
     and, in particular, the rank of the augmented tensor 
    $$T^{\Sigma}+T_{\mathcal{M}} \in (A\otimes {V})\otimes(B\otimes {V})\otimes\bigl((C\otimes{V}) \oplus \C^{{\bf m}}\bigr)$$
     is given by   
     \[
     \rk(T^{\Sigma} + T_{\mathcal{M}}) =\min \rk\bigl(T  \bmod (0,0,\langle w\rangle)\bigr)+{\bfm}.
     \]   
 \end{enumerate}
\end{proposition}
In order to prove \Cref{prop:induction_case1}, we need to construct a set of rank one matrices that span the desired space $\calM$. Let us assume that $\rk (w)=r$, so that we can write 
\[
w=\sum_{j=1}^ra_j\otimes b_j
\] 
and consider it as a matrix having the first $r$ elements of the main diagonal equal to $1$ and zeros everywhere else. That is, we can assume
\[
w=
\begin{bNiceArray}{cccc|cc}[margin]
\Block{3-4}<\LARGE>{\mathrm{id}_r} & & & & \Block{3-2}<\LARGE>{0} &\\
& & & & & \\
& & & & & \\
\hline
\Block{2-4}<\LARGE>{0} & & & & \Block{2-2}<\LARGE>{0} & \\
& & & & & \\
\end{bNiceArray}.    
\]
\subsection{Spaces of matrices spanning clones}
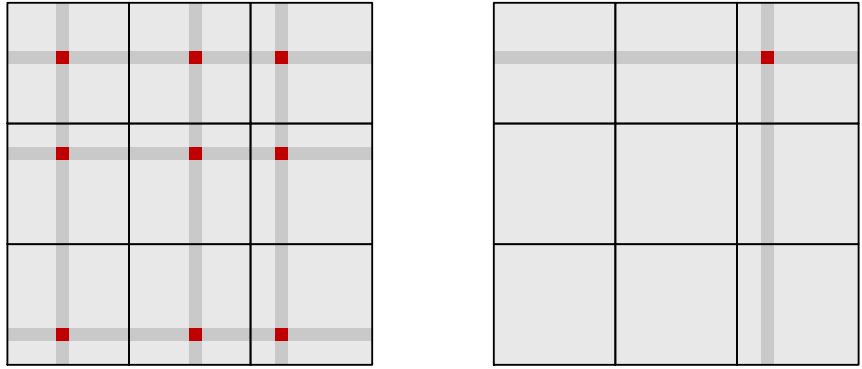
\begin{figure}[hbt!]
\centering
\begin{tikzpicture}[line cap=round,line join=round,scale=0.8]
\begin{scope}
\pgfsetlinewidth{0.7pt}

\pgfsetfillopacity{1}

\foreach \x in {0.8,3,4.4}{
\foreach \y in {0.4,3.4,5}{
\fill[color=gray!60] (0,\y) -- (6,\y) -- (6,\y+0.2) -- (0,\y+0.2) -- (0,\y);
\fill[color=gray!60] (\x,0) -- (\x,6) -- (\x+0.2,6) -- (\x+0.2,0) -- (\x,0);
}}

\pgfsetfillopacity{0.6}
\foreach \y in {0,2,4}{
\foreach \x in {0,2,4}{
\fill[color=gray!30] (\x,\y) -- (\x+2,\y) -- (\x+2,\y+2) -- (\x,\y+2) -- (\x,\y);
\draw[color=black] (\x,\y) -- (\x+2,\y) -- (\x+2,\y+2) -- (\x,\y+2) -- (\x,\y);
}}

\pgfsetfillopacity{1}

\foreach \x in {0.8,3,4.4}{
\foreach \y in {0.4,3.4,5}{
\fill[color=darkred] (\x,\y) -- (\x+0.2,\y) -- (\x+0.2,\y+0.2) -- (\x,\y+0.2) -- (\x,\y);
}}
\end{scope}
\begin{scope}[xshift=8cm]
\pgfsetlinewidth{0.7pt}

\pgfsetfillopacity{1}

\foreach \x in {4.4}{
\foreach \y in {5}{
\fill[color=gray!60] (0,\y) -- (6,\y) -- (6,\y+0.2) -- (0,\y+0.2) -- (0,\y);
\fill[color=gray!60] (\x,0) -- (\x,6) -- (\x+0.2,6) -- (\x+0.2,0) -- (\x,0);
}}

\pgfsetfillopacity{0.6}
\foreach \y in {0,2,4}{
\foreach \x in {0,2,4}{
\fill[color=gray!30] (\x,\y) -- (\x+2,\y) -- (\x+2,\y+2) -- (\x,\y+2) -- (\x,\y);
\draw[color=black] (\x,\y) -- (\x+2,\y) -- (\x+2,\y+2) -- (\x,\y+2) -- (\x,\y);
}}

\pgfsetfillopacity{1}

\foreach \x in {4.4}{
\foreach \y in {5}{
\fill[color=darkred] (\x,\y) -- (\x+0.2,\y) -- (\x+0.2,\y+0.2) -- (\x,\y+0.2) -- (\x,\y);
}}
\end{scope}
\end{tikzpicture}
\caption{Graphical representation of an example of a matrix in the set $\calF$ (on the left) and an example of a unit matrix in the set $\calU$ (on the right) in the case when $A=B=\bbC^3$. The only $1$ in the matrix on the right, represented in red, corresponds to a non-zero entry in a non-diagonal block of the matrix on the left.}
\label{fig:example_matrices_MPhi}
\end{figure}

Let us consider the values
\begin{equation}
\label{formula:initial_values}
\rho\coloneqq 2\bfa\bfb\bfc+1,\qquad \theta\coloneqq \lceil\log_2r\rceil,\qquad \sigma\coloneqq 2\rho^2 r.
\end{equation}
In particular, we have
$2^{\theta-1}< r\leq 2^{\theta}$ and
for any $\sigma\times\sigma$ matrix, the set of its $\sigma^2$ submatrices of size $1\times 1$ can be put in a one-to-one correspondence with the set $\{0,1,\dots,\sigma^2-1\}$. Indeed, any $s\in\{0,1,\dots,\sigma^2-1\}$ can be uniquely written as \[ s=\bigl(u(s,1)-1\bigr)\sigma+\bigl(u(s,2)-1\bigr),\] where $u(s,1)$, $u(s,2)\in\{1,\dots,\sigma\}$, that is, in base-$\sigma$.
\begin{definition}\label{def:matrix_M_Phi}
Let $\pi_i\colon\{1,\dots,\theta\}\to\{1,2\}$ be a function for every $i=1,\dots,r$, such that $\pi_i\neq \pi_j$ for any~$i\neq j$.  For any function 
$\Phi\colon\{1,\dots,\theta\}\to\{0,\dots,\sigma^2-1\}$,
the \textit{$\Phi$-tensor} associated to the set $\{\pi_i\}_{i=1,\dots,r}$ is a rank one tensor $M^{\Phi}\in A\otimes(\bbC^{\sigma})^{\otimes\theta}\otimes B\otimes (\bbC^{\sigma})^{\otimes\theta}$ defined as 
\begin{align*}
M^{\Phi}\coloneqq{} &\sum_{i,j=1}^r \biggl(a_i\otimes \bigotimes_{\gamma=1}^{\theta}e_{u(\Phi(\gamma),\pi_i(\gamma))}\otimes b_j\otimes \bigotimes_{\delta=1}^{\theta}e_{u(\Phi(\delta),3-\pi_j(\delta))}\biggr)\\
={}&\biggl(\sum_{i=1}^r a_i\otimes \bigotimes_{\gamma=1}^{\theta}e_{u(\Phi(\gamma),\pi_i(\gamma))}\biggr)\otimes \biggl(\sum_{j=1}^r b_j\otimes \bigotimes_{\delta=1}^{\theta}e_{u(\Phi(\delta),3-\pi_j(\delta))}\biggr).
\end{align*}
For every $i,j\in \{1,\dots,r\}$, the restriction
\begin{equation*}\label{formula:block}
M_{ij}^{\Phi}\coloneqq M^\Phi(a_i^*\otimes b_j^*)\in (\bbC^\sigma)^{\otimes \theta}\otimes(\bbC^{\sigma})^{\otimes \theta}
\end{equation*}
is called a \textit{block} of $M^{\Phi}$. If $i=j$, the $M_{ij}^\Phi$ is said to be a  \textit{diagonal block}.
\end{definition}
We denote by $\Sigma_{\theta}$ the canonical basis of $(\bbC^{\sigma})^{\otimes \theta}$, that is, 
\[
\Sigma_{\theta}\coloneqq\bigl\{e_{\alpha_1}\otimes\cdots\otimes e_{\alpha_{\theta}}\bigr\}_{\alpha_{1},\dots,\alpha_{\theta}\in\{1,\dots,\sigma\}}.
\]
When dealing with the tensor space $(\bbC^\sigma)^{\otimes \theta}\otimes(\bbC^{\sigma})^{\otimes \theta}$, we sometimes use  symbols $\varphi_{\alpha_1,\dots,\alpha_{\theta}}$ or $\psi_{\alpha_1,\dots,\alpha_{\theta}}$ to denote the element $e_{\alpha_1}\otimes\cdots\otimes e_{\alpha_{\theta}}$ in the first or in the second factor $(\bbC^\sigma)^{\otimes \theta}$, respectively.
We also define the set 
\[
\calF\coloneqq \Set{M^{\Phi}|\Phi\in {\{0,\dots,\sigma^2-1\}}^{\{1,\dots,\theta\}}},
\]
which is the set containing all the $\Phi$-tensors, and the set
\begin{equation} \label{formula:space_U}
\calU\coloneqq\Set{a_i\otimes\varphi \otimes b_j\otimes \psi|i\neq j,\,\varphi,\psi\in\Sigma_{\theta},\,\exists \Phi:M_{ij}^{\Phi}(\varphi\otimes\psi)\neq 0},    
\end{equation}
which is the set containing all those unit matrices that correspond to entries that are non-zero in at least one of the non-diagonal blocks in any of the matrices of $\calF$. Then we define the vector space
\begin{equation}\label{formula:space_M}
\calM\coloneqq \langle\calF\cup \calU\rangle\subset 
A\otimes(\bbC^\sigma)^{\otimes \theta}\otimes B\otimes(\bbC^{\sigma})^{\otimes \theta}.
\end{equation}
Any tensor $M^{\Phi}$ can be viewed as a rank one matrix divided into $\bfa\bfb$ blocks, each of which is a $\sigma^\theta\times\sigma^\theta$ matrix (see \Cref{fig:example_matrices_MPhi}). 

\begin{figure}[hbt!]
\centering
\begin{tikzpicture}[line cap=round,line join=round,scale=0.8]
\pgfsetlinewidth{0.7pt}

\pgfsetfillopacity{1}

\pgfsetfillopacity{0.6}
\foreach \y in {0}{
\foreach \x in {0,2}{
\fill[color=gray!30] (\x,\y) -- (\x+2,\y) -- (\x+2,\y+2) -- (\x,\y+2) -- (\x,\y);
}}

\foreach \y in {2}{
\foreach \x in {0,4}{
\fill[color=gray!30] (\x,\y) -- (\x+2,\y) -- (\x+2,\y+2) -- (\x,\y+2) -- (\x,\y);
}}

\foreach \y in {4}{
\foreach \x in {2,4}{
\fill[color=gray!30] (\x,\y) -- (\x+2,\y) -- (\x+2,\y+2) -- (\x,\y+2) -- (\x,\y);
}}

\pgfsetfillopacity{1}

\foreach \y in {0,2,4}{
\foreach \x in {4-\y}{
\fill[color=darkred] (\x,\y) -- (\x+2,\y) -- (\x+2,\y+2) -- (\x,\y+2) -- (\x,\y);

}}

\foreach \x in {0,0.2,0.8,3,4.4,4.6,5.6}{
\foreach \y in {0.4,0.8,2.2,2.8,3.4,5}{
\fill[color=darkred] (\x,\y) -- (\x+0.2,\y) -- (\x+0.2,\y+0.2) -- (\x,\y+0.2) -- (\x,\y);
}}

\foreach \x in {0.2,0.8,3.4,4.8}{
\foreach \y in {0.4,1,3.6,4.2,5.4}{
\fill[color=darkred] (\x,\y) -- (\x+0.2,\y) -- (\x+0.2,\y+0.2) -- (\x,\y+0.2) -- (\x,\y);
}}

\foreach \x in {0,0.4,0.8,1.2,2.4,2.6,3,3.2,3.4,3.8,4.4,4.8,5,5.4,5.8}{
\foreach \y in {0.6,0.8,1,1.2,2.2,2.6,3.2,3.6,3.8,4.8,5.2,5.8}{
\fill[color=darkred] (\x,\y) -- (\x+0.2,\y) -- (\x+0.2,\y+0.2) -- (\x,\y+0.2) -- (\x,\y);
}}

\foreach \x in {1.2,1.4,2.4,2.6,3,4.4,3.2}{
\foreach \y in {0.6,1.6,2.2,3.8,5.8}{
\fill[color=darkred] (\x,\y) -- (\x+0.2,\y) -- (\x+0.2,\y+0.2) -- (\x,\y+0.2) -- (\x,\y);
}}

\foreach \y in {0,2,4}{
\foreach \x in {0,2,4}{
\draw[color=black] (\x,\y) -- (\x+2,\y) -- (\x+2,\y+2) -- (\x,\y+2) -- (\x,\y);
}}

\begin{scope}[xshift=8cm]
\pgfsetlinewidth{0.7pt}

\pgfsetfillopacity{1}

\pgfsetfillopacity{0.6}
\foreach \y in {0}{
\foreach \x in {0,2}{
\fill[color=gray!30] (\x,\y) -- (\x+2,\y) -- (\x+2,\y+2) -- (\x,\y+2) -- (\x,\y);
\draw[color=black] (\x,\y) -- (\x+2,\y) -- (\x+2,\y+2) -- (\x,\y+2) -- (\x,\y);
}}

\foreach \y in {2}{
\foreach \x in {0,4}{
\fill[color=gray!30] (\x,\y) -- (\x+2,\y) -- (\x+2,\y+2) -- (\x,\y+2) -- (\x,\y);
\draw[color=black] (\x,\y) -- (\x+2,\y) -- (\x+2,\y+2) -- (\x,\y+2) -- (\x,\y);
}}

\foreach \y in {4}{
\foreach \x in {2,4}{
\fill[color=gray!30] (\x,\y) -- (\x+2,\y) -- (\x+2,\y+2) -- (\x,\y+2) -- (\x,\y);
\draw[color=black] (\x,\y) -- (\x+2,\y) -- (\x+2,\y+2) -- (\x,\y+2) -- (\x,\y);
}}

\pgfsetfillopacity{1}

\foreach \y in {0,2,4}{
\foreach \x in {4-\y}{
\fill[color=darkred] (\x,\y) -- (\x+2,\y) -- (\x+2,\y+2) -- (\x,\y+2) -- (\x,\y);
\draw[color=black] (\x,\y) -- (\x+2,\y) -- (\x+2,\y+2) -- (\x,\y+2) -- (\x,\y);
}}
\end{scope}
\end{tikzpicture}
\caption{Graphical representation of the matrix $\sum_{\, \Phi}M^{\Phi}$ shown on the left and the $\Sigma$-clone of the initial matrix $w$ shown on the right.}
\label{fig:example_sum_matrices_MPhi}
\end{figure}
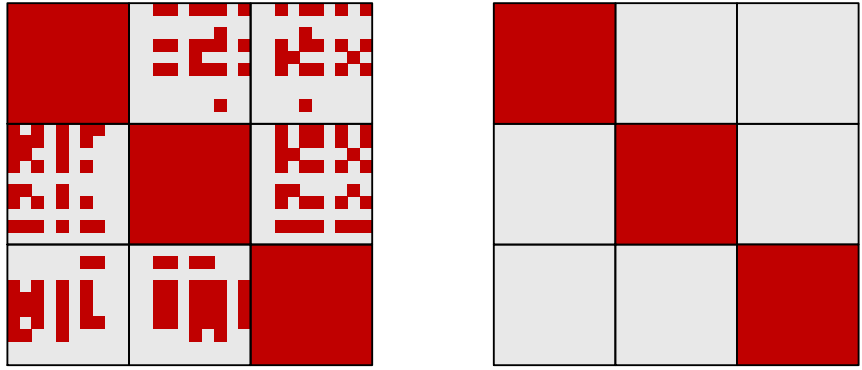

\begin{example}
\label[examples]{exam:second_example_M}
Let $r=3$ and let $A=\bbC^3$, $B= \bbC^4$, and $C=\bbC^2$. 
Take
$w=e_1\otimes e_1 + e_2\otimes e_2 + e_3\otimes e_3$
or, equivalently,
\[
w=\begin{pmatrix}
1&0&0&0\\
0&1&0&0\\
0&0&1&0
\end{pmatrix}.
\]
Then, we have
\[
\rho=49,\quad \theta=2,\quad \sigma=2\cdot 49^2\cdot 3=\numprint{14406}.
\]
Therefore, any $(a_i,b_j)$-block is a $\numprint{14406}^2\times\numprint{14406}^2$ matrix, that is, there are $12$ blocks formed by a square matrix with $(\numprint{207532836})^2$ entries. 
We define the three functions $\pi_1,\pi_2,\pi_3\colon\{1,2\}\to\{1,2\}$, such that
\begin{equation*}
\pi_1(1)=\pi_1(2)=1,\quad \pi_2(1)=\pi_2(2)=2,\quad \pi_3(1)=\pi_3(2)-1=1. 
\end{equation*}
We also consider the function 
\[
\Phi\colon\{1,2\}\to\{0,\dots,\numprint{207532835}\},
\]
defined by
$\Phi(1)=2\sigma+1$ and $\Phi(2)=4\sigma+3$.
We have
\begin{align*}
M^{\Phi}&=\bigl(e_1\otimes e_{u(\Phi(1),\pi_1(1))}\otimes e_{u(\Phi(2),\pi_1(2))}+e_2\otimes e_{u(\Phi(1),\pi_2(1))}\otimes e_{u(\Phi(2),\pi_2(2))}\\
&\hphantom{{}={}}+e_3\otimes e_{u(\Phi(1),\pi_3(1))}\otimes e_{u(\Phi(2),\pi_3(2))}\bigr)\otimes\bigl(e_1\otimes e_{u(\Phi(1),3-\pi_1(1))}\otimes e_{u(\Phi(2),3-\pi_1(2))}\\
&\hphantom{{}={}}+e_2\otimes e_{u(\Phi(1),3-\pi_2(1))}\otimes e_{u(\Phi(2),3-\pi_2(2))}+e_3\otimes e_{u(\Phi(1),3-\pi_3(1))}\otimes e_{u(\Phi(2),3-\pi_3(2))}\bigr)\\
&=\bigl(e_1\otimes e_{u(2\sigma+1,1)}\otimes e_{u(4\sigma+3,1)}+e_2\otimes e_{u(2\sigma+1,2)}\otimes e_{u(4\sigma+3,2)}+e_3\otimes e_{u(2\sigma+1,1)}\otimes e_{u(4\sigma+3,2)}\bigr)\\
&\hphantom{{}={}}\otimes\bigl(e_1\otimes e_{u(2\sigma+1,2)}\otimes e_{u(4\sigma+3,2)}+e_2\otimes e_{u(2\sigma+1,1)}\otimes e_{u(4\sigma+3,1)}+e_3\otimes e_{u(2\sigma+1,2)}\otimes e_{u(4\sigma+3,1)}\bigr)\\
&=(e_1\otimes e_{3}\otimes e_{5}+e_2\otimes e_{2}\otimes e_{4}+e_3\otimes e_{3}\otimes e_{4})\otimes(e_1\otimes e_{2}\otimes e_{4}+e_2\otimes e_{3}\otimes e_{5}+e_3\otimes e_{2}\otimes e_{5})\\
&=(e_1\otimes\varphi_{3,5}+e_2\otimes\varphi_{2,4}+e_3\otimes\varphi_{3,4})\otimes(e_1\otimes\psi_{2,4}+e_2\otimes\psi_{3,5}+e_3\otimes\psi_{2,5}),
\end{align*}
which can be viewed as a rank one matrix with $9$ blocks (see \Cref{fig:second_example_M}).
\end{example}
\begin{figure}[hbt!]
\[
\begin{bNiceArray}{ccc>{\columncolor{gray!30}}cccc|ccc>{\columncolor{gray!30}}cccc|ccc>{\columncolor{gray!30}}cccc}[first-row,code-for-first-row =\scriptscriptstyle  \rotate \color{blue},first-col,code-for-first-col=\scriptscriptstyle \color{blue},nullify-dots,margin,colortbl-like]
 & (1,1) & \Cdots &(3,4)  & (3,5) & (3,6)  & \Cdots & (\sigma,\sigma) & (1,1) & \Cdots &(2,3)  & (2,4) & (2,5)  & \Cdots & (\sigma,\sigma) & (1,1) & \Cdots &(3,3)  & (3,4) & (3,5)  & \Cdots & (\sigma,\sigma) \\
(1,1)  & 0 & \Cdots  & 0 & 0 & 0 & \Cdots & 0 & 0 & \Cdots & 0 & 0 & 0 & \Cdots & 0 & 0 & \Cdots  & 0 & 0 & 0 & \Cdots & 0 \\
\Vdots & \Vdots  & \Ddots & \Vdots & \Vdots & \Vdots & \Ddots & \Vdots & \Vdots  & \Ddots & \Vdots & \Vdots & \Vdots & \Ddots & \Vdots & \Vdots  & \Ddots & \Vdots & \Vdots & \Vdots & \Ddots & \Vdots \\
(2,3) & 0  & \Cdots & 0 & 0 & 0 & \Cdots & 0 & 0  & \Cdots & 0 & 0 & 0 & \Cdots & 0 & 0  & \Cdots & 0 & 0 & 0 & \Cdots & 0 \\
(2,4) & \rowcolor{gray!30} 0  & \Cdots & 0 & \color{darkred}1 & 0 & \Cdots & 0 & 0  & \Cdots & 0 & \color{darkred}1 & 0 & \Cdots & 0 & 0  & \Cdots & 0 & \color{darkred}1 & 0 & \Cdots & 0 \\
(2,5) & 0  & \Cdots & 0 & 0 & 0 & \Cdots & 0 & 0  & \Cdots & 0 & 0 & 0 & \Cdots & 0 & 0  & \Cdots & 0 & 0 & 0 & \Cdots & 0 \\
\Vdots & \Vdots  & \Ddots & \Vdots & \Vdots & \Vdots & \Ddots & \Vdots & \Vdots  & \Ddots & \Vdots & \Vdots & \Vdots & \Ddots & \Vdots & \Vdots  & \Ddots & \Vdots & \Vdots & \Vdots & \Ddots & \Vdots \\
(\sigma,\sigma) & 0 & \Cdots  & 0 & 0 & 0 & \Cdots & 0 & 0 & \Cdots & 0 & 0 & 0 & \Cdots & 0 & 0 & \Cdots  & 0 & 0 & 0 & \Cdots & 0 \\
\hline
(1,1)  & 0 & \Cdots  & 0 & 0 & 0 & \Cdots & 0 & 0 & \Cdots & 0 & 0 & 0 & \Cdots & 0 & 0 & \Cdots  & 0 & 0 & 0 & \Cdots & 0 \\
\Vdots & \Vdots  & \Ddots & \Vdots & \Vdots & \Vdots & \Ddots & \Vdots & \Vdots  & \Ddots & \Vdots & \Vdots & \Vdots & \Ddots & \Vdots & \Vdots  & \Ddots & \Vdots & \Vdots & \Vdots & \Ddots & \Vdots \\
(3,4) & 0  & \Cdots & 0 & 0 & 0 & \Cdots & 0 & 0  & \Cdots & 0 & 0 & 0 & \Cdots & 0 & 0  & \Cdots & 0 & 0 & 0 & \Cdots & 0 \\
(3,5) & \rowcolor{gray!30} 0  & \Cdots & 0 & \color{darkred}1 & 0 & \Cdots & 0 & 0  & \Cdots & 0 & \color{darkred}1 & 0 & \Cdots & 0 & 0  & \Cdots & 0 & \color{darkred}1 & 0 & \Cdots & 0 \\
(3,6) & 0  & \Cdots & 0 & 0 & 0 & \Cdots & 0 & 0  & \Cdots & 0 & 0 & 0 & \Cdots & 0 & 0  & \Cdots & 0 & 0 & 0 & \Cdots & 0 \\
\Vdots & \Vdots  & \Ddots & \Vdots & \Vdots & \Vdots & \Ddots & \Vdots & \Vdots  & \Ddots & \Vdots & \Vdots & \Vdots & \Ddots & \Vdots & \Vdots  & \Ddots & \Vdots & \Vdots & \Vdots & \Ddots & \Vdots \\
(\sigma,\sigma) & 0 & \Cdots  & 0 & 0 & 0 & \Cdots & 0 & 0 & \Cdots & 0 & 0 & 0 & \Cdots & 0 & 0 & \Cdots  & 0 & 0 & 0 & \Cdots & 0 \\
\hline
(1,1)  & 0 & \Cdots  & 0 & 0 & 0 & \Cdots & 0 & 0 & \Cdots & 0 & 0 & 0 & \Cdots & 0 & 0 & \Cdots  & 0 & 0 & 0 & \Cdots & 0 \\
\Vdots & \Vdots  & \Ddots & \Vdots & \Vdots & \Vdots & \Ddots & \Vdots & \Vdots  & \Ddots & \Vdots & \Vdots & \Vdots & \Ddots & \Vdots & \Vdots  & \Ddots & \Vdots & \Vdots & \Vdots & \Ddots & \Vdots \\
(2,4) & 0  & \Cdots & 0 & 0 & 0 & \Cdots & 0 & 0  & \Cdots & 0 & 0 & 0 & \Cdots & 0 & 0  & \Cdots & 0 & 0 & 0 & \Cdots & 0 \\
(2,5) & \rowcolor{gray!30} 0  & \Cdots & 0 & \color{darkred}1 & 0 & \Cdots & 0 & 0  & \Cdots & 0 & \color{darkred}1 & 0 & \Cdots & 0 & 0  & \Cdots & 0 & \color{darkred}1 & 0 & \Cdots & 0 \\
(2,6) & 0  & \Cdots & 0 & 0 & 0 & \Cdots & 0 & 0  & \Cdots & 0 & 0 & 0 & \Cdots & 0 & 0  & \Cdots & 0 & 0 & 0 & \Cdots & 0 \\
\Vdots & \Vdots  & \Ddots & \Vdots & \Vdots & \Vdots & \Ddots & \Vdots & \Vdots  & \Ddots & \Vdots & \Vdots & \Vdots & \Ddots & \Vdots & \Vdots  & \Ddots & \Vdots & \Vdots & \Vdots & \Ddots & \Vdots \\
(\sigma,\sigma) & 0 & \Cdots  & 0 & 0 & 0 & \Cdots & 0 & 0 & \Cdots & 0 & 0 & 0 & \Cdots & 0 & 0 & \Cdots  & 0 & 0 & 0 & \Cdots & 0 \\
\end{bNiceArray}
\]
\caption{Matrix $M^{\Phi}$ obtained in \Cref{exam:second_example_M}. The coordinates in blue represent the indices of the elements of the canonical basis of $(\bbC^{\sigma})^{\otimes \theta}=(\bbC^{\numprint{14406}})^{\otimes 2}$.
\label{fig:second_example_M}}
\end{figure}
By looking at \Cref{exam:second_example_M}, it is evident that the set of all the functions $\Phi$ covers all the possible entries of the diagonal blocks. By this fact, we obtain the following proposition, which describes how  the sum of all the matrices $M^{\Phi}$ is structured (see \Cref{fig:example_sum_matrices_MPhi}).
\begin{proposition}
\label{prop:matrices_MPhi_all_ones}
The diagonal blocks of the matrix 
\[\sum_{M^\Phi\in\calF}M^\Phi\] are matrices of all ones, that is,
\[
\sum_{M^\Phi\in\calF}M_{ii}^\Phi=\sum_{\alpha_1,\dots,\alpha_{\theta},\beta_1,\dots,\beta_{\theta}=1}^{\sigma}\varphi_{\alpha_1,\dots,\alpha_{\theta}}\otimes\psi_{\beta_1,\dots,\beta_{\theta}}, \qquad  i=1,\dots,r. 
\]
\end{proposition}
\begin{proof}
For any $i=1,\dots,r$, let 
\begin{equation}\label{formula:rank_one_diagonal_matrix}
a_i\otimes \varphi_{\alpha_1,\dots,\alpha_{\theta}}\otimes b_i\otimes \psi_{\beta_1,\dots,\beta_{\theta}}\in \langle a_i\rangle\otimes (\bbC^\sigma)^{\otimes \theta}\otimes \langle b_i\rangle\otimes(\bbC^{\sigma})^{\otimes \theta},
\end{equation}
and let $\pi_i\colon\{1,\dots,\theta\}\to\{1,2\}$
be the chosen function for the coefficient $i$, in the sense of \Cref{def:matrix_M_Phi}. We consider the function $\Phi$ defined by
\[
\Phi(\gamma)\coloneqq\begin{cases}
(\beta_\gamma-1)\sigma+(\alpha_{\gamma}-1),\quad &\text{if $\pi_i(\gamma)=1$},\\
(\alpha_\gamma-1)\sigma+(\beta_{\gamma}-1),\quad &\text{if $\pi_i(\gamma)=2$},
\end{cases}
\]
for every $\gamma=1,\dots,\theta$. 
In particular, it follows that
\begin{equation}\label{formula:equation_unit_component}
M_{ii}^{\Phi}(\varphi_{\alpha_1,\dots,\alpha_{\theta}}^*\otimes\psi_{\beta_1,\dots,\beta_{\theta}}^*)=M_{ii}^{\Phi}\bigl((e_{\alpha_1}^*\otimes\cdots \otimes e_{\alpha_\theta}^*)\otimes (e_{\beta_1}^*\otimes\cdots \otimes e_{\beta_\theta}^*)\bigr)=1.
\end{equation}
Hence, every rank one tensor as in~\eqref{formula:rank_one_diagonal_matrix} appears in at least one matrix $M^\Phi\in\calF$. 
Now we prove that the corresponding function $\Phi$ of the matrix $M^{\Phi}$ is unique. For this, consider two functions
$\Phi,\Phi'\colon \{1,\dots,\theta\}\to\{0,\dots,\sigma^2-1\}$
 such that $\Phi\neq \Phi'$. Then there exists a number $\gamma_{0}\in\{1,\dots,\theta\}$ such that
$
\Phi(\gamma_{0})\neq\Phi'(\gamma_{0}).
$
Therefore, we can write
\[
\Phi(\gamma_{0})=(\alpha_{\gamma_{0}}-1)\sigma+(\beta_{\gamma_{0}}-1),\qquad\Phi'(\gamma_{0})=(\alpha_{\gamma_{0}}'-1)\sigma+(\beta_{\gamma_{0}}'-1),
\]
for some $\alpha_{\gamma_{0}},\alpha_{\gamma_{0}}',\beta_{\gamma_{0}},\beta_{\gamma_{0}}'\in\{1,\dots,\sigma\}$ such that either $\alpha_{\gamma_{0}}\neq \alpha_{\gamma_{0}}'$ or $\beta_{\gamma_{0}}\neq \beta'_{\gamma_{0}}$. Since 
\[
\varphi_{\alpha_1,\dots,\alpha_{\theta}}\otimes\psi_{\beta_1,\dots,\beta_{\theta}}\neq \varphi_{\alpha_1',\dots,\alpha_{\theta}'}\otimes\psi_{\beta_1',\dots,\beta_{\theta}'},
\]
if equality \eqref{formula:equation_unit_component} holds,
we must have
\[
M_{ii}^{\Phi'}(\varphi_{\alpha_1,\dots,\alpha_{\theta}}^*\otimes \psi_{\beta_1,\dots,\beta_{\theta}}^*)=0.\qedhere\]
\end{proof}

\begin{corollary}\label{prop:Sigma_clone_W_calM}
Let $w^{\Sigma_{\theta}}$ be the $\Sigma_{\theta}$-clone of $w$. Then $w^{\Sigma_{\theta}}\in\calM$.
\end{corollary}
\begin{proof}
This is a direct consequence of \Cref{prop:matrices_MPhi_all_ones}, since by the definition of $\calM$ we can remove the rank one tensors in non-diagonal blocks, which are elements of $\calU\subset\calM$. 
\end{proof}

Now we present the following technical lemma, which is essential to obtain \Cref{lemma1}.
\begin{lemma}
\label{lem:permutations_rows_columns}
For any $M\in\calM$ and $i\neq j$, any non-diagonal block $M_{ij}=M(a_i^*\otimes  b_j^*)$ is a $\sigma\times\sigma$ block-diagonal matrix with blocks of size $\sigma^{\theta-1}\times\sigma^{\theta-1}$. That is,
\[
M_{ij}\in {\big\langle e_{\delta}\otimes e_{\delta}\otimes (\bbC^{\sigma})^{\otimes(\theta-1)}\otimes (\bbC^{\sigma})^{\otimes(\theta-1)}\big\rangle}_{\delta=1,\dots,\sigma}.
\]
\end{lemma}
\begin{proof}
We can consider linear combinations of the matrices in $\calF$ only because the union of their supports contains the support of any other matrix from $\calM$. 
We have
\[
M_{ij}^{\Phi}=e_{u(\Phi(1),\pi_i(1))}\otimes\cdots\otimes e_{u(\Phi(\theta),\pi_i(\theta))}\otimes e_{u(\Phi(1),3-\pi_j(1))}\otimes\cdots\otimes e_{u(\Phi(\theta),3-\pi_j(\theta))}.
\]
Since $i\neq j$, the selected functions $\pi_i$ and $\pi_j$ are distinct. Thus, there exists an element $\tau\in\{1,\dots,\theta\}$ such that $\pi_i(\tau)=3-\pi_j(\tau)$. 
In particular, we have
\[
u_{(\Phi(\tau),\pi_i(\tau))}=u_{(\Phi(\tau),3-\pi_j(\tau))}.
\]
Therefore, just permuting the order of the factors of the tensor product and putting the $\tau$-th and the $(\tau+\theta)$-th elements at the beginning, we may write
\begin{align*}
M_{ij}^{\Phi}&=e_{u(\Phi(\tau),\pi_i(\tau))}\otimes e_{u(\Phi(\tau),3-\pi_j(\tau))}\otimes \bigotimes_{\substack{1\leq \gamma\leq\theta\\ \gamma\neq \tau}} e_{u(\Phi(\gamma),\pi_i(\gamma))}\otimes\bigotimes_{\substack{1\leq \delta\leq\theta\\ \delta\neq \tau}} e_{u(\Phi(\delta),3-\pi_j(\delta))}\\
&=e_{u(\Phi(\tau),\pi_i(\tau))}\otimes e_{u(\Phi(\tau),\pi_i(\tau))}\otimes \bigotimes_{\substack{1\leq \gamma\leq\theta\\ \gamma\neq \tau}} e_{u(\Phi(\gamma),\pi_i(\gamma))}\otimes\bigotimes_{\substack{1\leq \delta\leq\theta\\ \delta\neq \tau}} e_{u(\Phi(\delta),3-\pi_j(\delta))}.
\end{align*}
In particular, setting $\delta_{\tau}\coloneqq u(\Phi(\tau),\pi_i(\tau))$, we have
\[
M_{ij}^{\Phi}\in e_{\delta_{\tau}}\otimes e_{\delta_{\tau}}\otimes (\bbC^{\sigma})^{\otimes(\theta-1)}\otimes (\bbC^{\sigma})^{\otimes(\theta-1)}.
\]
That is, for any given functions $\Phi_1,\dots,\Phi_m$ and any $\alpha_1,\dots,\alpha_m\in\bbC$, we have
\[
\sum_{k=1}^m\alpha_k M_{ij}^{\Phi_k}\in \bigoplus_{\delta=1}^{\sigma}\bigl( e_{\delta}\otimes e_{\delta}\otimes (\bbC^{\sigma})^{\otimes(\theta-1)}\otimes (\bbC^{\sigma})^{\otimes(\theta-1)}\bigr).\qedhere
\] 
\end{proof}
The proof of the following lemma is trivial. 
\begin{lemma}
\label{rem:partion_sets_S_D}
Let $S_1,\dots, S_r$ be disjoint sets of cardinality $\omega$, and let their union be partitioned into $\omega$ sets $D_1,\dots,D_\omega$. Consider $Z_{\alpha}\subset S_{\alpha}$  such that $|Z_{\alpha}|<\omega/r$ for every $\alpha=1,\dots,r$. 
Then, there exists a~number $\delta\in\mathbb{Z}_{>0}$,  such that $1\leq\delta\leq\omega$ and
\[
D_{\!\delta}\cap\biggl(\bigcup_{\alpha=1}^r Z_{\alpha}\biggr)=\varnothing.
\] 
\end{lemma}
Recalling by the initial hypothesis of  \eqref{formula:initial_values} that 
$\rho=2\bfa\bfb\bfc+1$ and $r=\rk(w)=\sigma/2\rho^2$,
we also have the following technical statement.
\begin{lemma}
\label{lem:M_0_rank_less_rho}
Let $\calM_0\subset \calM$ be a subset of cardinality $\bfc$. Then, one of the two following conditions holds:
\begin{enumerate}[label=(\arabic*), left= 0pt, widest=*,nosep]
\item $\calM_0$ contains a matrix of rank at least $\rho$;
\item there are $2r$ tensors $\varphi_1,\dots,\varphi_r,\psi_1,\dots,\psi_r\in\Sigma_{\theta}$
such that, for every $M\in\calM_0$,
\[
M\big|_{{\langle a_i^*\otimes\varphi_i^*\rangle}_{i=1}^r\otimes {\langle b_j^*\otimes\psi_j^*\rangle}_{j=1}^r}=\lambda\sum_{i=1}^ra_i\otimes\varphi_i\otimes b_i\otimes \psi_i
\]
for some $\lambda\in\bbC$. That is, it corresponds to a scalar $r\times r$ matrix.
\end{enumerate}
\end{lemma}
\begin{proof}
Let us assume that condition (1) is false, that is, every matrix in the set $\calM_0$ has rank at most $\rho-1=2\bfa\bfb\bfc$. For any $M\in\calM_0$, by \Cref{lem:permutations_rows_columns} we can write 
\[
M_{ij}=\sum_{\delta=1}^{\sigma} M_{ij\delta}\in\bigoplus_{\delta=1}^{\sigma}\bigl(e_{\delta}\otimes (\bbC^{\sigma})^{\otimes(\theta-1)}\otimes e_{\delta}\otimes (\bbC^{\sigma})^{\otimes(\theta-1)}\bigr),
\]
where 
\[
M_{ij\delta}\in e_{\delta}\otimes (\bbC^{\sigma})^{\otimes(\theta-1)}\otimes e_{\delta}\otimes (\bbC^{\sigma})^{\otimes(\theta-1)}
\]
for every $i,j,\delta=1,\dots,\sigma$ with $i\neq j$. 
By the hypothesis, for at most $(\rho-1)\sigma^{\theta-1}$ elements $\varphi\in\Sigma_{\theta}$,
\begin{equation}\label{formula:elements_block_satisfying_non_zero}
M_{ij}(\varphi^*)\neq 0.
\end{equation}
Indeed, otherwise there would be at least $\rho$ values $1\leq \delta_1,\dots,\delta_{\rho}\leq \sigma$ such that $M_{ij\delta_k}\neq 0$
for every $k=1,\dots,\rho$, implying that
$\rk(M_{ij})\geq\rho$.
For any $\alpha\in\{1,\dots,r\}$, we define the set 
\[
Z_{\alpha}\coloneqq\Set{\varphi\in\Sigma_{\theta}|\exists M\in\calM_0,\,\exists j\neq\alpha:M_{\alpha j}(\varphi^*)\neq 0}.
\]
In particular, since the number of elements satisfying formula \eqref{formula:elements_block_satisfying_non_zero} is at most $(\rho-1)\sigma^{\theta-1}$ for any $j\neq\alpha$
and $|\calM_0|=\bfc$, we have
\begin{equation}
\label{formula:upper_bound_rows}
|Z_\alpha|\leq(r-1)(\rho-1)\sigma^{\theta-1}\bfc<r\rho\sigma^{\theta-1}\bfc<\rho^2\sigma^{\theta-1}=\frac{\rho^2}{2\rho^2r}\sigma^{\theta}=\frac{\sigma^{\theta}}{2r}.
\end{equation}
Now, for every $\alpha=1,\dots,r$ we define the set
\[
S_\alpha\coloneqq\Set{a_{\alpha}\otimes \varphi\otimes b_{\alpha}\otimes \psi|\varphi, \psi\in\Sigma_{\theta}}
\]
of cardinality $\omega\coloneqq|S_{\alpha}|=\sigma^{2\theta}$.
By \Cref{prop:matrices_MPhi_all_ones}, the family
${\{D_{\Phi}\}}_{\Phi\in\calF}$, where
$D_{\Phi}\coloneqq\{M_{ii}^{\Phi}\}_{i=1,\dots,r}$ for any $\Phi\in\calF$,
can be interpreted, by considering their summands, as a partition of the set
\[
\calS\coloneqq\bigcup_{\alpha=1}^rS_\alpha.
\]
The expression of the linear combination of $\Phi$ has exactly one element of $S_{\alpha}$ for every $\alpha=1,\dots,r$. Moreover, since, for any $\alpha=1,\dots,r$,
\[
|Z_{\alpha}|< \frac{\sigma^{\theta}}{2r}< \frac{\sigma^{2\theta}}{r}=\frac{\omega}{r} 
\]
 and $|\calF|=\omega=\sigma^{2\theta}$,
\Cref{rem:partion_sets_S_D} tells us that there exists a function $\hat\Phi\in {\{0,\dots,\sigma^2-1\}}^{\{1,\dots,\theta\}}$
such that
\[
D_{\hat\Phi}\cap\biggl(\bigcup_{\alpha=1}^r Z_{\alpha}\biggr)=\varnothing.
\] 
That is, if
$M_{\alpha\alpha}^{\hat\Phi}=\varphi_{\alpha}\otimes\psi_{\alpha}$,
with $\varphi_{\alpha},\psi_{\alpha}\in\Sigma_{\theta}$ for every $\alpha=1,\dots,r$, then
\begin{equation}\label{formula:Mij}
M_{ij}\bigl(\varphi_{\alpha}^*\bigr)=0, \quad  M\in\calM_0, \quad  i\neq j.
\end{equation}
Now, for any $M\in\calM_0$, let us consider the restriction
\[
M\big|_{{\langle a_i^*\otimes\varphi_i^*\rangle}_{i=1}^r\otimes {\langle b_j^*\otimes\psi_j^*\rangle}_{j=1}^r}=\sum_{i,j=1}^r\lambda_{ij}a_i\otimes \varphi_i\otimes b_j\otimes \psi_j, \quad \lambda_{ij}\in\bbC.
\]
By relation \eqref{formula:Mij} we have that $\lambda_{ij}=0$ whenever $i\neq j$. Therefore, we can write
\begin{equation}
\label{formula:linear_combination_Mi}
M\big|_{{\langle a_i^*\otimes\varphi_i^*\rangle}_{i=1}^r\otimes {\langle b_j^*\otimes\psi_j^*\rangle}_{j=1}^r}=\sum_{i=1}^r\lambda_{i}a_i\otimes \varphi_i\otimes b_i\otimes \psi_i, \quad   \lambda_i\in\bbC.
\end{equation}
In particular, since for every $N\in\calU$ (defined in \eqref{formula:space_U}) and for every $i=1,\dots,r$ we have $N_{ii}=0$, the restriction \eqref{formula:linear_combination_Mi} is spanned only by restrictions of elements of $\calF$. But by \Cref{prop:matrices_MPhi_all_ones} we conclude that it is spanned only by the restriction of $M^{\hat\Phi}$, that is,
\[
M^{\Phi}\big|_{{\langle a_i^*\otimes\varphi_i^*\rangle}_{i=1}^r\otimes {\langle b_j^*\otimes\psi_j^*\rangle}_{j=1}^r}=
\lambda\sum_{i=1}^ra_i\otimes \varphi_i\otimes b_i\otimes \psi_i
\]
for some $\lambda\in\bbC$.
\end{proof}

\subsection{Minimum rank of modifications of clones}
The following proposition is a direct consequence of \Cref{lem:M_0_rank_less_rho} and involves the use of modifications of tensors we defined in \Cref{def: modofocation}.
\begin{proposition}
\label{cor_double_inequalities}
Let $T^{\Sigma_{\theta}}$ be the $\Sigma_{\theta}$-clone of $T$. Then
\[
\min\rk \bigl(T^{\Sigma_{\theta}}\bmod(0,0,\calM)\bigr)=\min\rk \bigl(T\bmod(0,0,\langle w\rangle)\bigr).
\]
\end{proposition}
\begin{proof}
The  inequality 
\begin{equation}\label{formula:inequality_rank}
\min\rk \bigl(T^{\Sigma_{\theta}}\bmod(0,0,\calM)\bigr)\leq \min\rk \bigl(T\bmod(0,0,\langle w\rangle)\bigr)
\end{equation}
is a direct consequence of \Cref{prop:Sigma_clone_W_calM}. 
To obtain the reverse inequality, let us first observe that, since $\rho=2\bfa\bfb\bfc+1$, by trivial inequality of \Cref{prop:trivial_lower_upper} we have that
\[
\rk(T^{\Sigma_{\theta}})=\rk(T)\leq\bfa\bfb=\frac{\rho-1}{2\bfc}\leq \frac{\rho-1}{2}=: s_{\rho}.
\]
By \Cref{prop:T(A)_contained_BC}, this means that the space 
\[
T^{\Sigma_{\theta}}\bigl(C^*\otimes{(\bbC^{\sigma})^{\otimes\theta}}^*\bigr)\subset A\otimes (\bbC^{\sigma})^{\otimes\theta}\otimes B\otimes (\bbC^{\sigma})^{\otimes\theta}
\]
is generated by at most $s_\rho$ rank one  elements from $A\otimes (\bbC^{\sigma})^{\otimes\theta}\otimes B\otimes (\bbC^{\sigma})^{\otimes\theta}$.
Now, recalling \Cref{defn:sigma_clone} and given
\begin{equation}\label{formula:defn:sigma_clone_T}
S_{\Sigma_{\theta}}\coloneqq\sum_{v\in\Sigma_{\theta}}v, \qquad T=\sum_{k=1}^{\bfc}f_k\otimes c_k,
\end{equation}
we have
\begin{equation}
\label{formula:clone_T_theta}
T^{\Sigma_{\theta}}=\sum_{k=1}^{\bfc}(f_{k}^{\Sigma_{\theta}})\otimes (c_{k}\otimes S_{\Sigma_{\theta}}).
\end{equation}
Let us  
consider any modification
\[
T_{\mathrm{mod}}^{\Sigma_{\theta}}\coloneqq T^{\Sigma_{\theta}}+\sum_{j=1}^{\bfc}M_{j}\otimes (c_{j}\otimes S_{\Sigma_{\theta}})=\sum_{j=1}^{\bfc}\bigl((f_j^{\Sigma_{\theta}})+M_{j}\bigr)\otimes (c_{j}\otimes S_{\Sigma_{\theta}}),
\]
where $M_j\in\langle\calM\rangle\subset A\otimes (\bbC^{\sigma})^{\otimes\theta}\otimes B\otimes (\bbC^{\sigma})^{\otimes\theta}$ for every $j=1,\dots,\bfc$.
First, we consider the case when there is an index $\delta$ such that $\rk (M_{\delta})\geq\rho$. 
Considering again the definition of $T$ in  formula 
\eqref{formula:defn:sigma_clone_T} and inequality \eqref{formula:inequality_rank}, by \Cref{prop:T(A)_contained_BC} there are $s_\rho$ elements $$g_1,\dots,g_{s_{\rho}}\in A\otimes (\bbC^{\sigma})^{\otimes\theta}\otimes B\otimes (\bbC^{\sigma})^{\otimes\theta}$$ of  rank one 
such that
$f_j^{\Sigma_{\theta}}\in\langle g_1,\dots,g_{s_{\rho}}\rangle$. If we had $\rk(T^{\Sigma_{\theta}}_{\textrm{mod}})\leq \rk(T^{\Sigma_{\theta}})$,
then again by \Cref{prop:T(A)_contained_BC}, we would have other $s_\rho$ rank one  elements $h_1,\dots,h_{s_{\rho}}\in A\otimes (\bbC^{\sigma})^{\otimes\theta}\otimes B\otimes (\bbC^{\sigma})^{\otimes\theta}$ such that
\[
f_j^{\Sigma_{\theta}}+M_j\in\langle h_1,\dots,h_{s_\rho}\rangle.
\]
This would imply that
$M_j\in \langle g_1,\dots,g_{s_{\rho}},h_1,\dots,h_{s_{\rho}}\rangle$,
which is not possible since $\rk(M_j)\geq\rho>\rho-1=2s_{\rho}$.
Therefore, in the case when $\rk(M_{\delta})\geq\rho$, we must have 
\[\rk(T^{\Sigma_{\theta}}_{\textrm{mod}})> \rk(T^{\Sigma_{\theta}}).\]

Now we may assume that $\rk(M_j)<\rho$ for every $j=1,\dots,\bfc$.
Then, by \Cref{lem:M_0_rank_less_rho} there are $2r$ tensors $\varphi_1,\dots,\varphi_r,\psi_1,\dots,\psi_r\in\Sigma_{\theta}$ such that, for every $k=1,\dots,\bfc$,
\[
M_k\big|_{{\langle a_i^*\otimes\varphi_i^*\rangle}_{i=1}^r\otimes {\langle b_j^*\otimes\psi_j^*\rangle}_{j=1}^r}=\lambda_k\sum_{\alpha=1}^ra_{\alpha}\otimes\varphi_{\alpha}\otimes b_{\alpha}\otimes \psi_{\alpha}.
\]
For every $i\in\{r+1,\dots,\bfa\}, j\in \{r+1,\dots,\bfb\}$, and 
 $k\in \{1,\dots,\bfc\}$ we  set $\varphi_i,\psi_j,\zeta_k\coloneqq e_1^{\otimes\theta}$.
We have
\begin{align*}
&T_{\mathrm{mod}}^{\Sigma_{\theta}}\big|_{{\langle a_i^*\otimes\varphi_i^*\rangle}_{i=1}^{\bfa}\otimes {\langle b_j^*\otimes\psi_j^*\rangle}_{j=1}^{\bfb}\otimes {\langle c_k^*\otimes\zeta_k^*\rangle}_{k=1}^{\bfc}}=
\sum_{k=1}^{\bfc}\Bigl(\bigl(f_k^{\Sigma_{\theta}}+M_k\bigr)\big|_{{\langle a_i^*\otimes\varphi_i^*\rangle}_{i=1}^r\otimes {\langle b_j^*\otimes\psi_j^*\rangle}_{j=1}^r}\Bigr)\otimes c_k\otimes \zeta_k\\
&\qquad=
\sum_{k=1}^{\bfc}\biggl(\sum_{i=1}^{\bfa}\sum_{j=1}^{\bfb}f_k(a_i^*,b_j^*)(a_i\otimes\varphi_i)\otimes (b_j\otimes\psi_j)+\lambda_k\sum_{\alpha=1}^r(a_{\alpha}\otimes\varphi_{\alpha})\otimes (b_{\alpha}\otimes \psi_{\alpha})\biggr)\otimes (c_k\otimes \zeta_k)\\
&\qquad=
\sum_{i=1}^{\bfa}\sum_{j=1}^{\bfb}\sum_{k=1}^{\bfc}f_k(a_i^*,b_j^*)(a_i\otimes\varphi_i)\otimes (b_j\otimes\psi_j)\otimes (c_k\otimes\zeta_k)\\
&\qquad\hphantom{{}={}}+\sum_{\alpha=1}^r(a_{\alpha}\otimes\varphi_{\alpha})\otimes (b_{\alpha}\otimes \psi_{\alpha})\otimes \sum_{k=1}^{\bfc}\lambda_k(c_k\otimes \zeta_k).
\end{align*}
 By identifying $a_i\otimes\varphi_i\simeq a_i$, $b_j\otimes\psi_j\simeq b_j$, and $c_k\otimes\zeta_k\simeq c_k$,
for every 
$i=1,\dots,\bfa$, $j= 1,\dots,\bfb$, and $k= 1,\dots,\bfc$, respectively,
we get
\begin{align*}
T_{\mathrm{mod}}^{\Sigma_{\theta}}\big|_{{\langle a_i^*\otimes\varphi_i^*\rangle}_{i=1}^{\bfa}\otimes {\langle b_j^*\otimes\psi_j^*\rangle}_{j=1}^{\bfb}\otimes {\langle c_k^*\otimes\zeta_k^*\rangle}_{k=1}^{\bfc}}&\simeq \sum_{i=1}^{\bfa}\sum_{j=1}^{\bfb}\sum_{k=1}^{\bfc}f_k(a_i^*,b_j^*)a_i\otimes b_j\otimes c_k+\sum_{\alpha=1}^ra_{\alpha}\otimes b_{\alpha}\otimes \sum_{k=1}^{\bfc}\lambda_kc_k\\
&=\sum_{k=1}^{\bfc}f_k\otimes c_k+w\otimes \sum_{k=1}^{\bfc}\lambda_k c_k=T+\sum_{k=1}^{\bfc}\lambda_k w\otimes c_k,
\end{align*}
which is a modification of $T$ by $\langle w\rangle$. This implies that
\[
\rk(T_{\mathrm{mod}}^{\Sigma_{\theta}})\geq \rk\Bigl(T_{\mathrm{mod}}^{\Sigma_{\theta}}\big|_{{\langle a_i^*\otimes\varphi_i^*\rangle}_{i=1}^{\bfa}\otimes {\langle b_j^*\otimes\psi_j^*\rangle}_{j=1}^{\bfb}\otimes {\langle c_k^*\otimes\zeta_k^*\rangle}_{k=1}^{\bfc}}\Bigr)\geq \min\rk \bigl(T\bmod(0,0,\langle w\rangle)\bigr),
\]
and, in particular,
\[
\min\rk \bigl(T^{\Sigma_{\theta}}\bmod(0,0,\calM)\bigr)\geq \min\rk \bigl(T\bmod(0,0,\langle w\rangle)\bigr).\qedhere
\]
\end{proof}
Now, we illustrate how \Cref{cor_double_inequalities} can be used to prove the general case of \Cref{prop:induction_case1}. 
\begin{proof}[Proof of \Cref{prop:induction_case1}]
We set $V\coloneqq(\bbC^{\sigma})^{\otimes\theta}$, with the basis $\Sigma\coloneqq\Sigma_{\theta}$, and we consider the space
$\calM$ defined as in \eqref{formula:space_M}. Let us consider the $\Sigma$-clone $T^{\Sigma}$ and the augmented tensor 
\begin{equation}
\label{formula:recursive_T}
T^{\Sigma}+T_{\calM}\in \bigl(A\otimes(\bbC^{\sigma})^{\otimes\theta}\bigr)\otimes \bigl(B\otimes(\bbC^{\sigma})^{\otimes\theta}\bigr)\otimes \Bigl(\bigl(C\otimes(\bbC^{\sigma})^{\otimes\theta}\bigr)\oplus\bbC^{\bfm}\Bigr).
\end{equation}
Since the set $\calM$ is composed only by rank one elements, point (1) of \Cref{lemma1} holds. Point (2), instead, is given by \Cref{prop:Sigma_clone_W_calM}. 
Moreover, setting $\bfm\coloneqq\dim\calM$, we get
\[
\rk(T^{\Sigma}+T_{\calM})=\min\rk\bigl(T^{\Sigma}\bmod(0,0,\calM)\bigr)+\bfm=\min\rk \bigl(T\bmod(0,0,\langle w\rangle)\bigr)+\bfm,
\]
where the first equality follows from \Cref{cor2} and the second one from \Cref{cor_double_inequalities}, which proves point (3) and concludes the proof.
\end{proof}
We are finally ready to present the proof of \Cref{lemma1}.
\begin{proof}[Proof of \Cref{lemma1}]
We proceed by induction on the sum of the dimensions of the spaces $W_A$, $W_B$, $W_C$. The first case is when two spaces are zero vector spaces and the third space is generated by one element, 
which is exactly the case of \Cref{prop:induction_case1}, and hence proved.
Now, we assume that
\[
W_A= W_A',\quad W_B= W_B',\quad W_C=  W_C'\oplus \langle w\rangle,
\]
where $w\not\in W_C'$ and, by the induction hypothesis,
the spaces $W_A'$, $W_B'$, and $W_C'$ satisfy the conditions of \Cref{lemma1}. 
In particular, there exists a vector space $V'$, with basis $\Sigma'$, and three spaces $\calM_A'$, $\calM_B'$, $\calM_C'$ spanned by elements of rank one, such that 
\begin{equation}\label{formula:inductive_hypothesis_2}
W_A'^{\Sigma'}\subset\calM_A',\quad W_B'^{\Sigma'}\subset\calM_B',\quad W_C'^{\Sigma'}\subset\calM_C',
\end{equation}
and, for any $T\in A\otimes B\otimes C$,
\[\rk(T^{\Sigma'} + T_{\mathcal{M}_A'} +  T_{\mathcal{M}_B'} + T_{\mathcal{M}_C'}) =\min \rk\bigl(T  \bmod (W_A',W_B',W_C')\bigr)+{\bf m}_A'+{\bf m}_B'+{\bf m}_C'.\] 
In particular, the equality 
\begin{equation}
\label{formula:inductive_hypothesis}
\min \rk\bigl(T^{\Sigma'}  \bmod (\calM_A',\calM_B',\calM_C')\bigr)=\min \rk\bigl(T  \bmod (W_A',W_B',W_C')\bigr)
\end{equation}
holds.
For 
simplicity, we denote
\[
A'\coloneqq A\otimes V',\quad B'\coloneqq B\otimes V',\quad C'\coloneqq C\otimes V',
\]
so that $T^{\Sigma'}\in A'\otimes B'\otimes C'$.
We apply \Cref{prop:induction_case1} to the element \[w^{\Sigma'}\in A'\otimes B'\subset (A'\oplus\bbC^{\bfm_A'})\otimes(B'\oplus \bbC^{\bfm_B'})\]
obtaining a vector space $V''$, with basis $\Sigma''$, and the vector space $\calM_C''$ satisfying the three conditions of 
\Cref{prop:induction_case1}.
We now show that 
$V\coloneqq V'\otimes V''$, with basis $\Sigma\coloneqq\Sigma'\otimes\Sigma''$, and
\[
\calM_A\coloneqq\calM_A'^{\Sigma''},\quad \calM_B\coloneqq\calM_B'^{\Sigma''},\quad
\calM_C\coloneqq\langle\calM_C'^{\Sigma''},\calM_C''\rangle
\]
satisfy the three conditions of \Cref{lemma1}.
 Since by \Cref{prop:induction_case1}  $\calM_C''$ is generated by rank one elements and, by the induction hypothesis, $\calM_A'$, $\calM_B'$, $\calM_C'$ 
are generated by rank one  elements as well, then so are $\calM_A$, $\calM_B$, $\calM_C$, which ensures condition (1).
We observe that the $\Sigma''$-clone of a $\Sigma'$-clone of a tensor $T$ is simply the $(\Sigma'\otimes \Sigma'')$-clone of $T$. 
Therefore, we have 
\begin{equation*}
{(w^{\Sigma'})}^{\Sigma''}=w^{\Sigma'\otimes \Sigma''}=w^{\Sigma}\in\calM_C'',
\end{equation*}
and, by inclusions \eqref{formula:inductive_hypothesis_2}, also
\begin{gather*}
W_A^{\Sigma}=W_A^{\Sigma'\otimes\Sigma''}={(W_A^{\Sigma'})}^{\Sigma''}\subset\calM_A'^{\Sigma''},\qquad
W_B^{\Sigma}=W_B^{\Sigma'\otimes\Sigma''}={(W_B^{\Sigma'})}^{\Sigma''}\subset\calM_B'^{\Sigma''},\\
(W_C')^{\Sigma}=(W_C')^{\Sigma'\otimes\Sigma''}={\bigl((W_C')^{\Sigma'}\bigr)}^{\Sigma''}\subset(\calM_C')^{\Sigma''},
\end{gather*}
which give condition (2). Finally, to obtain condition (3),
we first note that by \Cref{cor2} we have
\begin{align}\label{formula:middle_situation}
\nonumber&\rk\bigl(T^{\Sigma}+T_{\calM_A}+T_{\calM_B}+T_{\calM_C}\bigr)=\rk\Bigl(T^{\Sigma'\otimes\Sigma''}+T_{\calM_A'^{\Sigma''}}+T_{\calM_B'^{\Sigma''}}+T_{\langle\calM_C'^{\Sigma''},\,\calM_C''\rangle}\Bigr)\\
&\qquad=\min\rk\Bigl({(T^{\Sigma'})}^{\Sigma''}\bmod\bigl(\calM_A'^{\Sigma''},\calM_B'^{\Sigma''},\langle\calM_C'^{\Sigma''},\calM_C''\rangle\bigr)\Bigr)+\bfm_A+\bfm_B+\bfm_C.
\end{align}
Since any modification $T^{\Sigma'}_{\mathrm{mod}}$ of $T^{\Sigma'}$ by $\calM_A'$, $\calM_B'$, and $\calM_C'$ is still a tensor in the space $A'\otimes B'\otimes C'$, we can apply \Cref{cor_double_inequalities} and get the equality 
\begin{align*}
&\min \rk\bigl(T^{\Sigma}  \bmod (\calM_A'^{\Sigma''},\calM_B'^{\Sigma''},\langle\calM_C'^{\Sigma''},\calM_C''\rangle)\bigr)\\[1ex]
&\qquad=\min\Set{\min\rk\bigl((T^{\Sigma'}_{\mathrm{mod}})^{\Sigma''}\bmod(0,0,\calM_{C}'')\bigr)|T^{\Sigma'}_{\mathrm{mod}}\in T^{\Sigma'}\bmod(\calM_A',\calM_B',\calM_C')}\\[1ex]
&\qquad=\min\Set{\min\rk\bigl(T^{\Sigma'}_{\mathrm{mod}}\bmod(0,0,\langle w^{\Sigma'}\rangle)\bigr)|T^{\Sigma'}_{\mathrm{mod}}\in T^{\Sigma'}\bmod(\calM_A',\calM_B',\calM_C')}\\[1.5ex]
&\qquad=\min\rk\bigl(T^{\Sigma'}\bmod(\calM_A',\calM_B',\langle\calM_C',w^{\Sigma'}\rangle)\bigr).
\end{align*}
Proceeding analogously, we observe that any modification $T_{\mathrm{mod}}$ of $T$ by $\langle w\rangle$ is still a tensor in the space $A\otimes B\otimes C$. Therefore, we can use formula \eqref{formula:inductive_hypothesis} given by the induction hypothesis to get 
\begin{align*}
&\min\rk\bigl(T^{\Sigma'}\bmod(\calM_A',\calM_B',\langle\calM_C',w^{\Sigma'}\rangle)\bigr)\\[1ex]
&\qquad=\min\Set{\min\rk\bigl(T^{\Sigma'}_{\mathrm{mod}}\bmod(\calM_A',\calM_B',\calM_{C}')\bigr)|T_{\mathrm{mod}}\in T\bmod(0,0,\langle w\rangle)}\\[1ex]
&\qquad=\min\Set{\min\rk\bigl(T_{\mathrm{mod}}\bmod(W_A',W_B',W_{C}')\bigr)|T_{\mathrm{mod}}\in T\bmod(0,0,\langle w\rangle)}\\[2ex]
&\qquad=\min\rk\bigl(T\bmod(W_A',W_B',W_{C}'\oplus\langle w\rangle)\bigr)=\min\rk\bigl(T\bmod(W_A,W_B,W_C)\bigr).
\end{align*}
Substituting this value in formula \eqref{formula:middle_situation}, we have the equality 
\begin{equation}
\label{formula:l1}
\rk(T^{\Sigma} + T_{\mathcal{M}_A} +  T_{\mathcal{M}_B} + T_{\mathcal{M}_C})=\min \rk\bigl(T \bmod (W_A,W_B,W_C)\bigr)+{\bf m}_A+{\bf m}_B+{\bf m}_C,
\end{equation}
which gives us condition (3) and hence concludes the proof.
\end{proof}

\section{Modifications of blocks}\label{Sec5}
In this section, we provide the proof of  \Cref{lemma2}. To get this proof, we use some basic notions of differential geometry. Since we use both 
Euclidean and Zariski open sets, to avoid confusion, we  always specify which kind of topology we are considering. For our purposes,  recall that any Euclidean open set in $\bbC^m$ is  dense in Zariski topology.

For any $m\in\mathbb{Z}_{>0}$ and $\varepsilon\in\bbR_{>0}$, we define the value
\[
r_{m,\varepsilon}\coloneqq \biggl\lceil\frac{1+\varepsilon}{3}m^{2}\biggr\rceil
\]
and we 
can suppose that
$A_i=B_i=C_i=\C^m$ for $i=1,2$. Indeed, for our purposes it suffices to consider spaces up to isomorphism. In particular, we can describe any space using coordinates, by choosing, for instance, the standard basis of $\bbC^m$ and its dual basis in $({\bbC^m})^*$. With this assumptions, we can suppose that 
\[
I_1\otimes I_2=\Mat_{m\times m}(\bbC)={(\bbC^m)}^*\otimes \bbC^m
\]
for each $I = A, B, C$, so that any matrix can be viewed as a linear map in the corresponding space.
We also consider the 
space
\[
     \calV_{m,\varepsilon}  \coloneqq  (A_{1}\oplus B_{1}\oplus C_{1}\oplus A_{2}\oplus B_{2}\oplus C_{2})^{r_{m,\varepsilon}} = \C^{6mr_{m,\varepsilon}},
     \]
     depending on $m$ and $\eps$, and the space
     \[
     \calW_m  \coloneqq (A_{1}\oplus A_{2})\otimes (B_{1}\oplus B_{2})\otimes (C_{1}\oplus C_{2}) = (\C^m \oplus \C^m)^{\otimes 3},
\]
only depending on $m$.
Let us consider the map $\phi\colon \calV_{m, \varepsilon} \to \calW_m$ , defined by
\[
\phi(p) = \sum_{\alpha =1}^{r_{m,\eps}}(x^{\alpha},\xi^{\alpha})\otimes (y^{\alpha},\eta^{\alpha})\otimes (z^{\alpha},\tau^{\alpha}),
\]
for every element
$p=(x^{\alpha},y^{\alpha},z^{\alpha},\xi^{\alpha},\eta^{\alpha},\tau^{\alpha})_{\alpha = 1,\dots,r_{m,\eps}} \in \calV_{m,\varepsilon}$,
and the tensors of $(\bbC^m)^{\otimes 3}$ given by
\begin{gather*}
T^{1}_{\Phi_A,\Phi_B,\Phi_C}(p)\coloneqq\sum_{\alpha =1}^{r_{m,\eps}}  x^{\alpha} \otimes y^{\alpha} \otimes z^{\alpha} + \Phi_{A}(\xi^{\alpha}) \otimes y^{\alpha} \otimes z^{\alpha}+ x^{\alpha} \otimes \Phi_{B}(\eta^{\alpha}) \otimes z^{\alpha}  + x^{\alpha} \otimes y^{\alpha} \otimes \Phi_{C}(\tau^{\alpha}),\\
T^{2}_{\Phi_A,\Phi_B,\Phi_C}(p)\coloneqq\sum_{\alpha =1}^{r_{m,\eps}}  \xi^{\alpha} \otimes \eta^{\alpha} \otimes \tau^{\alpha} + \Phi_{A}(x^{\alpha}) \otimes \eta^{\alpha} \otimes \tau^{\alpha}+\xi^{\alpha} \otimes \Phi_{B}(y^{\alpha}) \otimes \tau^{\alpha} + \xi^{\alpha} \otimes \eta^{\alpha} \otimes \Phi_{C}(z^{\alpha}), 
\end{gather*}
for every  $\Phi_A,\Phi_B,\Phi_C\in\Mat_{m\times m}(\bbC)$.
Then we define the two subsets $V_{m,\varepsilon}^1$, $V_{m,\varepsilon}^2\subset V_{m,\eps}$ as
\[
V_{m,\eps}^i \coloneqq \Set{ p \in V_{m,\eps} | \min_{\Phi_A, \Phi_B, \Phi_C\in\Mat_{m\times m}(\bbC)} \rk \bigl(T^{i}_{\Phi_A,\Phi_B,\Phi_C}(p)\bigr) > \frac{r_{m,\eps}}{2}}
\]
for each $i=1,2$.
The following lemma, 
regarding the sets $V_{m,\eps}^1$ and $V_{m,\eps}^2$, is the key to proving \texorpdfstring{\Cref{lemma2}}{Lemma \ref{sec:foo}}.  
\begin{lemma}\label{lem: key_for_lem2}
There exist $\eps \in \bbR_{>0}$ and  $\bar m_{\eps} \in \mathbb{Z}_{>0}$ such that  $V_{m,\eps}^1\cap V_{m,\eps}^2\neq \varnothing$ for every $m>\bar m_{\eps}$.
\end{lemma}
Before presenting the proof of \Cref{lem: key_for_lem2}, we provide the proof of \Cref{lemma2}. The proof of \Cref{lem: key_for_lem2}, instead, is given in \Cref{section:existence_open}, after introducing the necessary notions. 
\begin{proof}[Proof of \Cref{lemma2}]
By \Cref{lem: key_for_lem2}, there exist $\eps \in \bbR_{>0}$ and $m\in\mathbb{Z}_{>0}$ such that 
\( V_{m,\eps}^1 \cap V_{m,\eps}^2 \neq \varnothing. \) 
Let us consider any point 
\( p \in V_{m,\eps}^1 \cap V_{m,\eps}^2 \) and let \( T = \phi(p) \in W_m \). Then, by the definition of $\phi$, we have $\rk(T)\leq r_{m,\eps}$. Now, for every $\alpha=1,\dots,r_{m,\eps}$ and every $i=1,\dots,m$ we write $x_i^{\alpha},y_i^{\alpha},z_i^{\alpha},\xi_i^{\alpha},\eta_i^{\alpha},\tau_i^{\alpha}$ to denote the $m$ coordinates of  $x^{\alpha},y^{\alpha},z^{\alpha},\xi^{\alpha},\eta^{\alpha},\tau^{\alpha}$, respectively. Let us consider the three blocks
\[
T_{211}=\sum_{\alpha=1}^{r_{m,\eps}}\xi^{\alpha}\otimes y^{\alpha}\otimes z^{\alpha},\qquad T_{121}=\sum_{\alpha=1}^{r_{m,\eps}}x^{\alpha}\otimes \eta^{\alpha}\otimes z^{\alpha},\qquad T_{112}=\sum_{\alpha=1}^{r_{m,\eps}}x^{\alpha}\otimes y^{\alpha}\otimes \tau^{\alpha}.
\]
Interpreting these tensors as linear maps, we have
\[
T_{211}(e_i^*)=\sum_{\alpha=1}^{r_{m,\eps}}\xi_i^{\alpha} (y^{\alpha}\otimes z^{\alpha}),\qquad T_{121}(e_i^*)=\sum_{\alpha=1}^{r_{m,\eps}}\eta_i^{\alpha}(x^{\alpha}\otimes z^{\alpha}),\qquad T_{112}(e_i^*)=\sum_{\alpha=1}^{r_{m,\eps}}\tau_i^{\alpha}(x^{\alpha}\otimes y^{\alpha}).
\]
In particular, we can write any modification of $T_{111}$ by $T_{211}(A_2^*)$ 
as 
\[
T_{111}+\sum_{i,j=1}^{m}\gamma_{ij}e_j\otimes \biggl(\sum_{\alpha=1}^{r_{m,\eps}}\xi_i^{\alpha} (y^{\alpha}\otimes z^{\alpha})\biggr)=T_{111}+\sum_{\alpha=1}^{r_{m,\eps}}\biggl(\sum_{i,j=1}^{m}\xi_i^{\alpha}\gamma_{ij}e_j\biggr)\otimes y^{\alpha}\otimes z^{\alpha}
\]
where $\gamma_{ij}\in\bbC$ for every $i,j=1,\dots,m$. Hence, we can write it as
\[
T_{111}+\sum_{\alpha=1}^{r_{m,\eps}}\Phi_A(\xi^{\alpha})\otimes y^{\alpha}\otimes z^{\alpha},
\]
where $\Phi_A\colon\bbC^m\to\bbC^m$ is the linear map defined on the elements of the standard basis by
\[
\Phi_A(e_i)\coloneqq \sum_{j=1}^{m}\gamma_{ij}e_j, \quad i=1,\dots,m.
\] 
Analogously, any modification of $T_{111}$ by $T_{121}(B_2^*)$ and $T_{112}(C_2^*)$  can be written as 
\[
T_{111}+\sum_{\alpha=1}^{r_{m,\eps}} x^{\alpha}\otimes \Phi_B(\eta^{\alpha})\otimes z^{\alpha},\qquad T_{111}+\sum_{\alpha=1}^{r_{m,\eps}} x^{\alpha}\otimes y^{\alpha}\otimes \Phi_C(\tau^{\alpha}),
\]
respectively, for some matrices $\Phi_B,\Phi_C\in\Mat_{m\times m}(\bbC)$. 
By analogous considerations for the set of modifications of $T_{222}$ by $T_{122}(A_1^*)$, $T_{212}(B_1^*)$, and $T_{221}(C_1^*)$, 
any element of the sets
\[T_{111} \bmod \bigl(T_{211}(A_2^*), T_{121}(B_2^*), T_{112}(C_2^*)\bigr),\qquad T_{222} \bmod \bigl(T_{122}(A_2^*), T_{212}(B_2^*), T_{221}(C_2^*)\bigr)\]
can be written as \[
T^1_{\Phi_A,\Phi_B,\Phi_C}(p),\qquad T^2_{\Phi_A',\Phi_B',\Phi_C'}(p),
\]
respectively, for some $\Phi_A,\Phi_B,\Phi_C,\Phi_A',\Phi_B',\Phi_C'\in\Mat_{m\times m}(\bbC)$. 
In particular, since $p \in V_{m,\eps}^1 \cap V_{m,\eps}^2$, it follows that
\begin{gather*}
\min \rk\Bigl(T_{111} \bmod \bigl(T_{211}(A_2^*), T_{121}(B_2^*), T_{112}(C_2^*)\bigr)\Bigr)
> \frac{r_{m,\eps}}{2},\\
\min \rk\Bigl(T_{222} \bmod \bigl(T_{122}(A_1^*), T_{212}(B_1^*), T_{221}(C_1^*)\bigr)\Bigr) > \frac{r_{m,\eps}}{2},
\end{gather*}
which concludes the proof.
\end{proof}

\subsection{Existence of open dense subsets}\label{section:existence_open}
 Throughout this subsection, we maintain the notation established 
at the beginning of \Cref{Sec5}. Additionally, for any real numbers $ \alpha,\beta,\gamma \in[0,1]\subset\bbR$, we define the two values
\[
\mu_{1}(\alpha,\beta,\gamma)  \coloneqq  \frac{\alpha+\beta+\gamma}{3}+(1-\alpha)(1-\beta)(1-\gamma),
\]
\begin{align}\label{LBnotation}
\nonumber      \mu_{2}(\alpha,\beta,\gamma) & \coloneqq (1-\alpha)(1-\beta)(1-\gamma) 
      +\alpha(1-\beta)(1-\gamma)+  \beta(1-\alpha)(1-\gamma)+\gamma(1-\alpha)(1-\beta)\vphantom{\dfrac{2}{3}} \\[1ex]
 \nonumber& \hphantom{{}\coloneqq{}} + \max \left\{ 0, \, \alpha\beta(1-\gamma)-\dfrac{2}{3}\min \{\alpha,\beta\}\right\}   +   \max \left\{ 0, \, \alpha\gamma(1-\beta)-\dfrac{2}{3}\min \{\alpha,\gamma\}\right\}  \\[1ex]
& \hphantom{{}\coloneqq{}} +\max \left\{ 0, \, \beta\gamma(1-\alpha)-\dfrac{2}{3}\min \{\beta,\gamma\}\right\},   
\end{align}
and we also set 
\begin{equation}\label{LBnotationmin}
 \mu \coloneqq \min\limits_{0\leq\alpha,\beta,\gamma\leq1}\max \bigl\{\mu_{1}(\alpha,\beta,\gamma),\mu_{2}(\alpha,\beta,\gamma)\bigr\}.   
\end{equation} 
Given arbitrary differential maps 
\begin{equation}\label{formula:psi_maps}
\psi_A\colon U_{m,\eps} \to A_1 \otimes A_2, \quad \psi_B\colon U_{m,\eps} \to B_1 \otimes B_2, \quad \psi_C\colon U_{m,\eps} \to C_1 \otimes C_2
\end{equation}
defined over any Euclidean open subset $U_{m,\eps}\subset \calV_{m,\eps}$,  
we write $\Psi_I^p$ to denote the image $\psi_{I}(p) \in I_1 \otimes I_2$ of any point $p = (x^{\alpha},y^{\alpha},z^{\alpha},\xi^{\alpha},\eta^{\alpha},\tau^{\alpha})_{\alpha =1,\dots,r_{m,\eps}}\in U_{m,\eps}$, for each $I = A, B, C$. 
\begin{remark}
   The notation above is needed to avoid confusion with each map $\psi_I$ and to make evident that $\psi_{I}(p) \coloneqq \Psi_I^p$
denotes indeed a matrix. 
\end{remark}
 We also define the two spaces
\begin{align}\label{formula:spaces_P_Q}
\nonumber\calP_{\psi_A,\psi_B,\psi_C}^p  &\coloneqq    \bigl\langle p_{A}\otimes y^{\alpha}\otimes z^{\alpha}\bigr\rangle_{p_A \in \mathrm{Im}\Psi_A^p}^{\alpha=1, \ldots, r_{m,\eps}}+\bigl\langle x^{\alpha}\otimes p_{B}\otimes z^{\alpha}\bigr\rangle_{p_B \in \mathrm{Im}\Psi_B^p}^{\alpha=1, \ldots, r_{m,\eps}}+ \bigl\langle x^{\alpha}\otimes y^{\alpha}\otimes p_{C}\bigr\rangle_{p_C \in \mathrm{Im}\Psi_C^p}^{\alpha=1, \ldots, r_{m,\eps}} \\[0.5ex]
\nonumber\calQ_{\psi_A,\psi_B,\psi_C}^p &\coloneqq   \bigl\langle  q_A \otimes  y^{\alpha}\otimes z^{\alpha}+q_A \otimes\Psi_{B}^p(\eta^{\alpha})\otimes z^{\alpha}+q_A \otimes y^{\alpha}\otimes \Psi_{C}^p(\tau^{\alpha})\bigr\rangle_{q_A\in \C^m}^{\alpha=1, \ldots, r_{m,\eps}} \\[0.5ex]
\nonumber  &\hphantom{{}\coloneqq{}}+  \bigl\langle x^{\alpha}\otimes q_B\otimes z^{\alpha}+ \Psi_{A}^p(\xi^{\alpha})\otimes q_B\otimes z^{\alpha}+x^{\alpha}\otimes q_B\otimes \Psi_{C}^p(\tau^{\alpha})\bigr\rangle_{q_B\in \C^m}^{\alpha=1, \ldots, r_{m,\eps}} \\[0.5ex]
 &\hphantom{{}\coloneqq{}}+ \bigl\langle x^{\alpha}\otimes y^{\alpha}\otimes q_C+ \Psi_{A}^p(\xi^{\alpha})\otimes y^{\alpha}\otimes q_C+x^{\alpha}\otimes \Psi_{B}^p(\eta^{\alpha})\otimes q_C\bigr\rangle_{q_C\in \C^m}^{\alpha=1, \ldots, r_{m,\eps}}.
\end{align}   
In order to complete our proof, we need to introduce the following real positive values, depending on $k,m\in\mathbb{Z}_{>0}$, which are defined as
\begin{gather*}
L^{-}_{k,m}\coloneqq \max \biggl\{\frac{9}{k+1}-\frac{(6k-2)m^{2}-(9k^{2}+9)m+(3k^{3}+9k+6)}{m^2(3m-3k-2)},\frac{2m^{2}+9m-6}{m^2(3m-2)}\biggr\},\\[1ex]
L^{+}_{k,m}\coloneqq \min \biggl\{ (2\mu-1)+\frac{2\bigl((2\mu-1)-9k-8\bigr)m^2-9m+6}{m^2(3m-2)}, \frac{2\mu-1}{3}+\frac{2\bigl((2\mu-1)-24\bigr)m^2-27m+18}{3m^2 (9m-2)} \biggr\}.
\end{gather*}
In the proof of \Cref{lem: key_for_lem2}, we  assume the existence of three maps as in \eqref{formula:psi_maps}. The following proposition concerns  
 lower bounds for the dimension of the sum of the related two spaces defined in \eqref{formula:spaces_P_Q}.
\begin{proposition}\label{df}
Let $k,m\in \mathbb{Z}_{>0}$, with $k\leq m-1$, let $\varepsilon\in\bbR_{>0}$, and let $\psi_A,\psi_B,\psi_C\colon U_{m,\eps}\to\Mat_{m\times m}(\bbC)$ be differentiable maps defined over a non-empty Euclidean open set $U_{m,\eps}\subset \calV_{m,\eps}$. If $L^{-}_{k,m}\leq \varepsilon$, then there exists a point  $\bar p\in U_{m,\eps}$ such that
\[\dim \bigl(\calP_{\psi_A,\psi_B,\psi_C}^{\bar p}+\calQ_{\psi_A,\psi_B,\psi_C}^{\bar{p}}\bigr)\geq m^2\bigl(m\mu-\max \{3k,m\varepsilon\}\bigr).
\]
\end{proposition}
The proof of \Cref{df}  requires certain  results, concerning lower bounds on dimensions. \Cref{sec:technical_lower_bounds} is fully dedicated to its proof. However, this result suffices to get the proof of \Cref{lem: key_for_lem2}, which is based on the existence of a non-empty Zariski open subset of the spaces $V_{m,\eps}^1$ and $V_{m,\eps}^2$. 
\begin{proof}[Proof of \Cref{lem: key_for_lem2}]
We prove 
that the set $V^{1}_{m,\varepsilon}$ contains a non-empty Zariski open set. The same result also holds for $V^{2}_{m,\varepsilon}$, and the procedure is completely analogous.
By contradiction, assume  that $V^{1}_{m,\varepsilon}$ does not contain a Zariski open set or, equivalently, 
that the complement 
$\calV_{m,\varepsilon}\setminus V^{1}_{m,\varepsilon}$ is Zariski dense in $\calV_{m,\varepsilon}$.
Let us consider the maps
\[
\pi_1\colon  \calV_{m,\varepsilon}\times  \bigl(\mathrm{Mat}_{m\times m}(\C)\bigr)^{3}\to \calV_{m,\varepsilon},\qquad F\colon  \calV_{m,\varepsilon}\times  \bigl(\mathrm{Mat}_{m\times m}(\C)\bigr)^{3}\to (\mathbb{C}^{m})^{\otimes 3},
\]
where $\pi_1$ is the projection to the first factor and the map $F$ is defined by
\[
F(p, \Psi_A, \Psi_B, \Psi_C) \coloneqq T_{\Psi_A,\Psi_B,\Psi_C}^1(p)
\]
for every $p\in V$ and for every $\Psi_A$, $\Psi_B$,  $\Psi_C\in\mathrm{Mat}_{m\times m}(\C)$. 
Let $Y$ be the preimage  of the affine cone over the $(r_{m,\eps}/2)$-th secant variety of the Segre variety $\Sigma_3^m=\Seg\bigl((\mathbb{P}\mathbb{C}^{m})^3\bigr)$ under $F$, that is,
\[Y \coloneqq F^{-1}\bigl( \hat{\sigma}_{\frac{r_{m,\eps}}{2}}(\Sigma_3^m)\bigr).\]
From the definition of $V^{1}_{m,\varepsilon}$, it follows that $ \calV_{m,\varepsilon}\setminus  V^{1}_{m,\varepsilon} \subseteq \pi_1(Y)$, which means that $\pi_1(Y)$ is dense in 
$\calV_{m,\varepsilon}$. Since 
$Y$ is an algebraic variety, and hence a constructible set, it follows from Chevalley’s theorem that \( \pi_1(Y) \) is also constructible as the image of $Y$
under a regular map (see, e.g., \cite[Theorem 4.19]{BerndMateusz}). Therefore, there exists a non-empty Zariski open subset \( U' \subset \pi_1(Y) \subset \calV_{m,\varepsilon} \). Now, 
considering the restriction
$\pi_1\colon Y \cap \pi_1^{-1}(U') \to U'$,
by the submersion property (see e.g.,~\cite[Theorem 4.26]{Lee13}), 
 there exists a differentiable section  \[ s\colon U_{m,\varepsilon}\to  \calV_{m,\eps}\times \bigl(\Mat_{m\times m}(\bbC)\bigr)^3, \] defined on a non-empty  Euclidean open subset \( U_{m,\eps} \subset U' \subset \calV_{m,\varepsilon}, \)
such that $s(U_{m,\eps})\subset Y\cap\pi_1^{-1}(U')$.

In particular, by restricting to $U_{m,\eps}\subset \calV_{m,\eps}$ we can write $s=(\id,\psi_A,\psi_B,\psi_C)$, 
where $\id\colon U_{m,\eps}\to U_{m,\eps}$ is the identity map, and $\psi_{I}\colon  U_{m, \varepsilon}\to \mathrm{Mat}_{m\times m}(\C)$ is a differentiable map for each $I=A,B,C$. Therefore, for any $p\in U_{m, \varepsilon}$, considering the map $f\coloneqq F\circ s\colon U_{m,\eps}\to(\bbC^{m})^{\otimes 3}$ we have
\begin{equation}\label{eq: f(U)}
 f(p)=F\bigl(p,\psi_{A}(p),\psi_{B}(p),\psi_{C}(p)\bigr) = F\bigl(p,\Psi_{A}^p,\Psi_{B}^p,\Psi_{C}^p\bigr) \in \hat{\sigma}_{\frac{r_{m,\eps}}{2}}(\Sigma_3^m), 
\end{equation}
 For any point $p\in U_{m,\eps}$, we define the map
$\tilde{f}_p\colon \calV_{m,\varepsilon} \to  (\mathbb{C}^{m})^{\otimes 3}$ 
 by 
\begin{align*}
\tilde{f}_p(v) =   \displaystyle\sum_{\alpha=1}^{r_{m,\eps}} &\Bigl( \Psi_{A}^p(v_{4}^{\alpha})\otimes y^{\alpha}\otimes z^{\alpha}+ x^{\alpha}\otimes \Psi_{B}^p(v_{5}^{\alpha})\otimes z^{\alpha}+ x^{\alpha}\otimes y^{\alpha}\otimes \Psi_{C}^p(v_{6}^{\alpha})+ v_{1}^{\alpha}\otimes  y^{\alpha}\otimes z^{\alpha} \\[1ex]
 &\, +v_{1}^{\alpha}\otimes\Psi_{B}^p(\eta^{\alpha})\otimes z^{\alpha}+v_{1}^{\alpha}\otimes y^{\alpha}\otimes \Psi_{C}^p(\tau^{\alpha})+  x^{\alpha}\otimes v_{2}^{\alpha}\otimes z^{\alpha}+ \Psi_{A}^p(\xi^{\alpha})\otimes v_{2}^{\alpha}\otimes z^{\alpha}\\[2ex]
 &\,+x^{\alpha}\otimes v_{2}^{\alpha}\otimes \Psi_{C}^p(\tau^{\alpha})+ x^{\alpha}\otimes y^{\alpha}\otimes v_{3}^{\alpha}+ \Psi_{A}^p(\xi^{\alpha})\otimes y^{\alpha}\otimes v_{3}^{\alpha}+x^{\alpha}\otimes \Psi_{B}^p(\eta^{\alpha})\otimes v_{3}^{\alpha}\Bigr),\end{align*}
for every point
$v =(v_{1}^{\alpha},\ldots,v_{6}^{\alpha})_{\alpha = 1,\dots,r_{m,\eps}}\in \calV_{m,\varepsilon}$.
This map coincides with $\dd{f_p}$, which is the differential of $f$ at the point $p$, on the 
subspace
\[D\coloneqq \Ker \dd(\psi_{A})_p\cap \Ker \dd(\psi_{B})_{p}\cap \Ker \dd(\psi_{C})_{p}\subset \calV_{m,\varepsilon}.\]
We then obtain the chain of inequalities
\begin{align*}
\dim \bigl(\tilde{f}_{p}(\calV_{m,\varepsilon})/\tilde{f}_{p}(D)\bigr) &\leq \dim (\calV_{m,\varepsilon}/D) \\
&\leq \dim \bigl(\calV_{m,\varepsilon}/\Ker \dd(\psi_{A})_{p}\bigr)+\dim (\calV_{m,\varepsilon}/\Ker \dd(\psi_{B})_{p}\bigr)+\dim (\calV_{m,\varepsilon}/\Ker \dd(\psi_{C})_{p}\bigr) \\
& \leq 3\dim \bigl(\mathrm{Mat}_{m\times m}(\C)\bigr)=3m^{2},
\end{align*}
from which 
it follows that
\begin{equation}\label{formula:inequality_differential}
\dim (\im\dd{f_p})\geq\dim (\im\dd{f_p}|_{D})=\dim (\im\tilde{f}_{p}|_{D})\geq \dim (\im \tilde{f}_{p})-3m^{2}.
\end{equation}
By definition, we have that 
\[\im \tilde{f}_{p}=\calP_{\psi_A,\psi_B,\psi_C}^p+\calQ_{\psi_A,\psi_B,\psi_C}^p,\]
for every $p\in U_{m,\eps}$.
Then, if $L_{k,m}^{-}\leq \eps$ for some $k\leq m-1$, by \Cref{df} there exists a point $\bar p\in U_{m,\eps}$ such that
\[
\dim\bigl(\calP_{\psi_A,\psi_B,\psi_C}^{\bar p}+\calQ_{\psi_A,\psi_B,\psi_C}^{\bar p}\bigr)\geq m^2\bigl(m\mu-\max \{3k,m\varepsilon\}\bigr).
\]
In particular, by inequality \eqref{formula:inequality_differential}, we have
\begin{equation}\label{dim_im_1}
\dim (\im\dd{f_{\bar p}})\geq m^2\bigl(m\mu-\max \{3k,m\varepsilon\}-3\bigr),
\end{equation}
On the other hand, by the definition of $f$ in formula \eqref{eq: f(U)}, we have
\begin{equation}\label{formula:inclusion_f(U)}
    f(U_{m,\eps}) \subset \hat{\sigma}_{\frac{r_{m,\eps}}{2}}(\Sigma_3^m).
\end{equation}
Now, by formula \eqref{formula:secant_variety_Segre}, we have 
\begin{equation}\label{formula:upper_bound_secant_variety}
\dim \bigl(\hat{\sigma}_{\frac{r_{m,\eps}}{2}}(\Sigma_3^m)\bigr)=\min\biggl\{\frac{r_{m,\eps}}{2}(3m-2),m^3\biggr\} \leq \frac{r_{m,\eps}}{2} (3m-2).
\end{equation}
In the case where $\eps< L_{k,m}^+$, we have the chain of inequalities
\begin{equation}
\label{dim_im_2}
 \dim (\im\dd{f_{\bar p}}) \leq \dim \bigl(\hat{\sigma}_{\frac{r_{m,\eps}}{2}}(\Sigma_3^m)\bigr) \leq \frac{r_{m,\eps}}{2} (3m-2)< m^2\bigl(m\mu-\max \{3k,m\varepsilon\}-3\bigr),
\end{equation}  
where the first two inequalities are given by  inclusion \eqref{formula:inclusion_f(U)} and the upper bound \eqref{formula:upper_bound_secant_variety}, respectively,
while the last inequality is obtained by \Cref{prop:epsCondition}.
According to \Cref{prop:2conditions}, we can determine 
$\eps \in \bbR_{>0}$  and \(k,\bar{m}_{\eps}\in\bbZ_{>0}\), with $k\leq \bar{m}_{\eps}-1$, such that $L_{k,m}^{-}\leq\varepsilon<L_{k,m}^{+}$ for every integer $m>\bar{m}_{\eps}$. For such $m,k\in\bbZ_{>0}$ and $\eps\in \bbR_{>0}$,
both inequalities \eqref{dim_im_1} and \eqref{dim_im_2} hold, leading to a contradiction. In particular, there exists a real value $\eps \in \bbR_{>0}$ and a positive integer $\bar {m}_{\eps}$ such that, for every $m>\bar{m}_{\eps}$, the set $V_{m,\eps}^1$ contains a non-empty Zariski open set. Repeating the same procedure for $V_{m,\eps}^2$, the statement follows.
\end{proof}
\begin{remark}\label{rem: approximatevalues}
   Based on computations in \texttt{MATLAB} \cite{MATLAB}, involving the discretization of the cube $[0,1]^3$ into subintervals of the length $0.001$ along each axis, we get $\mu\approx 0.52733$. Setting $\mu=0.52733$, computations  
   indicate that the minimal value of $m$ such that $L_{k,m}^{-}\leq\varepsilon<L_{k,m}^{+}$ is $m=\numprint{48352}$ with $k=328$. This results in the following range for $r_{m,\eps}$:  
\[
\numprint{790097248}\leq r_{m,\eps}\leq \numprint{790097406}.
\]
In recent work \cite[Remark 6.11]{Shi24}, Yaroslav Shitov demonstrates a counterexample to Strassen's additivity conjecture, involving tensors of order at most $10^{888}$. 
\end{remark}
\subsection{Technical lower bounds}\label{sec:technical_lower_bounds}
We start by proving the following   
lemma, for which we define the space
\[
\calO_r^{p_1}\coloneqq \bigl\langle q_A\otimes y^{\alpha}\otimes z^{\alpha}\bigr\rangle^{\alpha =1,\ldots,r}_{q_A\in \mathbb{C}^{m}} +\bigl\langle  x^{\alpha} \otimes q_B\otimes z^{\alpha}\bigr\rangle^{\alpha =1,\ldots,r}_{q_B\in \mathbb{C}^{m}} +\bigl\langle x^{\alpha} \otimes y^{\alpha}\otimes q_C\bigr\rangle^{\alpha =1,\ldots,r}_{q_C\in \mathbb{C}^{m}}
\]
for any point  $p_1=(x^{\alpha},y^{\alpha},z^{\alpha})_{\alpha =1,\dots,r}\in \bbC^{3mr}$ with $r\in\bbZ_{>0}$.
\begin{lemma}\label{lem:first_existence}
  Let $r\in\bbZ_{>0}$ be such that $r\geq \rmr_{\mathrm{gen}}^m$. Then there exists a non-empty Zariski open set $\calE_{m,r}\subset\bbC^{3mr}$ such that
$\calO_r^{p_1}=(\bbC^m)^{\otimes 3}$
  for every $p_1\in\calE_{m,r}$.
\end{lemma}
\begin{proof}
Let us consider the morphism $\tau\colon \mathbb{C}^{3mr}\to  (\mathbb{C}^{m})^{\otimes 3}$
defined by
\[
\tau(p_1)  =  \displaystyle{\sum_{\alpha =1}^{r}x^{\alpha}\otimes y^{\alpha}\otimes z^{\alpha}},
\]
for any $p_1=(x^{\alpha},y^{\alpha},z^{\alpha})_{\alpha=1,\dots,r}$.
Since $r\geq \rmr_{\mathrm{gen}}^m$, we have that the $r$-th secant variety of the Segre variety is the whole space, that is
\[
\sigma_{r}(\Sigma_{3}^m)=\mathbb{P}\bigl((\bbC^m)^{\otimes 3}\bigr).
\]
It follows that the map $\tau$ is dominant, that is, the image is a dense set.
 Thus, using \cite[Proposition III.10.4 and Corollary III.10.7]{Har77}, there exists a non-empty Zariski open subset $\calE_{m,r}\subset\bbC^{3mr}$ such that the differential $\dd{\tau_{p_1}}$ is surjective for every point $p_1=(x^{\alpha},y^{\alpha},z^{\alpha})_{\alpha =1,\dots,r}\in\calE_{m,r}$. Observing that $\dd \tau_{p_1}=\calO_{r}^{p_1}$ for any $p_1\in\bbC^{3mr}$, the statement follows.
\end{proof}

For any differentiable maps $\psi_{A},\psi_{B},\psi_{C}\colon U_{m,\varepsilon} \to \mathrm{Mat}_{m\times m}(\C)$, defined 
 over an open set $U_{m,\eps} \subseteq \calV_{m,\eps}$, and any
point $p\in U_{m,\eps}$, we define the values
 \[\delta_{\psi_A}^p\coloneqq \frac{\dim (\im\Psi_A^p)}{m}, \qquad \delta_{\psi_B}^p\coloneqq \frac{\dim (\im\Psi_B^p)}{m}, \qquad \delta_{\psi_C}^p\coloneqq \frac{\dim (\im\Psi_C^p)}{m}.\]
In particular, observe that $0\leq \delta_{\psi_A}^p, \delta_{\psi_B}^p, \delta_{\psi_C}^p \leq 1$.
For the statements later on, we define the value 
\begin{align*}
\Delta_{\psi_A,\psi_B,\psi_C}^p&\coloneqq{} \delta_{\psi_A}^p(1-\delta_{\psi_B}^p)(1-\delta_{\psi_C}^p)+\delta_{\psi_B}^p(1-\delta_{\psi_A}^p)(1-\delta_{\psi_C}^p)+\delta_{\psi_C}^p(1-\delta_{\psi_A}^p)(1-\delta_{\psi_B}^p)\\&\hphantom{{}\coloneqq{}}+(1-\delta_{\psi_A}^p)(1-\delta_{\psi_B}^p)(1-\delta_{\psi_C}^p)
=1-\delta_{\psi_A}^p\delta_{\psi_B}^p-\delta_{\psi_A}^p\delta_{\psi_C}^p-\delta_{\psi_B}^p\delta_{\psi_C}^p+2\delta_{\psi_A}^p\delta_{\psi_B}^p\delta_{\psi_C}^p\leq 1
\end{align*}
for any $p\in U_{m, \varepsilon}$. 
 Moreover, 
 for any $p\in U_{m,\eps}$, we define the space 
\begin{align*}
 \calQ'^p_{\psi_A,\psi_B,\psi_C}&\coloneqq \bigl\langle q_A\otimes \bigl(y^{\alpha}+\Psi_{B}^{p}(\eta^{\alpha})\bigr)\otimes \bigl(z^{\alpha}+\Psi_{C}^{p}(\tau^{\alpha})\bigr)\bigr\rangle^{\alpha =1,\ldots,r_{m,\eps}}_{q_A\in \mathbb{C}^{m}} \\
 &\hphantom{{}\coloneqq{}} +\bigl\langle  \bigl(x^{\alpha} +\Psi_{A}^{p}(\xi^{\alpha})\bigr)\otimes q_B\otimes \bigl(z^{\alpha}+\Psi_{C}^{p}(\tau^{\alpha})\bigr)\bigr\rangle^{\alpha =1,\ldots,r_{m,\eps}}_{q_B\in \mathbb{C}^{m}} \\
 &\hphantom{{}\coloneqq{}} +\bigl\langle \bigl(x^{\alpha}+\Psi_{A}^{p}(\xi^{\alpha})\bigr) \otimes \bigl(y^{\alpha}+\Psi_{B}^{p}(\eta^{\alpha})\bigr)\otimes q_C\bigr\rangle^{\alpha =1,\ldots,r_{m,\eps}}_{q_C\in \mathbb{C}^{m}}.
\end{align*}
In the following, we consider the two projections on the two components of $\calV_{m,\eps}$, 
denoted as
\[
\Pi_1\colon \calV_{m,\eps}\to (A_1\oplus B_1\oplus C_1)^{r_{m,\eps}}=\bbC^{3mr_{m,\eps}},\quad
\Pi_2\colon \calV_{m,\eps}\to (A_2\oplus B_2\oplus C_2)^{r_{m,\eps}}=\bbC^{3mr_{m,\eps}}.
\]
For any non-empty Euclidean open subset $U\subset \calV_{m,\eps}$ and any $q_2\in\Pi_2(U)$, we consider the fiber $$\Pi_2|_{U}^{-1}(q_2)=\Pi_2^{-1}(q_2)\cap U,$$ which is a non-empty Euclidean open set of $\Pi_2^{-1}(q_2)$. In particular, we define the set \[\Omega_{m,\eps}^{q_2}(U)\coloneqq\Pi_1\bigl(\Pi_2|_{U}^{-1}(q_2)\bigr)\simeq \Pi_2|_{U}^{-1}(q_2),\] which is a non-empty Euclidean open set of $\bbC^{3mr_{m,\eps}}$.
 
 The following lemma relies on some basic facts in differential geometry. We recall that we use both Euclidean and Zariski topology.
\begin{lemma}\label{lem:diff_geo}
Let $k,m\in \mathbb{Z}_{>0}$ with $k\leq m-1$. Let $\varepsilon\in\bbR_{>0}$ and  $\psi_A,\psi_B,\psi_C\colon U_{m,\eps}\to\Mat_{m\times m(\bbC)}$ be differentiable maps defined over a non-empty Euclidean open set $U_{m,\eps}\subset \calV_{m,\eps}$.
If $\eps\geq L_{k,m}^-$, then for any $q_2\in \Pi_2(U_{m,\eps})$ there exist three general linear maps $\pi_A,\pi_B,\pi_C\colon \bbC^m\to\bbC^{m-k}$ and a non-empty Euclidean open set $\Omega_{m,\eps}'\subset \Omega_{m,\eps}^{q_2}(U_{m,\eps})$ such that, for every $p_1'\in\Omega_{m,\eps}'$,
 \[
(\pi_A\otimes \pi_B\otimes \pi_C)\bigl( \calQ'^{(p_1',q_2)}_{\psi_A,\psi_B,\psi_C}\bigr)=(\bbC^{m-k})^{\otimes 3}.
 \]
\end{lemma}
\begin{proof}
Let $
q_2=(\xi^{\alpha},\eta^{\alpha},\tau^{\alpha})_{\alpha =1,\dots,r_{m,\eps}}\in\Pi_2(U_{m,\eps})$ and
let $g^{q_2}\colon \Omega_{m,\eps}^{q_2}(U_{m,\eps})\to\bbC^{3mr_{m,\eps}}$ be a differentiable map defined as
\[
g^{q_2}(p_1)=\begin{pNiceMatrix}[cell-space-limits = 2pt]
x^{\alpha}+\psi_A(p_1,q_2)(\xi^{\alpha})\\
y^{\alpha}+\psi_B(p_1,q_2)(\eta^{\alpha})\\
z^{\alpha}+\psi_C(p_1,q_2)(\tau^{\alpha})
\end{pNiceMatrix}_{\alpha=1,\dots,r_{m,\eps}}=
\begin{pNiceMatrix}[cell-space-limits = 2pt]
x^{\alpha}+\Psi_A^{(p_1,q_2)}(\xi^{\alpha})\\
y^{\alpha}+\Psi_B^{(p_1,q_2)}(\eta^{\alpha})\\
z^{\alpha}+\Psi_C^{(p_1,q_2)}(\tau^{\alpha})
\end{pNiceMatrix}_{\alpha=1,\dots,r_{m,\eps}}
\]
for every point $p_1=(x^{\alpha},y^{\alpha},z^{\alpha})_{\alpha=1,\dots,r_{m,\eps}}\in\Omega_{m,\eps}^{q_2}(U_{m,\eps})$. 
The differential of $g^{q_2}$ at $p_1$ is the linear map
$
\dd (g^{q_2})_{p_1}={\id}_{3mr_{m,\eps}}+L^{p_1},
$ 
where the second summand is the linear map $L^{p_1}\colon\bbC^{3mr_{m,\eps}}\to\bbC^{3mr_{m,\eps}}$, defined as
\[
L^{p_1}(q_1)\coloneqq \begin{NiceMatrix}[margin]
\dd(\psi_{A})_{(p_1,q_2)}(q_1,0)(\xi^{1})&\Block{3-1}{\quad\alpha=1} \\
\dd(\psi_{B})_{(p_1,q_2)}(q_1,0)(\eta^{1})& \\
\dd(\psi_{C})_{(p_1,q_2)}(q_1,0)(\tau^{1})& \\
\vdots & \\
\dd(\psi_{A})_{(p_1,q_2)}(q_1,0)(\xi^{r_{m,\eps}})&\Block{3-1}{\quad\,\, \alpha=r_{m,\eps}} \\
\dd(\psi_{B})_{(p_1,q_2)}(q_1,0)(\eta^{r_{m,\eps}})&\\
\dd(\psi_{C})_{(p_1,q_2)}(q_1,0)(\tau^{r_{m,\eps}})&
\CodeAfter
 \SubMatrix({1-1}{7-1})
 \SubMatrix{.}{1-1}{3-1}{\}}[xshift=4mm]
 \SubMatrix{.}{5-1}{7-1}{\}}[xshift=4mm]
\end{NiceMatrix}
\]
for every $q_1\in\bbC^{3mr_{m,\eps}}$. 
Consider the linear map 
\[
\bigl(\dd(\psi_A)_{(p_1,q_2)},\dd(\psi_B)_{(p_1,q_2)},\dd(\psi_C)_{(p_1,q_2)}\bigr)\colon\bbC^{3mr_{m,\eps}}\to\bigl(\Mat_{m\times m}(\bbC)\bigr)^3,
\] 
given by the Cartesian product of the differentials $\dd(\psi_A)_{(p_1,q_2)}$, $\dd(\psi_B)_{(p_1,q_2)}$, $\dd(\psi_C)_{(p_1,q_2)}$. Then we can write \[
L^{p_1}=\Lambda^{q_2}\circ \bigl(\dd(\psi_A)_{(p_1,q_2)},\dd(\psi_B)_{(p_1,q_2)},\dd(\psi_C)_{(p_1,q_2)}\bigr),
\]
where $\Lambda^{q_2}\colon(\Mat_{m\times m})^3\to\bbC^{3mr_{m,\eps}}$ is the linear map defined by
\[
\Lambda^{q_2}(M_1,M_2,M_3)=\bigl(M_1(\xi^1),M_2(\eta^1),M_3(\tau^1),\dots,M_1(\xi^{r_{m,\eps}}),M_2(\eta^{r_{m,\eps}}),M_3(\tau^{r_{m,\eps}})\bigr),
\]
for every $M_1,M_2,M_3\in\Mat_{m\times m}(\bbC)$. Since
\[
\rk\bigl(\dd(\psi_A)_{(p_1,q_2)},\dd(\psi_B)_{(p_1,q_2)},\dd(\psi_C)_{(p_1,q_2)}\bigr)\leq 3m^2,
\]
we also have $\rk(L^{p_1})\leq 3m^2$. So we obtain the relations
\begin{align*}
3mr_{m,\eps}&=\rk(\id_{3mr_{m,\eps}})\leq \rk(\id_{3mr_{m,\eps}}+L^{p_1})+\rk(-L^{p_1})\leq\rk\bigl(\dd {(g^{q_2})}_{p_1}\bigr)+3m^{2}
\end{align*} 
and hence, $\rk\bigl(\dd (g^{q_2})_{p_1}\bigr)\geq 3mr_{m,\eps}-3m^{2}$.
Now, for every $\alpha=1,\dots,r_{m,\eps}$ and each $I=A,B,C$, we consider the set \[
\calB_{\alpha,I}=\{e_1^{\alpha,I},\dots,e_m^{\alpha,I}\},
\] 
as a copy of the canonical basis of $\bbC^m$. 
 We can write
\[
\dd{(g^{q_2})}_{p_1}=
\sum_{I=A,B,C}\sum_{\alpha=1}^{r_{m,\eps}}\sum_{j=1}^{k_{\alpha,I}}\phi_j^{\alpha,I}\otimes e_j^{\alpha,I},
\]
for some $1\leq k_{\alpha,I}\leq m$,
$\phi_j^{\alpha,I}\in(\bbC^{m})^*$ for every $\alpha=1,\dots,r_{m,\eps}$ and each $I=A,B,C$. By definition
\[
\rk\bigl( \dd {(g^{q_2})}_{p_1}\bigr)=\sum_{I=A,B,C}\sum_{\alpha=1}^{r_{m,\eps}}k_{\alpha,I}.
\]
Let us consider the set $[r_{m,\eps}]\coloneqq\{1,\dots,r_{m,\eps}\}$, and let us define the sets 
\begin{gather*}
\Gamma_{k}\coloneqq \Set{(\alpha,I)\in[r_{m,\eps}]\times\{A,B,C\}\vphantom{\Bigl(}| k_{\alpha,I}< m-k},\\
\Gamma\coloneqq\Set{\alpha\in [r_{m,\eps}]\vphantom{\Bigl(} | \exists I\in\{A,B,C\}:(\alpha,I)\in \Gamma_{k}}.
\end{gather*} 
Denoting $d_{k}\coloneqq\lvert \Gamma_{k}\rvert$ and $d\coloneqq\lvert \Gamma\rvert$, we have that
\[
\rk\bigl(\dd{(g^{q_2})}_{p_1}\bigr)\leq (m-k-1)d_{k}+m(3r_{m,\eps}-d_{k}),
\]  
which implies that
\begin{equation}\label{formula:inequality_d}
d \leq d_{k}\leq \frac{3mr_{m,\eps}-\rk\bigl(\dd{(g^{q_2})}_{p_1}\bigr)}{k+1} \leq \frac{3m^{2}}{k+1}.
\end{equation}
Defining the two spaces $H_{\Gamma}\coloneqq\langle\calB_{\alpha,I}\rangle_{\alpha\in\Gamma}^{I=A,B,C}$ and $H'_{\Gamma}\coloneqq\langle\calB_{\alpha,I}\rangle_{\alpha\not\in\Gamma}^{I=A,B,C}$,
we have
$\bbC^{3mr_{m,\eps}}=H_{\Gamma}\oplus H_{\Gamma}'$ and
we can define the projection
\[
\pi_{\Gamma}'\colon \mathbb{C}^{3mr_{m,\eps}}\to H_{\Gamma}'\simeq\mathbb{C}^{3m(r_{m,\eps}-d)}.
\]
For any $(\alpha,I)\in[r_{m,\eps}]\times\{A,B,C\}$, 
we denote
\[
 \dd {(g^{q_2})}_{p_1}^{(\alpha,I)}\coloneqq \dd(g^{q_2})_{p_1}\big|_{(\bbC^{3mr_{m,\eps}})^*\otimes{\langle\calB_{\alpha,I}\rangle}}\in\Mat_{m\times 3mr_{m,\eps}}(\bbC)
\]
and, for any $(\alpha,I)\not\in\Gamma_{k}$, 
we define the open subset
\[
\Theta_{\alpha,I}\coloneqq \Set{M\in\mathrm{Mat}_{(m-k)\times m}(\mathbb{C})|\rk\bigl(M\circ\dd {(g^{q_2})}_{p_1}^{(\alpha,I)}\bigr)=m-k }\subset\mathrm{Mat}_{(m-k)\times m}(\bbC).
\]
By definition of $\Gamma_{k}$, we have
\[
\rk\bigl(\dd {(g^{q_2})}_{p_1}^{(\alpha,I)}\bigr)\geq m-k,
\]
which implies, in particular, that $\Theta_{\alpha,I}\neq\varnothing$. Therefore, defining the set
\[
\Gamma_I\coloneqq\Set{\alpha\in[r_{m,\eps}]| (\alpha,I)\not\in \Gamma_{k}}
\]
for every $I=A,B,C$, so that $[r_{m,\eps}]\setminus \Gamma=\Gamma_A\cap\Gamma_B\cap\Gamma_C$,
we can choose a general linear map
\[
\pi_I\in\bigcap_{\alpha\in \Gamma_I}\Theta_{\alpha,I}\neq\varnothing
\]
for each $I=A,B,C$. In particular, we can consider the projection
$
(\pi_{A},\pi_{B},\pi_{C})\colon\bbC^{3m}\to \bbC^{3(m-k)}$
and define the surjective map
\[
\Pi\coloneqq (\pi_{A},\pi_{B},\pi_{C})^{r_{m,\eps}-d}\circ\pi_{\Gamma}'\colon \mathbb{C}^{3mr_{m,\eps}}\to \bbC^{3(m-k)(r_{m,\eps}-d)}
\]
and the composition
$G^{q_2}\coloneqq \Pi\circ g^{q_2}\colon \Omega_{m,\eps}^{q_2}(U_{m,\eps})\to \mathbb{C}^{3(m-k)(r_{m,\eps}-d)}$.
By definition of $\Theta_{\alpha,I}$, the map $\dd (G^{q_2})_{p_1}$ is surjective. Therefore, by~\cite[Proposition~4.1]{Lee13}, there exists an Euclidean open subset $W_{m,\eps}\subset\Omega_{m,\eps}^{q_2}(U_{m,\eps})$ such that
\[
\dd (G^{q_2})_{q_1}=\dd(\Pi\circ g^{q_2})_{q_1}=\Pi\circ \dd(g^{q_2})_{q_1}=\Pi\circ(\id_{3mr_{m,\eps}}+L^{p_1})
\] 
is full rank for any $q_1\in W_{m,\eps}$, that is, 
\[
\rk\bigl(\dd{(G^{q_2})}_{q_1}\bigr)=3(m-k)(r_{m,\eps}-d).
\] 
In particular, by the property of submersions (see e.g.,~\cite[Theorem 4.26]{Lee13}), $G^{q_2}(W_{m,\eps})$ is a non-empty Euclidean open set and hence a dense subset in $\bbC^{3(m-k)(r_{m,\eps}-d)}$. 
By hypothesis, we have $\eps\geq L_{k,m}^{-}$, hence, by inequality \eqref{formula:inequality_d} and  \Cref{prop:epslowercond} we have
\[
r_{m,\eps}-d\geq r_{m,\eps}- \frac{3m^{2}}{k+1}\geq\biggl\lceil \frac{(m-k)^{3}}{3(m-k)-2}\biggr\rceil=\rmr_{\mathrm{gen}}^{m-k}.
\]
This implies, by \Cref{lem:first_existence}, the existence of a non-empty open Zariski set 
\[\calE_{m-k,r_{m,\eps}-d}\subset\bbC^{3(m-k)(r_{m,\eps}-d)}\] such that $\calO_{r_{m,\eps}-d}^{q_1'}=(\bbC^{m-k})^{\otimes 3}$ for every $q_1'\in \calE_{m-k,r_{m,\eps}-d}$. 
In particular, if we define the set
\[
\Omega_{m,\eps}'\coloneqq (G^{q_2})^{-1}\bigl(\calE_{m-k,r_{m,\eps}-d}\cap G^{q_2}(W_{m,\eps})\bigr)\subset W_{m,\eps},
\]
which is a non-empty Euclidean open subset, then we have that 
\[
\calO_{r_{m,\eps}-d}^{G^{q_2}(p_1')}=(\bbC^{m-k})^{\otimes 3}
\]
for every $p_1'\in\Omega_{m,\eps}'$.
In particular, if $p_1'=(x^{\alpha},y^{\alpha},z^{\alpha})_{\alpha=1,\dots,r_{m,\eps}}\in\Omega_{m,\eps}'$, then
\begin{align*}
\calO_{r_{m,\eps}-d}^{G^{q_2}(p_1')}&=\bigl\langle q_A\otimes \pi_{B}\bigl(y^{\alpha}+\Psi_B^{(p_1',q_2)}(\eta^{\alpha})\bigr)\otimes \pi_{C}\bigl(z^{\alpha}+\Psi_C^{(p_1',q_2)}(\tau^{\alpha})\bigr)\bigr\rangle_{q_A\in \mathbb{C}^{m-k}}^{\alpha \not \in \Gamma} \\[0.5ex]
 & \hphantom{{}={}}+\bigl\langle \pi_{A}\bigl(x^{\alpha}+\Psi_A^{(p_1',q_2)}(\xi^{\alpha})\bigr)\otimes q_B\otimes \pi_{C}\bigl(z^{\alpha}+\Psi_C^{(p_1',q_2)}(\tau^{\alpha})\bigr)\bigr\rangle_{q_B\in \mathbb{C}^{m-k}}^{\alpha \not \in \Gamma} \\[0.5ex]
 & \hphantom{{}={}}+\bigl\langle \pi_{A}\bigl(x^{\alpha}+\Psi_A^{(p_1',q_2)}(\xi^{\alpha})\bigr)\otimes \pi_{B}\bigl(y^{\alpha}+\Psi_B^{(p_1',q_2)}(\eta^{\alpha})\bigr)\otimes q_C\bigr\rangle^{\alpha \not \in \Gamma}_{q_C\in \mathbb{C}^{m-k}},\\
 &\subset (\pi_A\otimes \pi_B\otimes \pi_C)\bigl( \calQ'^{(p_1',q_2)}_{\psi_A,\psi_B,\psi_C}\bigr)=(\bbC^{m-k})^{\otimes 3}.\qedhere
\end{align*}
\end{proof}

In order to state the next result, we need to define another space depending on differential functions. For any differentiable maps $\psi_{A},\psi_{B},\psi_{C}\colon U_{m,\eps} \to \Mat_{m\times m}(\C)$ and for any $p\in U_{m,\eps}$, let $\calC_I^p$ be a complement of $\im\Psi_I^p$ in $\bbC^m$, that is, $ \mathbb{C}^{m} =\im\Psi_A^p\oplus \calC_{I}^p$, for each $I=A,B,C$. Then, we define the space
\[
\calR_{\psi_A,\psi_B,\psi_C}^{p,\calC_A^p,\calC_B^p,\calC_C^p}\coloneqq 
\bigl(\im\Psi_{A}^{p}\otimes \calC_{B}^p\otimes \calC_{C}^p\bigr)\oplus \bigl(\calC_{A}^p\otimes \im\Psi_{B}^{p}\otimes \calC_{C}^p\bigr)\oplus \bigl(\calC_{A}^p\otimes \calC_{B}^p\otimes\im\Psi_{C}^{p}\bigr)\oplus(\calC_A^p\otimes \calC_B^p\otimes \calC_{C}^p).
 \]
In the following lemma, we analyze a lower bound on the projection of the space $\calQ_{\psi_A,\psi_B,\psi_C}^p$ on $\calR_{\psi_A,\psi_B,\psi_C}^{p,\calC_A^p,\calC_B^p,\calC_C^p}$.
\begin{lemma}\label{lem:lower_bound_Q}
Let $k,m\in \mathbb{Z}_{>0}$ with $k\leq m-1$. Let $\varepsilon\in\bbR_{>0}$ and $\psi_A,\psi_B,\psi_C\colon U_{m,\eps}\to\Mat_{m\times m(\bbC)}$ be differentiable maps defined over a non-empty Euclidean open set $U_{m,\eps}\subset \calV_{m,\eps}$.
If $L_{k,m}^{-}\leq \varepsilon$, then for any $q_2\in \Pi_2(U_{m,\eps})$, there exists a non-empty Euclidean open set $\Omega_{m,\eps}'\subset \Omega_{m,\eps}^{q_2}(U_{m,\eps})$ such that
  \[
\dim\bigl(\pi_{\calR_p}(\calQ_{\psi_A,\psi_B,\psi_C}^p)\bigr)\geq m^3\Delta_{\psi_A,\psi_B,\psi_C}^p-3km^2
 \]
 for every $p_1'\in\Omega_{m,\eps}'$,  where $p=(p_1',q_2)$ and 
 \[\calR_p\coloneqq 
\calR_{\psi_A,\psi_B,\psi_C}^{p,\calC_A^p,\calC_B^p,\calC_C^p}\]
 for any complement $\calC_I^p$  of $\im\Psi_I^p$, for $I=A,B,C$, and $\pi_{\calR_p}\colon(\bbC^m)^{\otimes 3}\to\calR_{p}$ is the projection on~$\calR_p$.
\end{lemma}
\begin{proof}
By \Cref{lem:diff_geo}, for any $q_2\in\Pi_2(U_{m,\eps})$ there exist general linear maps $\pi_A,\pi_B,\pi_C\colon \bbC^m\to\bbC^{m-k}$ and a non-empty Euclidean open set $\Omega_{m,\eps}'\subset \Omega_{m,\eps}^{q_2}(U_{m,\eps})$ such that, for every $p_1'\in\Omega_{m,\eps}'$,
 \begin{equation}\label{formula:lem_Q_lower_bound}
(\pi_A\otimes \pi_B\otimes \pi_C)\bigl( \calQ'^{(p_1',q_2)}_{\psi_A,\psi_B,\psi_C}\bigr)=(\bbC^{m-k})^{\otimes 3}.
 \end{equation}
Let us consider any point $p=(p_1',q_2)\in U_{m,\eps}$ with $p_1'\in\Omega_{m,\eps}'$ and let us denote
\[
\calQ_p\coloneqq\calQ^{p}_{\psi_A,\psi_B,\psi_C},\qquad \calQ_p'\coloneqq\calQ'^{p}_{\psi_A,\psi_B,\psi_C}.
\]
We can suppose that the maps $\pi_A,\pi_B,\pi_C$
are surjective and that the images of $\im\Psi_I^{p}$ and $\calC_I^{p}$ have maximal dimension, that is, 
\[
\dim\bigl(\pi_I(\im\Psi_I^{p})\bigr)=\min\{\delta_{\psi_I}^pm,m-k\},\qquad \dim\bigl(\pi_I(\calC_I^{p})\bigr)=\min\{(1-\delta_{\psi_I}^p)m,m-k\},
\]
for each $I=A,B,C$. 
 We can consider a subspace $\calC_I'^p\subset \calC_I^p$ such that
\[
\calC_I'^p\simeq \pi_I(\calC_I'^p)=\pi_I(\calC_I^p)
\]
and the complement
$K_I^p\coloneqq \calC_I^p\cap\Ker\pi_I$,
which satisfies $\calC_I^p=\calC_I'^p\oplus K_I^p$. Analogously, we consider a subspace
$Z_I^p\subset\im\Psi_I^p$ such that
\[
Z_I^p\simeq \pi_I(Z_I^p)=\pi_I(\im\Psi_I^p)
\]
and the complement
$H_I^p\coloneqq \im\Psi_I^p\cap\Ker\pi_I$, 
so that $\im\Psi_I^p=H_I^p\oplus Z_I^p$.
Finally, by the surjectivity of $\pi_I$, there exists also a subspace $Z_I'\subset Z_I$ such that 
\begin{equation}\label{formula:m-c_direct_sum}
\bbC^{m-k}=
\im(\pi_I)=\pi_I(\calC_I^p)\oplus\pi_I(Z_I'^p)=\pi_I(\calC_I^p\oplus Z_I'^p),
\end{equation}
and hence a complement $Z_I''^p\subset Z_I^p$ such that $ Z_I^p=Z_I'^p\oplus Z_I''^p$. 
Observe that, if $\dim(\calC_I^p)\geq m-k$, then $Z_I'^p=0$.
Therefore, we have a decomposition of $\bbC^m$ for each $I=A,B,C$, given by
\[
\bbC^{m}=\calC_I^p\oplus \im\Psi_I^p=\calC_I'^p\oplus K_I^p\oplus H_I^p\oplus Z_I^p=\calC_I'^p\oplus K_I^p\oplus H_I^p\oplus Z_I'^p\oplus Z_I''^p.
\]
Observe that, if $\dim(\calC_I)\leq m-k$, then $K_I^p=0$.
Therefore, considering the space
\[
\calS_p'\coloneqq (\calC_A'^p\oplus Z_A'^p)\otimes (\calC_B'^p\oplus Z_B'^p)\otimes (\calC_C'^p\oplus Z_C'^p),
\] 
and the projection $\pi_{\calS'_p}\colon(\bbC^m)^{\otimes 3}\to \calS_p'$, 
we have $\bbC^m=(C_I'^p\oplus Z_I'^p)\oplus \Ker \pi_I$ for each $I=A,B,C$, and, in particular, $\calS_p' \simeq (\bbC^{m-k})^{\otimes 3}$.  Hence, by formula \eqref{formula:lem_Q_lower_bound}, 
\begin{equation}\label{formula:space_S}
\bigl((\pi_A\otimes\pi_B\otimes\pi_C)\circ \pi_{\calS_p'}\bigr)(\calQ_p')= (\pi_A\otimes\pi_B\otimes\pi_C)\bigl( \pi_{\calS_p'}(\calQ_p') \bigr) = (\pi_A\otimes\pi_B\otimes\pi_C)(\calQ_p')=(\bbC^{m-k})^{\otimes 3}.
\end{equation}
Now, let us define the subspace $\calR_p'\subset\calR_p$ as
\[
\calR_p'\coloneqq (Z_A'^p\otimes \calC_B'^p\otimes \calC_C'^p)\oplus (\calC_A'^p\otimes Z_B'^p\otimes \calC_C'^p)\oplus (\calC_A'^p\otimes \calC_B'^p\otimes Z_C'^p)\oplus (\calC_A'^p\otimes \calC_B'^p\otimes \calC_C'^p).
\]
Then, applying the map $(\pi_A\otimes \pi_B\otimes \pi_C)\colon (\bbC^m)^{\otimes 3}\to (\bbC^{m-k})^{\otimes 3}$, we obtain the space
\begin{align*}
(\pi_A\otimes\pi_B\otimes\pi_C)(\calR_p') &\coloneqq \bigl(\pi_A(Z_A')\otimes \pi_B(\calC_{B}')\otimes \pi_C(\calC_{C}')\bigr)\oplus \bigl(\pi_A(\calC_{A}')\otimes \pi_B(Z_B')\otimes \pi_C(\calC_{C}')\bigr)\\
&\hphantom{{}\coloneqq{}}\oplus \bigl(\pi_A(\calC_{A}')\otimes \pi_B(\calC_{B}')\otimes \pi_C(Z_C')\bigr)\oplus \bigl(\pi_A(\calC_{A}')\otimes \pi_B(\calC_{B}')\otimes \pi_C(\calC_{C}')\bigr).
\end{align*}
Since the map $\pi_A\otimes\pi_B\otimes\pi_C$ is an isomorphism when restricted to $\calS_p'$, considering the two projections
\[\pi_{\calR_p}\colon(\bbC^m)^{\otimes 3}\to\calR_p,\qquad \pi_{\calR_p'}\colon(\bbC^m)^{\otimes 3}\to\calR_p',\]
we have that the following diagram commutes:
\begin{equation}\label{formula:diagram_commuting}
\begin{tikzcd}
\arrow{d}{\pi_{\calS_p'}} (\bbC^{m})^{\otimes 3}   \arrow{r}{\pi_{\calR_p}}  & \calR_p\arrow{d}{\pi_{\calR_p'}}\\
\arrow[d,"\pi_A\otimes\pi_B\otimes\pi_C","\rotatebox{90}{\,\,$\widesim[2]{}$}"'] \calS_p'   \arrow{r}{\pi_{\calR_p'}}  & \calR_p'\arrow{d}{\pi_A\otimes\pi_B\otimes\pi_C}\\
 (\pi_A\otimes\pi_B\otimes\pi_C)\bigl((\bbC^{m})^{\otimes 3}\bigr)\arrow{r} {\pi_{\calR_p'\vphantom{'}}} &  (\pi_A\otimes\pi_B\otimes\pi_C)(\calR_p')
\end{tikzcd}
\end{equation}
It is easy to observe that $\pi_{\calR_p}(\calQ_p)=\pi_{\calR_p}(\calQ_p')$.
Thus, using equality \eqref{formula:space_S} and diagram \eqref{formula:diagram_commuting}
we have 
\begin{align}\label{formula:lower_bound_gap}
\nonumber \dim \bigl( \pi_{\calR_p}(\calQ_p)\bigr)  &=\dim \bigl( \pi_{\calR_p}(\calQ_p')\bigr)\geq \dim \bigl( \pi_{\calR_p'}(\calQ_p')\bigr) =\dim\Bigl(\pi_{\calR_p'}\bigl( \pi_{\calS_p'}(\calQ_p')\bigr)\Bigr)\\
\nonumber &\geq\dim \biggl((\pi_{A} \otimes \pi_{B} \otimes \pi_{C})\Bigl(\pi_{\calR_p'}\bigl( \pi_{\calS_p'}(\calQ_p')\bigr)\Bigr)\biggr)=\dim \biggl(\pi_{\calR_p'}\Bigl((\pi_{A} \otimes \pi_{B} \otimes \pi_{C})\bigl( \pi_{\calS_p'}(\calQ_p')\bigr)\Bigr)\biggr)\\ 
 &=\dim \Bigl(\pi_{\calR_p'}\bigl((\bbC^{m-k})^{\otimes 3}\bigr)\Bigr) =\dim \bigl((\pi_{A} \otimes \pi_{B} \otimes \pi_{C}) (\calR_p')\bigr)=\dim(\calR_p').
\end{align}
So, we only need to give an estimation of the dimension of the space $\calR_p'$. In order to do this, we need to distinguish between all the cases in which the dimension $\calC_I^p$ is lower than $m-k$ and those in which it is equal to $m-k$, for each $I=A,B,C$.
In the case where $\dim \calC_I'^p= m-k$ for every $I=A,B,C$, then we also have $Z_I'^p=0$ for every $I=A,B,C$.
We have
\[
\calR_p'=\calC_A'^p\otimes \calC_B'^p\otimes \calC_C'^p=(\bbC^{m-k})^{\otimes 3},
\]
which implies that
\[
\dim \calR_p'=(m-k)^3\geq m^3\Delta_{\psi_A,\psi_B,\psi_C}^p-3km^2.
\]
In the case where
\[
\dim \calC_A'^p< m-k,\qquad \dim \calC_B'^p=\dim \calC_C'^p=m-k,
\]
then we have
\[
\calR_p'=(Z_A'^p\otimes \calC_B'^p\otimes \calC_C'^p)\oplus (\calC_A'^p\otimes \calC_B'^p\otimes \calC_C'^p)=(Z_A'^p\oplus \calC_A'^p)\otimes \calC_B'^p\otimes \calC_C'^p=(\bbC^{m-k})^{\otimes 3},
\]
which gives the same lower bound as in the previous case. 
In the case where
\[
\dim \calC_A'^p< m-k,\quad \dim \calC_B'^p< m-k,\quad \dim \calC_C'^p=m-k,
\]
then we have that $Z_C'^p=0$, $\calC_A'^p=\calC_A^p$, and $\calC_B'^p=\calC_B^p$ so that
\[
\calR_p'= (Z_A'^p\otimes \calC_B^p\otimes \calC_C'^p)\oplus (\calC_A^p\otimes Z_B'^p\otimes \calC_C'^p)\oplus (\calC_A^p\otimes \calC_B^p\otimes \calC_C'^p).
\]
In particular, we have \[\dim \calC_I^p=(1-\delta_{\psi_I}^p)m,\qquad\dim Z_I'^p=\delta_{\psi_I}^pm-k\] for each $I=A,B$, 
which imply that
\begin{align*}
 &\dim (\calR_p')  = (m-k)\bigl(m(1-\delta_{\psi_A}^p)(\delta_{\psi_B}^pm-k)+m(1-\delta_{\psi_B}^p)(\delta_{\psi_A}^pm-k)+m^2(1-\delta_{\psi_A}^p)(1-\delta_{\psi_B}^p)\bigr)\\
 &\quad=(m-k)\Bigl(m^2\bigl(\delta_{\psi_A}^p(1-\delta_{\psi_B}^p)+\delta_{\psi_B}^p(1-\delta_{\psi_A}^p)+(1-\delta_{\psi_A}^p)(1-\delta_{\psi_B}^p)\bigr)-km\bigl((1-\delta_{\psi_A}^p)+(1-\delta_{\psi_B}^p)\bigr)\Bigr)\\
 &\quad=m^3(1-\delta_{\psi_A}^p\delta_{\psi_B}^p)-km^2\bigl(3-\delta_{\psi_A}^p\delta_{\psi_B}^p - \delta_{\psi_A}^p - \delta_{\psi_C}^p\bigr)+k^2m(2-\delta_{\psi_A}^p-\delta_{\psi_B}^p)\geq m^3\Delta_{\psi_A,\psi_B,\psi_C}^p-3km^2.
\end{align*}
Finally, in the case where
$\dim \calC_I'^p< m-k$ for every $I=A,B,C$,
we have $\calC_I'^p=\calC_I^p$, and hence, since
\[\dim \calC_I^p=(1-\delta_{\psi_I}^p)m,\qquad\dim Z_I'^p=\delta_{\psi_I}^pm-k\]
for each $I=A,B,C$, we have
\begin{align*}
\dim (\calR_p')&= m^2(\delta_{\psi_A}^pm-k)(1-\delta_{\psi_B}^p)(1-\delta_{\psi_C}^p)+m^2(1-\delta_{\psi_A}^p)(\delta_{\psi_B}^pm-k)(1-\delta_{\psi_C}^p)\\
&\hphantom{{}={}}+m^2(1-\delta_{\psi_A}^p)(1-\delta_{\psi_B}^p)(\delta_{\psi_C}^pm-k)+m^3(1-\delta_{\psi_A}^p)(1-\delta_{\psi_B}^p)(1-\delta_{\psi_C}^p)\\
&=m^3\Delta_{\psi_A,\psi_B,\psi_C}^p-km^2\bigl(3-\delta_{\psi_A}^p(2-\delta_{\psi_B}^p)-\delta_{\psi_B}^p(2-\delta_{\psi_C}^p)-\delta_{\psi_C}^p(2-\delta_{\psi_A}^p)\bigr)\\
&\geq m^3\Delta_{\psi_A,\psi_B,\psi_C}^p-3km^2.
\end{align*}
Therefore, by formula \eqref{formula:lower_bound_gap}, we obtain the lower bound 
 \[
\dim\bigl(\pi_{\calR_p}(\calQ_{\psi_A,\psi_B,\psi_C}^p)\bigr)\geq \dim (\calR_p')\geq m^3\Delta_{\psi_A,\psi_B,\psi_C}^p-3km^2,
 \]
 which concludes the proof. 
\end{proof}
\begin{remark}\label{rem:inaccuracy_Lan}
    In \cite[Section 2.5.3]{Lan19}, a similar lower bound is analyzed. 
Set $r_m\coloneqq (11/30)m^2$. For any point $p=(x^{\alpha},y^{\alpha},z^{\alpha},\xi^{\alpha},\eta^{\alpha},\tau^{\alpha})_{\alpha=1,\dots,r_m}\in\bbC^{6mr_m}$ and some fixed $\psi_A,\psi_B,\psi_C\colon \bbC^{3mr_m}\to\Mat_{m\times m}(\bbC)$, the function $f\colon\bbC^{6mr_m}\to(\bbC^{m})^{\otimes 3}$ is defined by
\[f(p)\coloneqq T_{\Psi_A^p,\Psi_B^p,\Psi_C^p}^1(p).\] 
With an abuse of notation, up to reordering the summands, the image of a point $p'=(0,0,0,\xi'^{\alpha},\eta'^{\alpha},\tau'^{\alpha})$ via the differential $\dd f_p$ of $f$ at $p$ is described in \cite[p.~35]{Lan19} as
\begin{align*}
\dd f_p(p')&=\sum_{\alpha=1}^{r_m}\Bigl(\bigl(\Psi_A^p(\xi_{\alpha}')+\Psi_A'^{p'}(\xi_{\alpha})\bigr)\otimes y^{\alpha}\otimes z^{\alpha}+x^{\alpha}\otimes\bigl(\Psi_B^p(\xi_{\alpha}')+\Psi_B'^{p'}(\xi_{\alpha})\bigr)\otimes z^{\alpha}\\
&\hphantom{{}={}}\qquad+x^{\alpha}\otimes y^{\alpha}\otimes \bigl(\Psi_C^p(\xi_{\alpha}')+\Psi_C'^{p'}(\xi_{\alpha})\bigr)\Bigr)
\end{align*}
for some linear maps $\psi_A',\psi_B',\psi_C'\colon\bbC^{3mr_m}\to\Mat_{m\times m}(\bbC)$. In particular, it is stated that, whenever $\sigma\coloneqq \delta_{\psi_A}^p+\delta_{\psi_B}^p+\delta_{\psi_C}^p>11/20$,
\[
\dim \Bigl(\dd f_p\bigl(\Pi_2(\bbC^{6mr_m})\cap\Ker\psi_A'\cap\Ker\psi_B'\cap \Ker\psi_C'\bigr)\Bigr)=\sigma m^3-3m^2>\frac {11}{20}m^3,
\]
for a sufficiently large $m$. However, {it is not so, as we now explain.} If we consider any point belonging to the above intersection $p'\in \Pi_2(\bbC^{6mr_m})\cap\Ker\psi_A'\cap\Ker\psi_B'\cap \Ker\psi_C'$,
we have
\[
\dd f_p(p')=\sum_{\alpha=1}^{r_m}\Psi_A^p(\xi_{\alpha}')\otimes y^{\alpha}\otimes z^{\alpha}+x^{\alpha}\otimes\Psi_B^p(\xi_{\alpha}')\otimes z^{\alpha}+x^{\alpha}\otimes y^{\alpha}\otimes \Psi_C^p(\xi_{\alpha}').
\]
Therefore, since $\rk(\Psi_I^p)=\delta_{\psi_I}^p$ for each $I=A,B,C$, we have
\[
\dim\bigl(\im \dd f_p\bigr)\leq \delta^p_{\psi_A}mr_m+\delta^p_{\psi_A}mr_m+\delta^p_{\psi_A}mr_m=\sigma mr_m=\frac{11}{30}\sigma m^3,
\]
which is, for sufficiently large $m$, less than $\sigma m^3-3m^2$.
\end{remark} 
\begin{lemma}\label{projectionQ} 
Let $k,m\in \mathbb{Z}_{>0}$ with $k\leq m-1$.  Let $\varepsilon\in\bbR_{>0}$ and $\psi_A,\psi_B,\psi_C\colon U_{m,\eps}\to\Mat_{m\times m(\bbC)}$ be differentiable maps defined over a non-empty Euclidean open set $U_{m,\eps}\subset \calV_{m,\eps}$.
If $L_{k,m}^{-}\leq \varepsilon$, then for any $q_2\in \Pi_2(U_{m,\eps})$, there exists a non-empty Zariski open set $\calE_{m,\eps}'\subset \Omega_{m,\eps}^{q_2}(U_{m,\eps})$ such that 
 \[\pi_{\calC}(\calQ_{\psi_a,\psi_B,\psi_C}^p)=\calC_A^p\otimes \calC_B^p\otimes \calC_{C}^p\]
 for every $p_1' \in \calE_{m,\eps}'$, where $p = (p_1', q_2)$, and any complement $\calC_I^p$  of $\im\Psi_I^p$, for each $I=A,B,C$, where
 \[\pi_{\calC}: (\mathbb{C}^m)^{\otimes 3} \to \calC_A^p \otimes \calC_B^p \otimes \calC_C^p\]
 is the projection on $\calC_A^p\otimes \calC_B^p\otimes \calC_{C}^p$.
\end{lemma}
\begin{proof}
Since $L_{k,m}^-\leq \eps$, it follows by \Cref{prop:secondepslowercond} that
$r_{m,\eps}\geq\rmr_{\mathrm{gen}}^m$.
Therefore, by 
\Cref{lem:first_existence}, there exists a non-empty Zariski open set $\calE_{m,r_{m,\eps}}\subset\bbC^{3mr_{m,\eps}}$, such that $\calO_{r_{m,\eps}}^{p_1}=(\bbC^m)^{\otimes 3}$ for every $p_1\in\calE_{m,r_{m, \eps}}$. 
The set 
\[\calE_{m,\eps}'\coloneqq \calE_{m,r_{m,\eps}}\cap \Omega_{m,\eps}^{q_2}(U_{m,\eps})\] is a non-empty Zariski open set in $\Omega_{m,\eps}^{q_2}(U_{m,\eps})$. Moreover, since by definition we have
\[
\pi_{\calC}\bigl(\calQ_{\psi_a,\psi_B,\psi_C}^{(p_1,q_2)}\bigr)=\pi_{\calC}(\calO_{r_{m,\eps}}^{p_1})=\pi_{\calC}\bigl((\bbC^m)^{\otimes 3}\bigr),
\]
for any $p_1\in\calE_{m,\eps}'$, the statement follows.
\end{proof}

We are finally ready to provide the last proof of this section.
\begin{proof}[Proof of \Cref{df}]
Let $\calA_{m,\eps}\subset \bbC^{3mr_{m,\eps}}$ be the set of  points $p_1=(x^{\alpha},y^{\alpha},z^{\alpha})_{\alpha=1,\dots,r_{m,\eps}}\in \bbC^{3mr_{m,\eps}}$ such that
 \[
\dim(\langle y^{\alpha}\otimes z^{\alpha}\rangle_{\alpha=1,\dots,r_{m,\eps}})=\dim(\langle x^{\alpha}\otimes z^{\alpha}\rangle_{\alpha=1,\dots,r_{m,\eps}})=\dim(\langle x^{\alpha}\otimes y^{\alpha}\rangle_{\alpha=1,\dots,r_{m,\eps}})=\min\{r_{m,\eps},m^2\}.
 \]
 Clearly, $\calA_{m,\eps}$ is a non-empty Zariski open set and  the set $\calA_{m,\eps}'\coloneqq \calA_{m,\eps}\cap\Omega_{m,\eps}^{q_2}(U_{m, \eps})$  is a non-empty Zariski open set of $\Omega_{m,\eps}^{q_2}(U_{m, \eps})$, for any $q_2\in \Pi_2(U_{m,\eps})$. Now, for any $p_1'\in\calA_{m,\eps}'$, denoting $p=(p_1',q_2)$, 
 let $\calC_I^p$ be a complement of $\im\Psi_I^p$ for each $I=A,B,C$.  We investigate the intersections of $\calP_{\psi_A,\psi_B,\psi_C}^p$ with the spaces \[D_A^p\coloneqq \calC_{A}^p\otimes \im\Psi_{B}^p\otimes \im\Psi_{C}^p,\quad D_B^p\coloneqq\im\Psi_{A}^p\otimes \calC_{B}^p\otimes \im\Psi_{C}^p,\quad D_C^p\coloneqq\im\Psi_{A}^p\otimes \im\Psi_{B}^p\otimes \calC_{C}^p,\]
respectively.
 Let us consider the spaces \[
E_{A}^p\coloneqq\bigl\langle w_{A}\otimes y^{\alpha}\otimes z^{\alpha}\bigr\rangle^{\alpha=1,\dots,r_{m,\eps}}_{w_{A}\in \im(\Psi_{A}^p)},\quad E_{B}^p\coloneqq\bigl\langle x^{\alpha}\otimes w_{B}\otimes z^{\alpha}\bigr\rangle^{\alpha=1,\dots,r_{m,\eps}}_{w_{B}\in \im(\Psi_{B}^p)},\quad E_{C}^p\coloneqq\bigl\langle x^{\alpha}\otimes y^{\alpha}\otimes w_{C}\bigr\rangle^{\alpha=1,\dots,r_{m,\eps}}_{w_{C}\in \im(\Psi_{C}^p)},\]
and the spaces
\[
T_A^p\coloneqq\im\Psi_{A}^p\otimes \mathbb{C}^{m}\otimes \mathbb{C}^{m},\quad T_B^p\coloneqq \mathbb{C}^{m}\otimes\im\Psi_{B}^p\otimes \mathbb{C}^{m},\quad T_C^p\coloneqq\mathbb{C}^{m}\otimes \mathbb{C}^{m}\otimes \im\Psi_{C}^p.
\]
Now, since $E_A^p+D_C^p\subset T_A^p$ and 
\[
\dim E_{A}^p=m\delta_{\psi_A}^p\cdot \min\{r_{m,\eps},m^2\},
\]
then, we have 
\begin{align}\label{formula:first_lower_bound_first_inequality}
\nonumber\dim (D_C^p\cap E_{A}^p)&=\dim D_C^p+\dim E_{A}^p-\dim (D_C^p+ E_{A}^p)\geq \dim D_C^p+\dim E_{A}^p-\dim T_A^p\\[1ex]
\nonumber&=m^3\delta_{\psi_A}^p\delta_{\psi_B}^p(1-\delta_{\psi_C}^p)+m\delta_{\psi_A}^p\cdot\min\biggl\{\biggl\lceil\frac{1+\eps}{3}m^2\biggr\rceil,m^2\biggr\}-m^3\delta_{\psi_A}^p\\
&\geq m^{3}\delta_{\psi_A}^p\left( \delta_{\psi_B}^p(1-\delta_{\psi_C}^p)+\min\biggl\{\frac{1+\eps}{3},1\biggr\}-1 \right)\geq m^{3}\delta_{\psi_A}^p\left( \delta_{\psi_B}^p(1-\delta_{\psi_C}^p)- \frac{2}{3} \right).
\end{align} 
Proceeding analogously with the spaces $E_{B}^p$ and $T_{B}^p$, we obtain  the inequality
\begin{equation}\label{formula:first_lower_bound_second_inequality}
\dim (D_C^p\cap E_{B}^p)\geq m^{3}\delta_{\psi_B}^p\left( \delta_{\psi_A}^p(1-\delta_{\psi_C}^p)- \frac{2}{3} \right).
\end{equation}
By using inequalities \eqref{formula:first_lower_bound_first_inequality} and \eqref{formula:first_lower_bound_second_inequality} together, we have
\begin{align}\label{formula:inequality1}
\nonumber&\dim \bigl(\calP_{\psi_A,\psi_B,\psi_C}^p\cap D_{C}^p\bigr)   \geq \max\bigl\{\dim (D_C^p \cap E_{A}^p),\dim (D_C^p \cap E_B^p)\bigr\}\\[1ex]
\nonumber&\quad\geq m^{3} \max \left\{\delta_{\psi_A}^p\left( \delta_{\psi_B}^p(1-\delta_{\psi_C}^p)- \frac{2}{3} \right),\delta_{\psi_B}^p\left( \delta_{\psi_A}^p(1-\delta_{\psi_C}^p)- \frac{2}{3} \right) \right\} \\
&\quad=m^{3}\left( \delta_{\psi_A}^p\delta_{\psi_B}^p(1-\delta_{\psi_C}^p)- \frac{2}{3}\min \{\delta_{\psi_A}^p,\delta_{\psi_B}^p\} \right).
\end{align}
By permuting the roles of $A,B,C$, 
we get the analogous inequalities
\begin{align}\label{formula:inequality2}
&\dim \bigl(\calP_{\psi_A,\psi_B,\psi_C}^p\cap D_{B}^p\bigr)   \geq m^{3}\biggl( \delta_{\psi_A}^p\delta_{\psi_C}^p(1-\delta_{\psi_B}^p)- \frac{2}{3}\min \{\delta_{\psi_A}^p,\delta_{\psi_C}^p\} \biggr)\\
\label{formula:inequality3}&\dim \bigl(\calP_{\psi_A,\psi_B,\psi_C}^p\cap D_{A}^p\bigr)   \geq m^{3}\biggl( \delta_{\psi_B}^p\delta_{\psi_C}^p(1-\delta_{\psi_A}^p)- \frac{2}{3}\min \{\delta_{\psi_B}^p,\delta_{\psi_C}^p\} \biggr)
\end{align}
Now, by \Cref{lem:lower_bound_Q}, there exists a non-empty Euclidean open set $\Omega_{m,\eps}'\subset \Omega_{m,\eps}^{q_2}(U_{m,\eps})$ such that, if $p_1'\in \Omega_{m,\eps}'$,
\[
\dim\bigl(\pi_{\calR_p}(\calQ_{\psi_A,\psi_B,\psi_C}^p)\bigr)\geq m^3\Delta_{\psi_A,\psi_B,\psi_C}^p-3km^2.
\]
If we consider the restriction of the projection
$\pi_{\calR_p}\colon \calP_{\psi_A,\psi_B,\psi_C}^p+\calQ_{\psi_A,\psi_B,\psi_C}^p\to \calR_{\psi_A,\psi_B,\psi_C}^{p,\calC_A^p,\calC_B^p,\calC_C^p}$, it is straightforward to observe that
\begin{gather*}
(\calP_{\psi_A,\psi_B,\psi_C}^p\cap D_A^p)\oplus (\calP_{\psi_A,\psi_B,\psi_C}^p\cap D_B^p)\oplus (\calP_{\psi_A,\psi_B,\psi_C}^p\cap D_C^p)\subset\Ker(\pi_{\calR_p}),\\[0.5ex]
\pi_{\calR_p}(\calQ_{\psi_A,\psi_B,\psi_C}^p)\subset\pi_{\calR_p}(\calP_{\psi_A,\psi_B,\psi_C}^p+\calQ_{\psi_A,\psi_B,\psi_C}^p).
\end{gather*}
Therefore, if $p_1'\in\Omega_{m,\eps}'\cap\calA_{m,\eps}'$,
by inequalities \eqref{formula:inequality1}, \eqref{formula:inequality2}, and \eqref{formula:inequality3}, we obtain
\begin{align}\label{formula:first_bound}
\nonumber&\dim (\calP_{\psi_A,\psi_B,\psi_C}^p+\calQ_{\psi_A,\psi_B,\psi_C}^p) \geq  \dim (\calP_{\psi_A,\psi_B,\psi_C}^p\cap D_A^p)+\dim (\calP_{\psi_A,\psi_B,\psi_C}^p\cap D_B^p)\\[2ex]
\nonumber&\quad \hphantom{{}={}}+\dim (\calP_{\psi_A,\psi_B,\psi_C}^p\cap D_C^p) + \dim\bigl(\pi_{\calR_p}(\calQ_{\psi_A,\psi_B,\psi_C}^{p})\bigr) \\[1ex] 
\nonumber&\quad\geq   m^{3}\biggl( \max \biggl\{ 0, \delta_{\psi_A}^p\delta_{\psi_B}^p(1-\delta_{\psi_C}^p)-\frac{2}{3}\min \{\delta_{\psi_A}^p,\delta_{\psi_B}^p\}\biggr\}+ \max \biggl\{ 0, \delta_{\psi_A}^p\delta_{\psi_C}^p(1-\delta_{\psi_B}^p)-\frac{2}{3}\min \{\delta_{\psi_A}^p,\delta_{\psi_C}^p\}\biggr\}  \\[1ex]
\nonumber&\quad \hphantom{{}={}}+\max \left\{ 0, \delta_{\psi_B}^p\delta_{\psi_C}^p(1-\delta_{\psi_A}^p)-\frac{2}{3}\min \{\delta_{\psi_B}^p,\delta_{\psi_C}^p\}\right\} \biggr)+m^3\Delta_{\psi_A,\psi_B,\psi_C}^p-3m^2k\\[1ex]
&\quad=   m^{3}\mu_{1}(\delta_{\psi_A}^p,\delta_{\psi_B}^p,\delta_{\psi_C}^p)-3m^{2}k.
\end{align}
Considering any basis $\{v_{I,j_I}\}_{j_I=\delta_{\psi_I}^pm+1,\dots,m}$ of $\calC_I^p$, for any $I=A,B,C$, we have 
\[\Bigl\langle\overline{v_{A,j_{A}}\otimes y^{\alpha}\otimes z^{\alpha}},\overline{x^{\alpha}\otimes v_{B,j_{B}}\otimes z^{\alpha}}, \overline{x^{\alpha}\otimes y^{\alpha}\otimes v_{C,j_{C}}}\Bigr\rangle_{j_A,j_B,j_C}^{\alpha=1,\dots,r_{m,\eps}}=(\mathbb{C}^{m})^{\otimes 3} /\calP_{\psi_a,\psi_B,\psi_C}^p.\] 
Therefore, we get
\[
\dim\bigl((\mathbb{C}^{m})^{\otimes 3} /\calP_{\psi_a,\psi_B,\psi_C}^p\bigr)\leq mr_{m,\eps}(3-\delta_{\psi_A}^p-\delta_{\psi_B}^p-\delta_{\psi_C}^p),\]
which is equivalent to 
\begin{equation}\label{formula:inequality4}
\dim \calP_{\psi_a,\psi_B,\psi_C}^p\geq m^{3}- mr_{m,\eps}(3-\delta_{\psi_A}^p-\delta_{\psi_B}^p-\delta_{\psi_C}^p)\geq m^{3}\biggl(1- \frac{1+\varepsilon}{3}(3-\delta_{\psi_A}^p-\delta_{\psi_B}^p-\delta_{\psi_C}^p)\biggr).
\end{equation}
Defining the space
\[
\calL_p\coloneqq D_A^p\oplus D_B^p\oplus D_C^p\oplus \bigl(\im\Psi_{A}^{p}\otimes \calC_{B}^p\otimes \calC_{C}^p\bigr)\oplus \bigl(\calC_{A}^p\otimes \im\Psi_{B}^{p}\otimes \calC_{C}^p\bigr)\oplus \bigl(\calC_{A}^p\otimes \calC_{B}^p\otimes\im\Psi_{C}^{p}\bigr),
\]
we have $\calP_{\psi_a,\psi_B,\psi_C}^p\subseteq \calL_p$. By \Cref{projectionQ}, there exists a non-empty Zariski open set $\calE_{m,\eps}'\subset\Omega_{m,\eps}^{q_2}(U_{m,\eps})$, such that, if $p_1'\in\calE_{m,\eps}'$, then
\[
\pi_{\calC}(\calQ_{\psi_A,\psi_B,\psi_C}^p)=\calC_A^p\otimes\calC_B^p\otimes \calC_C^p.
\]
If we consider the restriction of the projection 
$\pi_{\calC}\colon \calP_{\psi_A,\psi_B,\psi_C}^p+\calQ_{\psi_A,\psi_B,\psi_C}^p\to \calC_A^p\otimes\calC_B^p\otimes \calC_C^p$, it is straightforward to observe that
\[
(\calP_{\psi_A,\psi_B,\psi_C}^p\cap \calL_p)\subset\Ker(\pi_{\calC}),\qquad
\pi_{\calC}(\calQ_{\psi_A,\psi_B,\psi_C}^p)\subset\pi_{\calC}(\calP_{\psi_A,\psi_B,\psi_C}^p+\calQ_{\psi_A,\psi_B,\psi_C}^p).
\]
Thus, if $p_1'\in\calE_{m,\eps}'$, by inequality \eqref{formula:inequality4}, we have
\begin{align}\label{formula:second_bound}
\nonumber&\dim (\calP_{\psi_a,\psi_B,\psi_C}^p+\calQ_{\psi_a,\psi_B,\psi_C}^p)  \geq \dim (\calP_{\psi_a,\psi_B,\psi_C}^p\cap \calL_p) +\dim \pi_{\calC}(\calQ_{\psi_a,\psi_B,\psi_C}^p)\\[2ex]
\nonumber&\qquad= \dim (\calP_{\psi_a,\psi_B,\psi_C}^p)+\dim( \calC_A^p\otimes\calC_B^p\otimes\calC_C^p)\\[1ex]
 \nonumber&\qquad \geq m^{3}- \frac{1+\varepsilon}{3}m^{3}(3-\delta_{\psi_A}^p-\delta_{\psi_B}^p-\delta_{\psi_C}^p)+m^3(1-\delta_{\psi_A}^p)(1-\delta_{\psi_B}^p)(1-\delta_{\psi_C}^p) \\
 \nonumber&\qquad = m^{3}\biggl( \frac{\delta_{\psi_A}^p+\delta_{\psi_B}^p+\delta_{\psi_C}^p}{3}+(1-\delta_{\psi_A}^p)(1-\delta_{\psi_B}^p)(1-\delta_{\psi_C}^p) \biggr)+m^{3}\varepsilon\biggl(\frac{\delta_{\psi_A}^p+\delta_{\psi_B}^p+\delta_{\psi_C}^p}{3}-1\biggr) \\
	 \nonumber&\qquad \geq  m^{3}\biggl( \frac{\delta_{\psi_A}^p+\delta_{\psi_B}^p+\delta_{\psi_C}^p}{3}+(1-\delta_{\psi_A}^p)(1-\delta_{\psi_B}^p)(1-\delta_{\psi_C}^p) \biggr)-m^{3}\varepsilon\\[1ex]
     &\qquad\geq m^{3}\cdot \mu_{2}(\delta_{\psi_A}^p,\delta_{\psi_B}^p,\delta_{\psi_C}^p)-m^{3}\varepsilon.
\end{align}
Now, since $\Omega_{m,\eps}'\cap\calA_{m,\eps}'$ is a non-empty Euclidean open set of $\Omega_{m,\eps}^{q_2}(U_{m,\eps})$ and $\calE_{m,\eps}'$ is a non-empty Zariski open set of $\Omega_{m,\eps}^{q_2}(U_{m,\eps})$, then also the intersection $\Omega_{m,\eps}'\cap\calA_{m,\eps}'\cap \calE_{m,\eps}'$ is a non-empty Euclidean open set of $\Omega_{m,\eps}^{q_2}(U_{m,\eps})$. In particular, there exists a point $\bar{p}_1'\in\Omega_{m,\eps}'\cap\calA_{m,\eps}'\cap \calE_{m,\eps}'$ such that, setting $\bar{p}=(\bar{p}_1',q_2)$, by inequalities \eqref{formula:first_bound} and \eqref{formula:second_bound} together,
we have 
\begin{align*}
\dim \bigl(\calP_{\psi_A,\psi_B,\psi_C}^{\bar{p}}+\calQ_{\psi_A,\psi_B,\psi_C}^{\bar{p}}\bigr) & \geq \max \{m^{3}  \mu_{1}(\delta_{\psi_A}^{\bar{p}},\delta_{\psi_B}^{\bar{p}},\delta_{\psi_C}^{\bar{p}})-3km^{2},m^{3} \mu_{2}(\delta_{\psi_A}^{\bar{p}},\delta_{\psi_B}^{\bar{p}},\delta_{\psi_C}^{\bar{p}})-m^3\varepsilon\} \\
 & \geq m^{3}\cdot \max \{\mu_{1}(\delta_{\psi_A}^{\bar{p}},\delta_{\psi_B}^{\bar{p}},\delta_{\psi_C}^{\bar{p}}),\mu_{2}(\delta_{\psi_A}^{\bar{p}},\delta_{\psi_B}^{\bar{p}},\delta_{\psi_C}^{\bar{p}})\} -\max \{3km^{2},m^{3}\varepsilon\} \\
 & \geq m^2\bigl(m\mu-\max \{3k,m\varepsilon\}\bigr),
\end{align*}
where the last inequality follows from the definition of $\mu$, given in formula \eqref{LBnotationmin}.
\end{proof}

\section*{Appendix A}
 \label{SecA} \addcontentsline{toc}{section}{Appendix~A}
  \stepcounter{section}
  \renewcommand{\thesection}{A}

\begin{proposition}\label{prop: mu>1/2}
 The inequality $2\mu-1 >0$ holds. 
\end{proposition}
\begin{proof}
We prove the required 
inequality by showing that for every $\alpha,\beta,\gamma\in[0,1]\subset\bbR$, we have
\begin{equation}\label{formula:two_inequalities}
\mu_i(\alpha,\beta,\gamma) > \frac12
\end{equation}
either for $i=1$ or for $i=2$.
First, notice that  for all the six points in the set \[S \coloneqq \left\{\biggl(0,\frac12,1\biggr),\biggl(0,1,\frac12\biggr),\biggl(\frac12,1,0\biggr), \biggl(1,\frac12, 0\biggr),\biggl(\frac12,0,1\biggr),\biggl(1,0,\frac12\biggr)\right\}\] the second inequality holds. 
Thus, it is enough to show the claim for $(\alpha, \beta, \gamma) \in [0,1]^3 \setminus S$. In fact, for any such triple the statement holds even for a weaker version of inequality \eqref{formula:two_inequalities} for $i=2$, namely,
 \begin{equation}\label{LB2}
(1-\alpha)(1-\beta)(1-\gamma) +\alpha(1-\beta)(1-\gamma) + \beta(1-\alpha)(1-\gamma)+\gamma(1-\alpha)(1-\beta) > \frac{1}{2}.
\end{equation}
Let us rewrite inequalities \eqref{formula:two_inequalities} for $i=1$ and \eqref{LB2} by {setting} $a \coloneqq 1-\alpha$, $b \coloneqq 1 - \beta$, and $c \coloneqq 1 - \gamma$.
Inequality \eqref{formula:two_inequalities} for $i=1$ is then
    \begin{equation}\label{first_ine}
    \frac{a+b+c}{3} - abc - \frac12 < 0,
    \end{equation}
    while inequality \eqref{LB2} corresponds to
    \begin{equation}\label{second_ine}
    ab + bc + ac - 2abc - \frac12 > 0.
    \end{equation}
Since both of them are symmetric in $a,b,c$, we may assume $a \geq b \geq c$. Then, if
    $ab > 1/2$ we have \[ab + bc + ca - 2abc - \frac12 > c(b + a - 2ab) > c(2\sqrt{ab} - 2ab) \geq 0,\] where the last inequality follows from $ab\leq 1$, proving \eqref{second_ine}. So, let us suppose
    that $ab \leq 1/2$. We claim that the \eqref{first_ine} holds.
To prove the latter we fix $b$ and $c$ and consider a univariate function $f_{b,c}\colon[b,1]\to\bbR$, consisting on the first member of inequality \eqref{first_ine}, that is, $$f_{b,c}(a) \coloneqq \frac{a+b+c}{3} - abc - \frac12 = a \left( \frac{1}{3} - bc\right) + \frac{b+c}{3} - \frac12.$$
In particular, inequality \eqref{first_ine} holds if and only if $f_{b,c}(a)<0$. Its derivative at any point $a\in[0,1]$ is given by \[f_{b,c}'(a) = \frac{1}{3} - bc.\] 
Now we investigate the maximal values of $f_{b,c}$, which can be considered in three different cases, depending on the sign of $f_{b,c}'$.
If $bc\geq 1/3$, then the function is decreasing in $[b,1]$, hence the maximum value is at $a=b$. In particular,  since $ab\leq 1/2$, we have $b^2\leq 1/2$ and hence, for any $a\geq b$,
$$f_{b,c}(a) \leq f_{b,c}(b) = \frac{2b}{3} + \frac{c}{3} - b^2c - \frac12 \leq \frac{2b}{3} + \frac{c}{3} - \frac{c}{3}-\frac12 \leq \frac{2}{3\sqrt{2}} - \frac12 < 0,$$
which gives inequality \eqref{first_ine}.
It remains to consider the case where $bc < 1/3$, for which the maximum is at $a=1$ and 
\[
f_{b,c}(1) = b\biggl( \frac{1}{3} - c \biggr) + \frac{c}{3} - \frac{1}{6}.
\]
Proceeding analogously, for any fixed $c\in[0,1]$, we define a function $g_c\colon[0,1]\to\bbR$, defined $$g_c(b) \coloneqq  b\biggl( \frac{1}{3} - c \biggr) + \frac{c}{3} - \frac{1}{6}.$$ 
Then, we have \[g_c'(b) = \frac{1}{3}-c\] and, considering that $b^2\leq ab\leq 1/2$, we have the following cases:
\begin{enumerate}[label=\textbullet, widest=*,nosep,parsep=4pt,wide=0pt]
    \item if $c <1/3$, $g_c$ is increasing and hence, for every $b\leq 1/2$, 
    $$g_c(b)\leq g_c\biggl(\frac12\biggr) = -\frac{c}{6}< 0;$$
    \item if $c \geq 1/3$, $g_c$ is decreasing and hence, for every $b\geq c$, 
    $$g_c(b)\leq g_c(c) = \frac{2c}{3}-c^2 - \frac{1}{6} < 0.$$
\end{enumerate}
This proves inequality \eqref{first_ine} also when $ab\leq 1/2$ and hence concludes the proof.
\end{proof}

\begin{remark}
Based on computational calculations in \texttt{MATLAB} \cite{MATLAB} involving the discretization of the cube $[0,1]^3$ into subintervals of the length 0.001 along each axis, we get $\mu \approx 0.52733$.
\end{remark}

\begin{proposition}\label{prop:epslowercond}
    Let $k,m\in\bbZ_{>0}$ be such that $k\leq m-1$ and let $\eps\in\bbR_{>0}$. If $L_{k,m}^-\leq \eps$, then
    \[
r_{m,\eps}\geq \frac{3m^2}{k+1}+\biggl\lceil\frac{(m-k)^3}{3(m-k)-2}\biggr\rceil.
    \]
\end{proposition}
\begin{proof}
    It is sufficient to notice that the stronger inequality
    \[
\frac{1+\eps}{3}m^2\geq \frac{3m^2}{k+1}+\frac{(m-k)^3}{3(m-k)-2}+1
    \]
    is equivalent to 
    \[
\eps\geq 
\frac{9}{k+1}-\frac{(6k-2)m^{2}-(9k^{2}+9)m+(3k^{3}+9k+6)}{m^2(3m-3k-2)},
    \] which is satisfied since $L_{k,m}^-\leq \eps$.
\end{proof}
\begin{proposition}\label{prop:secondepslowercond}
    Let $k,m\in\bbZ_{>0}$ be such that $k\leq m-1$ and let $\eps\in\bbR_{>0}$. If $L_{k,m}^-\leq \eps$, then
    \[
r_{m,\eps}\geq \biggl\lceil\frac{m^3}{3m-2}\biggr\rceil.
    \]
\end{proposition}

\begin{proof}
    It is sufficient to notice that the stronger inequality
    \[
\frac{1+\eps}{3}m^2\geq \frac{m^3}{3m-2}+1
    \]
    is equivalent to 
    \[
\eps\geq 
\frac{2m^2+9m-6}{m^2(3m-2)},
    \] which is satisfied since $L_{k,m}^-\leq \eps$.
\end{proof}

\begin{proposition}\label{prop:epsCondition}
    Let $k,m\in\bbZ_{>0}$ be such that $k\leq m-1$ and let $\eps\in\bbR_{>0}$. If $\eps<L_{k,m}^+$, then
$$\frac{(3m-2)r_{m,\eps}}{2} < m^{2}\bigl(m\mu - \max \{3k,m\varepsilon\} - 3\bigr).$$  
\end{proposition}
\begin{proof}
    We have to consider the two cases where either $\varepsilon<3k/m$ or $\varepsilon\geq 3k/m$. These are, respectively, 
    $$\biggl(\frac{3m-2}{2}\biggr)\biggl\lceil\frac{(1+\varepsilon)m^2}{3}\biggr\rceil  < \bigl(\mu m - 3k - 3\bigr)m^{2},\qquad \biggl(\frac{3m-2}{2}\biggr)\biggl\lceil\frac{(1+\varepsilon)m^2}{3}\biggr\rceil  < \bigl(\mu m - \eps m - 3\bigr)m^{2}.$$ 
    In particular, the statement is obtained just by solving the stronger  
    inequalities 
    \[\biggl(\frac{3m-2}{2}\biggr)\biggl(\frac{(1+\varepsilon)m^2}{3}+1\biggr)  < \bigl(m\mu - 3k - 3\bigr)m^{2},\quad \biggl(\frac{3m-2}{2}\biggr)\biggl(\frac{(1+\varepsilon)m^2}{3}+1\biggr)  < \bigl(m\mu - m\eps - 3\bigr)m^{2},\]
    which are equivalent to
    \[
\frac{(1+\varepsilon)m^2}{3}< \frac{2\mu m^3-6km^2-6m^2-3m+2}{3m-2},\qquad \frac{(1+\varepsilon)m^2}{3}<\frac{2\mu m^3-2\eps m^3-6m^2-3m+2}{3m-2},
    \]
    respectively.
The first inequality is equivalent to
\[
\varepsilon< \frac{3(2\mu-1)m^3-18km^2-16m^2-9m+6}{3m^3-2m^2}=2\mu-1+\frac{2\bigl((2\mu-1)-9k-8\bigr)m^2-9m+6}{3m^3-2m^2},
\]
while the second one is given by
\[
 \eps<\frac{3(2\mu-1)m^3-16m^2-9m+6}{9m^3-2m^2}=\frac{2\mu-1}{3}+\frac{2\bigl((2\mu-1)-24\bigr)m^2-27m+18}{3(9m^3-2m^2)},
\]
which are both satisfied since $\eps < L_{k,m}^+$.
\end{proof}

\begin{proposition}\label{prop:2conditions}
There exist $\eps \in \mathbb{R}_{>0}$ and \(k, \bar{m}_{\eps}\in\bbZ_{>0}\), with $k\leq \bar{m}_{\eps}-1$, such that $L_{k,m}^{-}\leq\varepsilon<L_{k,m}^{+}$ for every integer $m>\bar{m}_{\eps}$.
\end{proposition}
\begin{proof}

By \Cref{prop: mu>1/2} we have $2\mu - 1 > 0$, hence there exists $k \in \mathbb{Z}_{>0}$ satisfying $$\frac{9}{k+1} < \frac{2\mu-1}{3}.$$ 
For any such $k$ we have $$L^-_{k,m} = \max\left\{h_1(m)+\frac{9}{k+1}, h_2(m)\right\}, \ \ \ \ \ L^+_{k,m} = \min\left\{(2\mu-1) + h_3(m), \frac{2\mu-1}{3} + h_4(m)\right\}$$ with $h_i \in o(1)$ for $i = 1,2,3,4$. Thus, for any $\eps\in\bbR_{>0}$ such that
\[
\frac{9}{k+1} < \eps < \frac{2\mu-1}{3},
\]
there exists a positive integer $\bar{m}_{\eps} \geq k+1$ such that for every $m > \bar{m}_{\eps}$ we have \[L^-_{k,m} = h_1(m) + \frac{9}{k+1} \leq \eps < \frac{2\mu-1}{3}+h_4(m) = L^+_{k,m}.
\qedhere
\]
\end{proof}

\bibliographystyle{amsalpha}

\addcontentsline{toc}{section}{References}
\bibliography{Strassen.bib}

\vspace*{-3ex}

\end{document}